\documentclass[12pt]{amsart}
\usepackage{a4}
\usepackage{amsfonts}
\usepackage{amssymb}
\usepackage[latin1]{inputenc}

\author{Fran\c cois Court\`es} 
\title[Distinction of the Steinberg representation III]{Distinction of the Steinberg representation III: the tamely ramified case}
\address{Universit\'e de Poitiers\\Laboratoire de Math\'ematiques et Applications\\UMR 7348 du CNRS\\T\'el\'eport 2\\Boulevard Marie et Pierre Curie\\86962 Futuroscope Chasseneuil Cedex}
\email{courtes@math.univ-poitiers.fr}
\begin{document}
\begin{abstract}
Let $F$ be a nonarchimedean local field, let $E$ be a Galois quadratic extension of $F$ and let $G$ be a quasisplit group defined over $F$; a conjecture by Dipendra Prasad states that the Steinberg representation $St_E$ of $G(E)$ is then $\chi$-distinguished for a given unique character $\chi$ of $G(F)$, and that $\chi$ occurs with multiplicity $1$ in the restriction of $St_E$ to $G(F)$. In the first two papers of the series, Broussous and the author have proved the Prasad conjecture when $G$ is $F$-split and $E/F$ is unramified; this paper deals with the tamely ramified case, still with $G$ $F$-split.
\end{abstract}
\maketitle
\vskip 1cm
\font\teuf=eufm10
\font\seuf=eufm7
\font\sseuf=eufm6
\newfam\euffam
\textfont\euffam=\teuf
\scriptfont\euffam=\seuf
\scriptscriptfont\euffam=\sseuf
\def \got{\fam\euffam}
\def \mth{\mathbb}
\newtheorem{theo}{Theorem}[section]
\newtheorem{prop}[theo]{Proposition}
\newtheorem{lemme}[theo]{Lemma}
\newtheorem{cor}[theo]{Corollary}
\newtheorem{defi}[theo]{Definition}
\newtheorem{conj}{Conjecture}

\section{Introduction}

Let $F$ be a nonarchimedean local field with finite residual field, let $E$ be a Galois quadratic extension of $F$ and let $G$ be a reductive group defined over $F$. Let $G_E$ (resp. $G_F$) be the group of $E$-points (resp. $F$-points) of $G$ and let $\pi$ be a smooth representation of $G_E$; we say that $\pi$ is {\em distinguished} with respect to the symmetric space $G_E/G_F$ if the space $Hom_{G_F}(\pi,1)$, where $1$ is the one-dimensional trivial representation of $G_F$, is nontrivial. This article deals with the important particular case of the distinction of the Steinberg representation of $G_E$.

In \cite{pr0}, Dipendra Prasad has proved that when $G=GL_2$, the Steinberg representation $St_E$ of $G_E$ is not distinguished with respect to $G_E/G_F$; on the other hand, if we set $\chi=\varepsilon_{E/F}\circ det$, where $\varepsilon_{E/F}$ is the norm character of $E^*/F^*$, the space $Hom_{G_F}(St_E,\chi)$ happens to be of dimension $1$. For that reason, the definition of distinguishedness will be extended  the following way: let $\chi$ be any character of $G_F$; we say that $\pi$ is {\em $\chi$-distinguished} with respect to $G_E/G_F$ if $Hom_{G_F}(\pi,\chi)$ is nontrivial.

In \cite{pr}, Prasad has stated a conjecture about the distinction of the Steinberg representation which generalizes his result of \cite{pr0}; the conjecture, as initially stated, concerns quasisplit groups, but can be extended to any connected reductive group (see \cite{pr2}). Let $G^{ad}$ be the adjoint group $G/Z$, where $Z$ is the center of $G$, and let $\chi_{ad}$ be some given character of the group $G^{ad}_F$ of $F$-points of $G^{ad}$, called the Prasad character (see \cite{pr} for the definition of the Prasad character in the case of a $F$-quasisplit group, and \cite{pr2} for its extension to the general case. Note that since this article only deals with $F$-split groups, we will simply use \cite{pr} as a reference for its definition.) We then have:

\begin{conj}[Prasad]
The representation $St_E$ is $\chi_{ad}$-distinguished with respect to $G^{ad}(F)$, and $Hom_{G^{ad}(F)}(St_E,\chi_{ad})$ is one-dimensional. Moreover, $St_E$ is not $\chi'$-distinguished for any character $\chi'$ of $G^{ad}(F)$ distinct from $\chi_{ad}$.
\end{conj}

It is not hard to see that the above conjecture is equivalent to the same one with $G_{ad}$ replaced with $G$ and $\chi_{ad}$ with the Prasad character $\chi$ of $G_F$. The result has been proved for $G=GL_n$ and $F$ of characteristic $0$ by Anandavardhanan and Rajan (\cite{ar}), and more recently by Matringe for $G$ being an inner form of $GL_n$ and $F$ of characteristic different from $2$ (\cite{mat}). It has also been proved for any $F$-split $G$ by Broussous and the author (\cite{bc} and \cite{cou2}) when $E/F$ is unramified;  the present article deals with the tamely ramified case. More precisely, we prove the following results, which are the respective analogues of \cite[theorems 1 and 2]{bc}: let $\chi$ be the Prasad character of $G_F$ relative to $E/F$; we have:

\begin{theo}\label{th1}
Assume $G$ is split over $F$ and $E/F$ is totally and tamely ramified. The Steinberg representation $St_E$ of $G_E$ is then $\chi$-distinguished with respect to $G_F$.
\end{theo}

\begin{theo}\label{th2}
With the same hypotheses, the character $\chi$ occurs with multiplicity at most $1$ in the restriction to $G_F$ of $St_E$, and $St_E$ is not $\chi'$-distinguished for any character $\chi'$ of $G_F$ distinct from $\chi$.
\end{theo}

By the previous remarks we do not lose any generality by assuming that $G$ is semisimple and adjoint. To make proofs clearer, we even assume that $G$ is simple, the general case of semisimple groups being easy to deduce from the simple case.

The proof uses the model of the Steinberg representation that was already used in \cite{bc}: the Steinberg representation can be viewed as the space of smooth harmonic cochains over the set of chambers of the Bruhat-Tits building of $G_E$, with $G_E$ acting on it via its natural action twisted by a charater $\varepsilon$ (defined in section $3$), whose restriction to $G_F$ happens to be trivial when $E/F$ is ramified (proposition \ref{epstriv}). To prove theorem \ref{th1}, we thus only need to exhibit a $(G_F,\chi)$-equivariant linear form on that space, as well as a test vector for that form.
This is done in subsections $7.2$ and $7.3$. We prove the convergence of our linear form and the existence of a test ector with the help of the Poincar\'e series of affine Weyl groups (see \cite[section 3]{mcdo}), which allows us to get rid of the condition on $q$ we had to impose in \cite{bc}: the trick should work in the unramified case as well, which would lead to a simpler proof than the one given in \cite{bc} and \cite{cou2}. The author thanks Paul Broussous and Dipendra Prasad for suggesting him to use these series.

 To prove theorem \ref{th2}, as in \cite[section 6]{bc}, we prove the equivalent result that the space of $G_{F,der}$-invariant harmonic cochains on the building $X_E$, where $G_{F,der}$ is the derived group of $G_F$, is of dimension at most $1$ (sections $5$ and $6$). We will proceed by induction on the set $Ch_E$ of chambers of $X_E$, as in \cite{bc}, but since it turns out that contrary to the unramified case, the support of our harmonic cochains is not the whole set $Ch_E$, the induction we use here is quite different from the one of \cite{bc}.

We start by partitioning the set of chambers of $X_E$ into $F$-anisotropy classes the following way: set $\Gamma=Gal(E/F)$. For every chamber $C$, there exists a $\Gamma$-stable apartment $A$ of $X_E$ containing $C$ and an $E$-split $F$-torus $T$ attached to $A$ (proposition \ref{apst}, see also \cite{ds}); $A$ and $T$ are not unique, but the $F$-anisotropy class of $T$ does not depend on the choice of $A$ (corollary \ref{anisdim}), and we define the $F$-anisotropy class of $C$ as that class. Our goal is to prove theorem \ref{th2} with the help of an induction on these classes.

Contrary to the unramified case, the building $X_F$ of $G_F$ is not a subcomplex of the building $X_E$ of $G_E$, but if we consider their respective geometric realizations $\mathcal{B}_F$ and $\mathcal{B}_E$, the former is still the set of $\Gamma$-stable points of the latter, at least when $E/F$ is tamely ramified, and we can thus consider the set $Ch_\emptyset$ of chambers of $X_E$ whose geometric realization is contained in $\mathcal{B}_F$; that set is obviously $G_F$-stable, but in the ramified case, it contains more than one $G_F$-orbit of chambers. We thus first have to prove that the restrictions of our $G_{F,der}$-invariant harmonic cochains to $Ch_\emptyset$ are entirely determined by their value on some given element of $Ch_\emptyset$.

It quickly turns out that we have to treat the case of groups of type $A_{2n}$ separately from the other cases. In the case of type $A_{2n}$, the $G_{F,der}$-invariant harmonic cochains are identically zero on $Ch_\emptyset$ outside a particular orbit of chambers that we call $Ch_c$ (corollary \ref{chsol1}). We then use an induction (similar in its basic idea to the one of \cite[section 6]{bc}, but technically quite different) to prove that these harmonic cochains are entirely determined by their constant value on $Ch_c$, which proves theorem \ref{th2} in this case (corollary \ref{th2a2n}). In the proof of theorem \ref{th1} in the case of a $q$ large enough, our linear form $\lambda$ has its support on $Ch_c$, and our test vector is the Iwahori-spherical vector $\phi_C$ relative to some given chmnber $C$ in $Ch_c$; we also compute explicitly the value of $\lambda(\phi_C)$ (proposition \ref{tva2n}).

In the case of groups of type other than $A_{2n}$, the $G_{F,der}$-invariant harmonic cochains are identically zero on the whole set $Ch_\emptyset$ (corollary \ref{chsol1} again). In fact, it turns out that we can prove with our induction that these cochains are identically zero on the whole set $Ch_E$ outside a unique $F$-anisotropy class $Ch_a$, on which the induction fails; that class corresponds to the $E$-split tori of $G$ whose $F$-anisotropic component is of maximal dimension (corollary \ref{suppaut}); we thus may use as a starting point for a new induction the subset $Ch_a^0$ of the elements of $Ch_a$ which contain a $\Gamma$-fixed facet of $X_E$ of the greatest possible dimension; we prove in a similar way as in \cite[section 6]{bc} that the $G_{F,der}$-invariant harmonic cochains are entirely determined by their values on $Ch_a^0$ (corollary \ref{gtransa2}), then we check that the space of the restrictions to $Ch_a^0$ of our $G_{F,der}$-invariant harmonic cochains is of dimension at most $1$ (section $6$). That part of the proof is rather technical because $Ch_a^0$ does not consist of one single $G_{F,der}$-orbit in general; it is also the reason why, to prove theorem \ref{th1}, the test vector we choose in section $7.3$ is not an Iwahori-spherical vector. (Note that at the end of the paper (corollary \ref{testv}), we prove that an Iwahori-spherical vector attached to some given element of $Ch_a^0$ works as well, but using it as a test vector in the first place leads to a more complicated proof.)

The model used in this article and the previous ones can probably be used as well for the remaining cases. For groups whose $E$-rank and $F$-rank are the same, the induction should work the same way. For groups whose $E$-rank and $F$-rank are different, the induction has to be modified to take into acccount the fact that the apartments of $\mathcal{B}_F$ are now proper affine subspaces of the apartments of $\mathcal{B}_E$, but the same basic principle still applies.

The author also expects it to be possible to use the same model and a pretty similar proof to prove the Prasad conjecture in the wildly ramified case as well, but in that case, additional technical problems arise. The main two are the following ones: firstly, it is not true anymore that every chamber of $X_E$ is contained in a $\Gamma$-stable apartment; that problem can be adressed by considering, for chambers which do not satisfy that condition, $\Gamma$-stable parts of apartments instead of whole apartments, but we still need to extend the result of proposition \ref{gtransa1} to these bad chambers. Secondly, in the tamely ramified case, the geometric realizations of the $\Gamma$-fixed subspaces of such apartments are always contained in $\mathcal{B}_F$; this is not true anymore when $E/F$ is wildly ramified, which makes dealing with the values of the harmonic cochains on $Ch_a^0$ even more complicated than it already is in the tamely ramified case.

This paper is organized as follows. In section $2$, we define the notations we use throughout the paper. In section $3$, we give the definition of the Prasad character $\chi$, and we check that the $\chi$-distinction of the Steinberg representation is equivalent to the $\chi$-distinction of the natural representation of $G_E$ on the space of the smooth harmonic cochains over its Bruhat-Tits building $X_E$. In section $4$, we separate the set of chambers of $X_E$ into $F$-anisotropy classes. In section $5$, we determine the support of the $G_{F,der}$-invariant harmonic cochains, and we prove theorem \ref{th2} in the case of a group of type $A_{2n}$; for other types, we reduce the problem to a similar assertion over $Ch_a^0$. In section $6$, we deal with $Ch_a^0$ and finish the proof of theorem \ref{th2} for groups of type different from $A_{2n}$. In section $7$, finally, we prove theorem \ref{th1}.

\section{Notations}
Let $F$ be a nonarchimedean local field with discrete valuation and finite residual field. Let $E$ be a ramified Galois quadratic extension of $F$; $E/F$ is totally ramified, and is tamely ramified if and only if the residual characteristic $p$ of $F$ is odd.

Set $\Gamma=\rm{Gal}(E/F)$,; we denote by $\gamma$ its nontrivial element. We denote by $N_{E/F}: x\mapsto x\gamma(x)$ the norm application $x\mapsto x\gamma(x)$ from $E$ to $F$.

Let $\mathcal{O}_F$ (resp. $\mathcal{O}_E$) be the ring of integers of $F$ (resp. $E$), and let $\mathfrak{p}_F$ (resp. $\mathfrak{p}_E$) be the maximal ideal of $\mathcal{O}_F$ (resp. $\mathcal{O}_E$). Let $k_F=\mathcal{O}_F/\mathfrak{p}_F$ (resp. $k_E=\mathcal{O}_E/\mathfrak{p}_E$) be the residual field of $F$ (resp. $E$); since $E/F$ is totally ramidied, $k_E$ and $k_F$ are canonically isomorphic. Let $q=q_E=q_F$ be their common cardinality.

Let $\varpi_E$  be a uniformizer of $E$, and set $\varpi_F=N_{E/F}(\varpi_E)$. Since $E/F$ is totally ramified, $\varpi_F$ is a uniformizer of $F$.

Let $v=v_F$ be the normalized valuation on $F$ extended to $E$; we have $v(F)={\mth{Z}}\cup\{+\infty\}$ and $v(E)=\frac 12{\mth{Z}}\cup\{+\infty\}$.

Let $G$ be a connected reductive group defined and split over $F$. We fix a $F$-split maximal torus $T_0$ of $G$ and a Borel subgroup $B_0$ of $G$ containing $T_0$; $B$ is then $F$-split too. Let $\Phi$ be the root system of $G$ relative to $T_0$; in the sequel we assume $\Phi$ is irreducible. Let $\Phi^+$ be the set of positive roots of $\Phi$ corresponding to $B_0$, let $\Delta$ be the set of simple roots of $\Phi^+$ and let $\alpha_0$ be the highest root of $\Phi^+$. We also denote by $\Phi^\vee$ the set of coroots of $G/T_0$, and by $W$ the Weyl group of $\Phi$.

A Levi subgroup $M$ of $G$ is {\em standard} (relatively to $T_0$ and $B_0$) if $T_0\subset M$ and $M$ is a Levi component of some parabolic subgroup of $G$ containing $B_0$. A root subsystem $\Phi'$ of $\Phi$ is a {\em Levi subsystem} if it is the root system of some Levi subgroup of $G$ containing $T_0$; $\Phi'$ is {\em standard} if that Levi subgroup is standard, or in other words if $\Phi'$ is generated by some subset of $\Delta$.

For every algebraic extension $F'$ of $F$ and every algebraic group $L$ defined over $F'$, we denote by $L_{F'}$ the group of $F'$-points of $L$.

For every algebraic extension $F'$ of $F$, let $X_{F'}$ be the Bruhat-Tits building of $G_{F'}$: $X_{F'}$ is a simplicial complex whose dimension is, since $G$ is $F$-split, the semisimple rank $d$ of $G$. We have a set inclusion $X_F\subset X_E$ compatible with the action of $G_F$, but contrary to the unramified case, that inclusion is not simplicial. (Note that there exist isomorphisms of simplicial complexes between $X_E$ and $X_F$, but these isomorphisms are neither canonical nor useful for our purpose.) For that reason, we work most of the time with the geometric realization $\mathcal{B}_F$ (resp. $\mathcal{B}_E$) of $X_F$ (resp. $X_E$).

We have an inclusion $\mathcal{B}_F\subset\mathcal{B}_E$, and for every $x\in X_F$, $x$ has the same geometric realization in both $\mathcal{B}_F$ and $\mathcal{B}_E$. Once again, the inclusion is not simplicial: a facet of $\mathcal{B}_F$ is usually the (disjoint) union of several facets of $\mathcal{B}_E$ of various types. Moreover, when $E/F$ is tamely ramified, $\mathcal{B}_F$ is precisely the set of $\Gamma$-stable points of $\mathcal{B}_E$; this is not true when $E/F$ is wildly ramified.

For every facet $D$ of $X_E$ (resp. $X_F$), we denote by $R(D)$ its geometric realization in $\mathcal{B}_E$ (resp. $\mathcal{B}_F$). Similarly, if $A$ is an apartment of $X_E$ (resp. $X_F$), we denote by $R(A)$ its geometric realization in $\mathcal{B}_E$ (resp. $\mathcal{B}_F$). Note that $D$ can be a facet of both $X_E$ and $X_F$ at the same time only if it is a vertex, and $A$ cannot be an apartment of both $X_E$ and $X_F$ at the same time, hence there is no ambiguity with the notation.

Since $G_E$ and $G_F$ have the same semisimple rank, every apartment $\mathcal{A}$ of $\mathcal{B}_F$ is also an apartment of $\mathcal{B}_E$. Note that the apartments $A_E$ of $X_E$ and $A_F$ of $X_F$ whose geometric realization is $\mathcal{A}$ are different; we though have the (nonsimplicial) set equality $A_F=A_E\cap X_F$. We denote by $\mathcal{A}_0$ the apartment of $\mathcal{B}_F$ (and also of $\mathcal{B}_E$) associated to $T_0$, and by $A_{0,E}$ (resp. $A_{0,F}$) the apartment of $X_E$ (resp. $X_F$) whose geometric realization is $\mathcal{A}_0$.

For every subset $S$ of $\mathcal{B}_E$, let $K_{S,E}$ (resp $K_{S,F}$) be the connected fixator of $S$ in $G_E$ (resp. $G_F$); this is an open compact subgroup of $G_E$ (resp. $G_F$). If $D$ is a facet of $X_E$ (resp. $X_F$), we also write $K_{D,E}$ (resp. $K_{D,F}$) for $K_{R(D),E}$ (resp. $K_{R(D),F}$). If now $X$ is any subset of $X_E$ (resp. $X_F$), we define $K_{X,E}$ (resp. $K_{X,F}$) as the intersection of the $K_{x,E}$ (resp. $K_{x,F}$), $x\in X$; it is easy to check that this definition is consistent with the previous one when $X$ is a facet. Finally, if $T$ is a maximal torus of $G$ defined over $E$ (resp. $F$), we denote by $K_{T,E}$ (resp. $K_{T,F}$) the maximal compact subgroup of $T_E$ (resp. $T_F$); it is easy to check that if $A_E$ (resp. $A_F$) is the apartment of $X_E$ (resp. $X_F$) associated to $T$, we have $K_{T,E}=K_{A_E,E}$ (resp. $K_{T,F}=K_{A_F,F}$).

We say that a vertex $x$ of $X_E$ (resp. $X_F$) is {\em $E$-special} (resp. {\em $F$-special}) if $x$ is a special vertex of $X_E$ (resp. $X_F$), or in other words, if the root system of the reductive quotient $K_{x,E}/K_{x,E}^0$ (resp. $K_{x,F}/K_{x,F}^0$) relative to some maximal torus, where $K_{x,E}^0$ (resp. $K_{x,F}^0$) is the pro-unipotent radical of $K_{x,E}$ (resp. $K_{x,F}$), is the full root system $\Phi$ of $G_E$ (resp. $G_F$). Special vertices always exist (see \cite[\S 3, cor. to proposition 11]{bou5} for example). We also say that a vertex of $\mathcal{B}_E$ (resp. $\mathcal{B}_F$) is $E$-special (resp. $F$-special) if it is the geometric realization of some $E$-special (resp. $F$-special) vertex of $E$ (resp. $F$).

It is easy to prove that every $F$-special vertex of $X_F$ is also $E$-special, but the converse is not true: $E$-special vertices of $X_F$ are not necessarily $F$-special, and some $E$-special vertices of $X_E$ do not even belong to $X_F$.

We fix once for all a $F$-special vertex $x_0$ of $A_{0,E}$. We can identify $\mathcal{A}_0$ with the ${\mth{R}}$-affine space $(X_*(T)/X_*(Z))\otimes{\mth{R}}$, where $Z$ is the center of $G$, by setting the origin at $x_0$; the elements of $\Phi$ are then identified, via the standard duality product $<.,.>$ between $X^*(T)$ and $X_*(T)$, with affine forms on $\mathcal{A}_0$, and the walls of $\mathcal{A}_0$ as an apartment of $\mathcal{B}_F$ (resp. $\mathcal{B}_E$) are the hyperplanes satisfying an equation of the form $\alpha(x)=c$, with $\alpha\in\Phi$ and $c\in{\mth{Z}}$ (resp. $c\in\frac 12{\mth{Z}}$). Moreover, every facet $D$ of $A_{0,F}$ (resp. $A_{0,E}$) is determined by a function $f_D$ from $\Phi$ to ${\mth{Z}}$ (resp. $\frac 12{\mth{Z}}$) the following way: for every $\alpha\in\Phi$, $f_D(\alpha)$ is the smallest element of ${\mth{Z}}$ (resp. $\frac 12{\mth{Z}}$) which is greater or equal to $\alpha(x)$ for every $x\in R(D)$. If $D$ is a facet of $A_{0,F}$ (resp. $A_{0,E}$), $f_D$ satisfies the following properties:
\begin{itemize}
\item $f_D$ is a concave function, or in other words:
\begin{itemize}
\item for every $\alpha\in\Phi$, $f(\alpha)+f(-\alpha)\geq 0$;
\item for every $\alpha,\beta\in\Phi$ such that $\alpha+\beta\in\Phi$, $f(\alpha+\beta)\leq f(\alpha)+f(\beta)$.
\end{itemize}
\item for every $\alpha\in\Phi$, $f(\alpha)+f(-\alpha)\leq 1$ (resp. $\frac 12$);
\item if $D$ is a $F$-special (resp. $E$-special) vertex, then for every $\alpha\in\Phi$, $f(\alpha)+f(-\alpha)=0$. If $D$ is a chamber of $X_F$ (resp. $X_E$), then for every $\alpha\in\Phi$, $f(\alpha)+f(-\alpha)=1$ (resp. $\frac 12$).
\end{itemize}

Note that if $D$ is a $E$-special vertex of $X_E$ belonging to $X_F$ but not $F$-special, the functions $f_D$ attached to $D$ as a facet of respectively $X_E$ and $X_F$ are different. For these particular vertices, we have to denote by respectively $f_{D,E}$ and $f_{D,F}$ these two functions. In all other cases, either $D$ is a facet of only one of the two buildings or the concave functions are identical, and there is then no ambiguity with the notation $f_D$.

We denote by $C_{0,F}$ the chamber of $X_F$ such that $K_{C_{0,F}}$ is the standard Iwahori subgroup of $G_F$ (relative to $T_0$, $\Phi^+$ and $x_0$), or in other words the chamber of $\mathcal{A}^0$ whose associated concave function $f_{C_{0,F}}$ is defined by $f(\alpha)=0$ (resp. $f(\alpha)=1$) for every positive (resp. negative) $\alpha$. We also set $\mathcal{C}_{0,F}=R(C_{0,F})$.

For every $\alpha\in\Phi$, let $U_\alpha$ be the root subgroup of $G$ attached to $\alpha$, and let $\phi_\alpha$ be the valuation on $U_{\alpha,E}$ defined the following way: for every $u\in U_{\alpha,E}$, $\phi_\alpha(u)$ is the largest element of $\frac 12{\mth{Z}}$ such that $u$ fixes the half-plane $\alpha(x)\leq \phi_\alpha$ of $\mathcal{A}_0$ pointwise. (By convention, we have $\phi_\alpha(1)=+\infty$.) Obviously, the valuation on $U_{\alpha,F}$ defined in a similar way is just the restriction of $\phi_\alpha$ to $U_{\alpha,F}$, hence there is no ambiguity in the notation. The quadruplet $(G,T_0,(U_\alpha)_{\alpha\in\Phi},(\phi_\alpha)_{\alpha\in\Phi})$ is a valued root datum in the sense of Bruhat-Tits (see \cite[I. 6.2]{bt}).

Now we give the definition of the harmonic cochains that we will be using throughout the whole paper. Let $Ch_E$ be the set of chambers of $X_E$, and let $\mathcal{H}(X_E)$ be the vector space of harmonic cochains on $Ch_E$, or in other words the space of applications from $Ch_E$ to ${\mth{C}}$ satisfying the following condition (called the {\em harmonicity condition}): for every facet $D$ of codimension $1$ of $X_E$, we have:
\[\sum_{C\in Ch_E,D\subset C}f(C)=0.\]
The group $G_E$ acts naturally on $\mathcal{H}(X_E)$ by $g.f:C\mapsto f(g^{-1}C)$. For every subgroup $L$ of $G_E$, we denote by $\mathcal{H}(X_E)^L$ the subpace of $L$-invariant elements of $\mathcal{H}(X_E)$. We also denote by $\mathcal{H}(X_E)^\infty$ the subspace of smooth elements of $\mathcal{H}(X_E)$, which is the union of the $\mathcal{H}(X_E)^K$, with $K$ running over the set of open compact subgroups of $G$.

\section{The characters $\chi$ and $\varepsilon$}

Let $\chi$ be the character of $G_F$ defined the following way: let $\rho$ be the half-sum of the elements of $\Phi^+$. By \cite[\S I, proposition 29]{bou}, for every element $\alpha^\vee\in\Phi^\vee$, $<\rho,\alpha^\vee>$ is an integer, hence $<2\rho,\alpha^\vee>$ is even; we deduce from this that for every quadratic character $\eta$ of $F^*$, the character $\eta\circ 2\rho$ of $(T_0)_F$ is trivial on the subgroup of $(T_0)_F$ generated by the images of the $\alpha^\vee$, which is the group $(T_0)_F\cap G_{F,der}$, where $G_{F,der}$ is the derived group of $G_F$; $\eta\circ 2\rho$ then extends in a unique way to a quadratic character of $(T_0)_FG_{F,der}=G_F$; it is easy to check that such a character does not depend on the choice of $T_0$, $B_0$ and $\Phi^+$.

Let $\varepsilon_{E/F}$ be the quadratic character of $F^*$ associated to the extension $E/F$: for every $x\in F^*$, $\varepsilon_{E/F}(x)=1$ if and only if $x$ is the norm of an element of $E^*$. Let $\chi$ be the character $\varepsilon_{E/F}\circ 2\rho$ extended to $G_F$.

\begin{prop}
The character $\chi$ of $G_F$ is the Prasad character of $G_F$ relative to the extension $E/F$.
\end{prop}

According to \cite[section 2]{cou2}, the Prasad character is of the form $\varepsilon_{E/F}\circ\chi_0$ for some $\chi_0\in X^*(G)$, and we deduce from \cite[lemma 3.1]{cou2} that $\chi_0$ is trivial if and only if $\rho\in X^*(T)$. On the other hand, since $\varepsilon_{E/F}$ is of finite order, $\varepsilon_{E/F}\circ\chi_0$ factors through a subgroup of finite index $G_0$ of $G_F$, and in particular the proposition holds when the quotient $G_F/G_0$ is cyclic. By \cite[plates I to IX, (VIII)]{bou}, that condition is satisfied as soon as $\Phi$ is not of type $D_d$ with $d$ even, 

Assume then $\Phi$ is of type $D_d$, with $d=2n$ being even. By \cite[plate IV, (VII)]{bou}, we have:
\[\rho=\left(\sum_{i=1}^{2n-2}(2ni-\frac{i(i-1)}2\alpha_i)\right)+\frac{n(2n-1)}2(\alpha_{2n-1}+\alpha_{2n}).\]
When $n$ is even, $\rho$ belongs to $X^*(T)$ and $\chi$ is then trivial, hence the proposition holds again. Assume now $n$ is odd. Then by \cite[section 5]{cou2} again, we have for every $g\in G_F$:
\[\chi(g)=\varepsilon_{E/F}\circ(\alpha_{2n-1}+\alpha_{2n})(g),\]
and using the above expression of $\rho$, we obtain, given that $\varepsilon_{E/F}$ is quadratic:

\[\varepsilon_{E/F}\circ 2\rho(g)=\varepsilon_{E/F}\circ\sum_{i=1}^{2n-2}(4ni-{i(i-1)})\alpha_i(g)+\varepsilon_{E/F}\circ n(2n-1)(\alpha_{2n-1}+\alpha_{2n})(g)\]
\[=\varepsilon_{E/F}\circ(\alpha_{2n-1}+\alpha_{2n})(g).\]
Hence $\chi$ and $\varepsilon_{E/F}\circ\chi_0$ are equal, as desired. $\Box$

Note that, since we are dealing with a ramified extension here, the subgroup $G_0$ of $G_F$ we are using in the above proof is not the same as in \cite{cou2}, but this is of no importance: once we are reduced to a finite group, that group, up to a canonical isomorphism, depends only on $\Phi$ and not on $E$ and $F$, and the proof works exactly the same way in the ramified and unramified cases.

Let now $\varepsilon$ be the character of $G_E$ defined the following way: let $g$ be an element of $G_E$ and let $C$ be a chamber of $X_E$. Since $X_E$ is labellable (see for example \cite[IV, proposition 1]{br}), there exists a canonical bijection $\lambda$ between the vertices of $C$ and the vertices of $gC$, and the application $x\mapsto g\lambda^{-1}(x)$ is then a permutation of the set of vertices of $gC$. We set $\varepsilon(g)$ to be the signature of that permutation; it is easy to check (see \cite[lemma 2.1 (i) and (ii)]{bc}) that $\varepsilon$ is actually a character of $G_E$ and that it does not depend on the choice of $C$.

Let $(\pi_E,\mathcal{H}(X_E)^\infty)$ be the representation of $G_E$ defined the following way: for every $g\in G_E$ and every $f\in\mathcal{H}(X_E)^\infty$, we have:
\[\pi_E(g)f:C\in Ch_E\longmapsto\varepsilon(g)f(g^{-1}C).\]
By \cite[proposition 3.2]{bc}, the representation $(\pi_E,\mathcal{H}(X_E)^\infty)$ of $G_E$ is equivalent to $St_E\otimes\varepsilon$. On the other hand, when $E/F$ is ramified, we have:

\begin{prop}\label{epstriv}
The character $\varepsilon$ is trivial on $G_F$.
\end{prop}

Let $K_{T_0,F}$ be the maximal compact subgroup of $(T_0)_F$, and let $X_{T_0,F}$ be the subgroup of $T_0$ whose elements are the $\xi(\varpi_F)$, with $\xi\in X_*(T_0)$. From the decomposition $F^*=\varpi^{\mth{Z}}\mathcal{O}_F^*$ of $F^*$, we deduce the following decomposition of $(T_0)_F$:
\[(T_0)_F=K_{T_0,F}X_{T_0,F}.\]
Since $G_F=G_{F,der}(T_0)_F$, we finally obtain the following decomposition:
\[G_F=G_{F,der}K_{T_0,F}X_{T_0,F}.\]
Now consider the restriction of the character $\varepsilon$ to $G_F$. Since $G_{F,der}$ is contained in $G_{E,der}$, $\varepsilon$ is trivial on $G_{F,der}$; since $K_{T_0,E}$ fixes every chamber of $(A_0)_E$ pointwise, $\varepsilon$ is also trivial on that group, and in particular on $K_{T_0,F}$; finally, $X_{T_0,F}$ is generated by the $\xi(\varpi_F)$, $\xi\in X_*(T_0)$; since $\varpi_F$ is the product of $\varpi_E^2$ with some element $x$ of $O_E^*$, for every $\xi\in X_*(T_0)$, we have $\xi(\varpi_F)=\xi(\varpi_E)^2\xi(x))$, and since $\xi(x)\in K_{T_0,E}$ and $\varepsilon$ is quadratic and trivial on $K_{T_0,E}$, we obtain $\varepsilon(\xi(\varpi_F))=1$. Therefore, $\varepsilon$ is trivial on $X_{T_0,F}$, hence on $G_F$ and the proposition is proved. $\Box$

\begin{cor}
The restriction to $G_F$ of the representation $\pi'_E$ given by the natural action of $G_E$ on $\mathcal{H}(X_E)$ is isomorphic to the restriction of $St_E$.
\end{cor}

\begin{cor}
For every character $\chi$ of $G_F$, $Hom_{G_F}(St_E,\chi)$ and $Hom_{G_F}(\pi'_E,\chi)$ are canonically isomorphic.
\end{cor}

This last corollary proves that when $E/F$ is ramified, the $\chi$-distinctions of $St_E$ and $\pi'_E$ with respect to $G_E/G_F$ are two equivalent problems. For that reason, in the sequel, we work with $\pi'_E$ instead of $St_E$.

\section{The anisotropy class of a chamber}

In this section, we classify the chambers of $X_E$ according to the $F$-anisotropy classes of $E$-split $F$-tori of $G$, at least when $E/F$ is tamely ramified.

First we have to prove that for every chamber $C$, there exists a $E$-split maximal $F$-torus of $G$ such that $C$ is contained in the apartment of $X_E$ associated to $T$; this is an immediate consequence of the following result, which is the tamely ramified equivalent of \cite[Lemma A.2]{bc}:

\begin{prop}\label{apst}
Assume $E/F$ is tamely ramified. Let $C$ be any chamber of $X_E$; there exists a $\Gamma$-stable apartment of $X_E$ containing both $C$ and $\gamma(C)$.
\end{prop}

This is simply a particular case of \cite[proposition 3.8]{ds}. $\Box$

Note that the result of \cite{ds} is also valid when $E/F$ is unramified, but only when the residual characteristic of $F$ is odd; this is the reason why we used a different proof for \cite[lemma A.2]{bc}, which works for any $F$.

Note also that the above proposition is not true when $E/F$ is wildly ramified. As a counterexample, consider a $\Gamma$-stable chamber $C$ of $X_E$ whose geometric realization is not contained in $\mathcal{B}_F$; such chambers actually exist when $E/F$ is wildly ramified. Let $A$ be a $\Gamma$-stable apartment of $X_E$ containing $C$; since $\Gamma$ fixes a chamber of $A$, it fixes $A$ pointwise, which implies that $A$ is associated to some $F$-split torus of $G$, and we must then have $R(A)\subset\mathcal{B}_F$; since $R(C)\subset R(A)$ is already not contained in $\mathcal{B}_F$ by hypothesis, we reach a contradiction.

We now classify $E$-split $F$-tori of $G$ according to the roots of $G$ intervening in their anisotropic component. Recall that two elements $\alpha$ and $\beta$ of $\Phi$ are said to be {\em strongly orthogonal} if they are orthogonal (or in other words, if $<\alpha,\beta^\vee>=0$) and $\alpha+\beta$ is not an element of $\Phi$. First we prove some lemmas.

\begin{lemme}\label{negsorth}
Assume $\alpha$ and $\beta$ are strongly orthogonal. Then $-\alpha$ and $\beta$ are also strongly orthogonal.
\end{lemme}

If $\alpha$ and $\beta$ are orthogonal, then $-\alpha$ and $\beta$ are orthogonal as well. Moreover, let $s_\alpha\in W$ be the reflection associated to $\alpha$; we have $s_\alpha(\alpha+\beta)=-\alpha+\beta$, and since $\alpha+\beta\not\in\Phi$, $-\alpha+\beta$ cannot belong to $\Phi$ either and the lemma is proved. $\Box$

\begin{lemme}\label{stoshort}
Let $\alpha,\beta$ be two elements of $\Phi$. If $\alpha$ and $\beta$ are orthogonal and at least one of them is long, then they are strongly orthogonal.
\end{lemme}
(By convention, if $\Phi$ is simply-laced, all of its elements are considered long.)

It is easy to check (it is nothing else than the good old Pythagorean theorem) that when $\alpha$ and $\beta$ are orthogonal, $\alpha+\beta$ is strictly longer than either of them. Hence since $\Phi$ is reduced, $\alpha+\beta$ can be a root only if $\alpha$ and $\beta$ are both short. The lemma follows. $\Box$

\begin{lemme}\label{sigma1}
The following assertions are equivalent:
\begin{itemize}
\item there exists $w\in W$ such that $w(\alpha)=-\alpha$ for every $\alpha\in\Phi$;
\item there exists a subset $\Sigma$ of $\Phi$ whose cardinality is the rank $d$ of $\Phi$ and such that two distinct elements of $\Sigma$ are always strongly orthogonal.
\end{itemize}
Moreover, when $\Sigma$ exists, it is unique up to conjugation by an element of $W$.
\end{lemme}

Assume $w\in W$ is such that $w(\alpha)=-\alpha$ for every $\alpha\in\Phi$. We prove the first implication by induction on the rank $d$ of $\Phi$; we prove in addition that, if $\Sigma$ satisfies the conditions of the second assertion, we have:
\[w=\prod_{\alpha\in\Sigma}s_\alpha,\]
where for every $\alpha$, $s_\alpha$ is the reflection associated to $\alpha$. Note that since the elements of $\Sigma$ are all orthogonal to each other, the $s_\alpha$ commute, hence the above product can be taken in any order.

The case $d=0$ is trivial: assume $d>0$. Let $\alpha_0$ be the highest root in $\Phi^+$; by \cite[proposition 25 (iii)]{bou}, $\alpha_0$ is always a long root. Consider the elementary reflection $s_{\alpha_0}\in W$ associated to $\alpha_0$; the set $\Phi_{\alpha_0}$ of roots $\beta$ of $\Phi$ such that $s_{\alpha_0} w(\beta)=-\beta$ is precisely the set of elements of $\Phi$ which are orthogonal to $\alpha_0$, hence strongly orthogonal to $\alpha_0$ by lemma \ref{stoshort}. Moreover, $\Phi_{\alpha_0}$ is a closed and symmetrical subset of $\Phi$, hence a root subsystem of $\Phi$, of rank strictly smaller than $d$, and for every $\beta\in\Phi^+$, we have $s_{\alpha_0} w(\beta)=-\beta+<\beta,\alpha_0^\vee>\alpha_0$, which is negative if and only if $\beta\in\Phi_{\alpha_0}$; we can thus apply the induction hypothesis to $\Phi_{\alpha_0}$ and $s_{\alpha_0}w$ to obtain a subet $\Sigma'$ of $\Phi_{\alpha_0}$ satisfying the conditions of the second assertion (relatively to $\Phi_{\alpha_0}$ and such that we have:
\[s_{\alpha_0} w=\prod_{\beta\in\Sigma'}s_\beta.\]
Finally, we set $\Sigma=\Sigma'\cup\{\alpha_0\}$.

Note that $\Phi_{\alpha_0}$ may be reducible; in such a case, we apply the induction hypothesis to each one of its irreducible components and take as $\Sigma'$ the union of the sets of roots we obtain that way, given that two elements of $\Phi_{\alpha_0}$ which belong to different irreducible components are always strongly orthogonal.

It only remains to check that $\Sigma$ contains $d$ elements. Since these elements must be linearly independent, $\Sigma$ cannot contain more than $d$ of them. Assume it contains less than $d$ elements; there exists then $\beta\in\Phi$ which is not a linear combination of elements of $\Sigma$. On the other hand, it is easy to check (for example by decomposing it into a sum of terms of the form $s_\alpha(\beta')-\beta'$ (which is a multiple of $\alpha$), with $\alpha\in\Sigma$ and $\beta'\in\Phi$) that $w(\beta)-\beta$ is a linear combination of elements of $\Sigma$; we then cannot have $w(\beta)=-\beta$, hence a contradiction.

Conversely, let $\Sigma$ be a subset of $\Phi$ satisfying the conditions of the second assertion; set:
\[w=\prod_{\alpha\in\Sigma}s_\alpha.\]
Since the elements of $\Sigma$ are all orthogonal to each other, we must have $w(\alpha)=-\alpha$ for every $\alpha\in\Sigma_T$. Moreover, since the cardinality of $\Sigma_T$ is $d$ and its elements are linearly independent, they generate $X^*(T)\otimes{\mth{Q}}$ as a ${\mth{Q}}$-vector space, and every element of $\Phi$ is then  a linear combination of them, which implies, since $w$ extends to a linear automorphism of $X^*(T)\otimes{\mth{Q}}$, that we have $w(\alpha)=-\alpha$ for every $\alpha\in\Phi$, as required. $\Box$

By \cite[plates I to IX, (XI)]{bou}, the conditions of the above proposition are satisfied for every $\Phi$ except the following ones:
\begin{itemize}
\item $\Phi$ if type $A_d$ with $d>1$;
\item $\Phi$ of type $D_d$ with $d$ odd;
\item $\Phi$ of type $E_6$.
\end{itemize}

Now we separate the $E$-split maximal tori of $G$ into $F$-anisotropy classes. The reductive subgroups $L_0$ and $L$ of $G$ that we introduce in proposition \ref{acltorus} and its proof will be of some use later (see section $6$).

Assume $E/F$ is tamely ramified. Let $\mathcal{A}$ be a $\Gamma$-stable apartment of $\mathcal{B}_E$. Since $E/F$ is tamely ramified, $\mathcal{A}^\Gamma$ is contained in $\mathcal{B}_F$; by \cite[I. proposition 2.8.1]{bt}, there exists an apartment $\mathcal{A}'$ of $\mathcal{B}_F$ containing $\mathcal{A}^\Gamma$, and by eventually conjugating $\mathcal{A}$ by a suitable element of $G_F$, we can assume $\mathcal{A}'=\mathcal{A}_0$. Let $T$ be the $E$-split maximal torus of $G_E$ associated to $\mathcal{A}$; the $F$-split component $T_s$ of $T$ is then contained in $T_0$.

\begin{prop}\label{acltorus}
Let $a$ be the dimension of the $F$-anisotropic component of $T$. With the above hypotheses, there exists a unique (up to conjugation) subset $\Sigma_T$ of $\Phi$, of cardinality $a$, such that:
\begin{itemize}
\item $T$ is $G_F$-conjugated to some maximal torus of $G$ contained in the reductive subgroup $L_0$ of $G$ generated by $T_0$ and the root subgroups $U_{\pm\alpha}$, $\alpha\in\Sigma_T$, and $F$-elliptic in $L_0$;
\item if $\alpha,\beta\in\Sigma_T$, then $\alpha$ and $\beta$ are strongly orthogonal.
\end{itemize}
Conversely, for every $\Sigma\subset\Phi$ satisfying the second condition, there exists an $E$-split maximal torus $T$ of $G$ defined over $F$ such that we can choose $\Sigma_T=\Sigma$.
\end{prop}

Let $\mathcal{A}^\Gamma$ be the affine subspace of $\Gamma$-fixed points of $\mathcal{A}$; since $T_0$ contains the split component of $T$, every facet of maximal dimension of $\mathcal{A}^\Gamma$ is contained in the closure of some chamber of $\mathcal{A}^0$. Let $\mathcal{D}$ be such a facet; by eventually conjugating $T$ by a suitable element of $G_F$, we can assume that $\mathcal{D}$ is contained in the closure of $R(C_{0,F})$. 

Moreover, $T$ is contained in the centralizer $Z_G(T_s)$ of $T_s$ in $G$, hence if $\Sigma_T$ exists, we can assume that the root subgroups $U_{\pm\alpha}$, $\alpha\in\Sigma_T$, are also all contained in $Z_G(T_s)$. Hence by replacing $G$ by $Z_G(T_s)/T_s$, we can assume that $T$ is $F$-anisotropic, which implies that $\mathcal{D}$ is a vertex of $\mathcal{B}_E$ contained in $\mathcal{B}_F$. (Note that $\mathcal{D}$ is not necessarily a vertex of $\mathcal{B}_F$.) The existence of a subset $\Sigma_T$ of $\Phi$ of cardinality $d$ satisfying the strong orthogonality condition is then a consequence of lemma \ref{sigma1}, but we still have to prove that such a $\Sigma_T$ satisfies the first condition as well.

Since $T$ is $E$-split, there exists $g\in G_E$ such that $gTg^{-1}=T_0$; the conjugation by $g^{-1}$ sends then $\Phi$ to the root system of $G$ relative to $T$. Since $\mathcal{A}^\Gamma$ consists of a single point, the action of the nontrivial element $\gamma$ of $\Gamma$ on $\mathcal{A}$ is the central symmetry relative to that point. This means in particular that for every $\alpha\in\Phi$, $\gamma(Ad(g^{-1})\alpha)=-Ad(g^{-1})\alpha$.

Let $L_0$ be the subgroup of $G$ generated by $T_0$ and the root subgroups $U_{\pm\alpha}$, $\alpha\in\Sigma_T$. Set $L=gL_0g^{-1}$; $L$ is then the subgroup of $G$ generated by $T$ and the root subgroups $gU_{\pm\alpha}g^{-1}$, $\alpha\in\Sigma_T$. This group is a closed $E$-split reductive subgroup of $G$ of type $(A_1)^d$; moreover, for every $\alpha\in\Sigma_T$, since $\gamma(Ad(g^{-1})\alpha)=-Ad(g^{-1})\alpha$, we have $\gamma(gU_\alpha g^{-1})=gU_{-\alpha}g^{-1}$; we deduce from this that $L$ is $\Gamma$-stable, hence defined over $F$. To prove the first assertion of the proposition, we only have to prove that $L$ and $L_0$ are $G_F$-conjugates.

We first prove the following lemma:

\begin{lemme}\label{lsplit}
The group $L$ is $F$-split.
\end{lemme}

Since the elements of $\Sigma_T$ are all strongly orthogonal to each other, $L$ is $F$-isogeneous to the direct product of $d$ semisimple and simply-connected groups of type $A_1$, namely the groups generated by the $U_{\pm Ad(g)\alpha}$ for every $\alpha\in\Sigma$; moreover, since for every $\alpha\in\Sigma$, $\gamma$ swaps $Ad(g)\alpha$ and $-Ad(g)\alpha$, every such component is $\Gamma$-stable. On the other hand, by \cite[17.1]{spr}, there are exactly two simply-connected groups of type $A_1$ defined over $F$: the split group $SL_2$, and its unique nonsplit $F$-form, whose group of $F$-points is isomorphic to the group of the norm $1$ elements of the unique quaternionic division algebra over $F$ (these groups are the only inner $F$-forms of $SL_2$ by \cite[proposition 17.1.3]{spr}, and by the remark made at the beginning of \cite[17.1.4]{spr}, $SL_n$ can have outer forms only if $n\geq 3$). Let $F'$ be the unique quadratic unramified extension of $F$; these groups are both $F'$-split, which proves that $L$ must be $F'$-split as well.

Let $T'$ be a maximal $F'$-split $F$-anisotropic torus of $G$ contained in $L$ and let $K_{T',F}$ be the maximal compact subgroup of $T'_F$. By \cite[theorem 3.4.1]{deb}, there exists a pair $(K,{\mth{T}}')$, with $K$ being a maximal parahoric subgroup of $G_F$ and ${\mth{T}}'$ being a maximal $k_F$-torus in the quotient ${\mth{G}}=K/K^0$, $k_F$-anisotropic modulo the center of ${\mth{G}}$, such that $K_{T',F}\subset K$ and ${\mth{T}}'$ is the image of $K_{T',F}$ in ${\mth{G}}$; moreover, the dimension of the $k_F$-anisotropic component of ${\mth{T}}'$ is the same as the dimension of the $F$-anisotropic component of $T'$, which implies that ${\mth{G}}$ is of semisimple rank $d$ and ${\mth{T}}'$ is $k_F$-anisotropic.

Consider now the image ${\mth{L}}$ of $L_F\cap K$ in ${\mth{G}}$; ${\mth{L}}$ is the group of $k_F$-points of a reductive $k_F$-group $k_F$-isogeneous to the direct product of $d$ $k_{F'}$-split simply-connected groups of type $A_1$. Since by \cite[1.17]{car}, every group over a finite field is quasisplit, and since the only quasisplit simply-connected group of type $A_1$ over any field is $SL_2$, which is split, ${\mth{L}}$ is isogeneous to a $k_F$-split group, hence  is $k_F$-split itself and contains a $k_F$-split maximal torus ${\mth{T}}''$. Let $I$ be an Iwahori subgroup of $G_F$ contained in $K$ whose image in ${\mth{G}}$ contains ${\mth{T}}''$; considering the Iwahori decomposition of $I$ (or alternatively, using \cite[theorem 3.4.1]{deb} again), we see that there exists a maximal torus $T''$ of $G$ whose maximal compact subgroup $K_{T''}$ is contained in $I$ and such that ${\mth{T}}''$ is the image of $K_{T''}$ in ${\mth{G}}$, and $T''$ must then be $F$-split. Hence $L$ is $F$-split, as desired. $\Box$

Now we go back to the proof of proposition \ref{acltorus}. We prove $L$ is $G_F$-conjugated to $L_0$, and also the unicity of $\Sigma_T$ up to conjugation. By eventually conjugating $L$ by some element of $G_F$, we can assume that it contains $T_0$; $L$ is then generated by $T_0$ and the $U_{\pm\alpha}$, with $\alpha$ belonging to some set $\Sigma'$ satisfying the strong orthogonality condition, and $L$ and $L_0$ are $G_E$-conjugated by some element $n$ of the normalizer of $T_0$ in $G_E$, which implies that $\Sigma_T$ and $\Sigma'$ are $W$-conjugates. Moreover, since $G$ is $F$-split, it is possible to choose $n$ as an element of $G_F$, hence $L$ and $L_0$ are $G_F$-conjugates as well and the first assertion of proposition \ref{acltorus} is proved.

Now we prove the second assertion. Let $\Sigma$ be any subset of $\Phi^+$ such that every $\alpha\neq\beta\in\Sigma$ are strongly orthogonal. The reductive subgroup $L$ of $G$ generated by $T_0$ and the $U_{\pm\alpha}$, $\alpha\in\Sigma$, is then of type $A_1^a$, where $a$ is the cardinality of $\Sigma$; by quotienting $L$ by its center and considering separately every one of its irreducible components, we are reduced to the case where $L$ is a simple group of type $A_1$, hence isogeneous to $SL_2$; according to a well-known result about $SL_2$, since $E/F$ is quadratic and separable, $L$ contains a $1$-dimensional $E$-split $F$-anisotropic torus, as required (for example, when $E/F$ is tamely ramified, the group of elements of $SL_2$ of the form $\left(\begin{array}{cc}a&b\\-\varpi b&a\end{array}\right)$, where $\varpi$ is an uniformizer of $F$ which is the square of some uniformizer of $E$). $\Box$

More generally, since every $E$-split $F$-torus $T$ of $G$ is $G_F$-conjugated to some torus $T'$ whose $F$-split component is contained in $T_0$, by the previous proposition, we can attach to $T$ a subset $\Sigma_T$ of $\Phi^+$, defined up to conjugation, which is the subset attached to $T'$ by that proposition. The class of $\Sigma_T$ is called the {\em $F$-anisotropy class} of $T$.

Note that, although the $F$-anisotropy classes are parametred by the conjugacy classes of subsets of strongly orthogonal elements of $\Phi$, in the sequel, by a slight abuse of notation, we will often designate an $F$-anisotropy class by one of the representatives of the corresponding conjugacy class; more precisely, we will say "the $F$-anisotropy class $\Sigma$" instead of "the $F$-anisotropy class corresponding to the conjugacy class of subsets of strongly orthogonal elements of $\Phi$ which contain $\Sigma$".

Note also that, as we will see later, two $E$-split $F$-tori belonging to the same $F$-anisotropy class are not necessarily $G_F$-conjugates; we though have the following result:

\begin{prop}\label{fpconj}
Assume $E/F$ is tamely ramified. Let $T,T'$ be two $E$-split maximal $F$-tori of $G$ belonging to the same $F$-anisotropy class $\Sigma$ and let $\mathcal{A}$ (resp. $\mathcal{A}'$) be the $\Gamma$-stable apartment of $\mathcal{B}_E$ associated to $T$ (resp. $T'$). Then the affine subspaces $\mathcal{A}^\Gamma$ and $\mathcal{A}'^\Gamma$ are $G_{F,der}$-conjugates.
\end{prop}

Since $E/F$ is tamely ramified, $\mathcal{A}^\Gamma$ and $\mathcal{A}'^\Gamma$ are contained in $\mathcal{B}_F$, and by eventually conjugating $T$ and $T'$ by elements of $G_{F,der}$, we can assume that they are both contained in $\mathcal{A}_0$; they are then conjugated by some element $n$ of the normalizer of $T_0$ in $G_{E,der}$. Moreover, since $T_0$ is $F$-split, every element of the Weyl group of $G/T_0$ admits representatives in $G_F$, and even in $G_{F,der}$ since the Weyl groups of $G_F$ and $G_{F,der}$ are the same; hence by eventually conjugating $T$ again, we may assume $n\in T_0\cap G_{E,der}$. Finally, we have $T_0\cap G_{E,der}=(K_{T_0,E}\cap G_{E,der})(X_{T_0,E}\cap G_{E,der})$, where $X_{T_0,E}$ is the subgroup of $T_0$ generated by the $\xi(\varpi_E)$, $\xi\in X_*(T_0)$, and $K_{T_0,E}$ fixes $\mathcal{A}^\Gamma$ pointwise; we thus may assume that $n$ belongs to $X_{T_0,E}\cap G_{E,der}$, which is, since $G_{der}$ is simply-connected, the subgroup of $X_{T_0,E}$ generated by the $\alpha^\vee(\varpi_E)$, $\alpha^\vee\in\Phi^\vee$. In such a case, the split components of $T$ and $T'$ are both contained in $T_0$ and conjugated by an element of $T_0$, hence identical; we thus can assume that $\Sigma$ is contained in the root subsystem of the elements of $\Phi$ whose restriction to that common split component is trivial. 

We now prove the result with $n$ being of the form $\alpha^\vee(\varpi_E)$  for some $\alpha^\vee$; the general case follows by an easy induction. If  $<\beta,\alpha^\vee>=0$ for every $\beta\in\Sigma$, then $\mathcal{A}^\Gamma=\mathcal{A}'^\Gamma$ and there is nothing to prove. If $<\beta,\alpha^\vee>$ is odd for some $\beta\in\Sigma$, then either $\mathcal{A}^\Gamma$ or $\mathcal{A}'^\Gamma$, say for example $\mathcal{A}^\Gamma$, is contained in some hyperplane of $\mathcal{A}_0$ which is a wall in $\mathcal{B}_F$ and whose associated roots are $\pm\beta$; on the other hand, if $L$ is the reductive subgroup of $G$ associated to $T$ as in proposition \ref{acltorus} and if $L_\beta$ is the subgroup of $L$ generated by the root subgroups $U_{\pm\beta}$, $T\cap L_\beta$ is then split, hence $T$ is of anisotropy class strictly contained in $\Sigma$, which leads to a contradiction. Hence $<\beta,\alpha^\vee>$ must be even for every $\beta\in\Sigma$.

Assume now $<\beta,\alpha^\vee>$ is even for every $\beta\in\Sigma$ and nonzero for at least one $\beta$; that nonzero $<\beta,\alpha^\vee>$ must then be equal to $\pm 2$. As a consequence, there exists a wall $\mathcal{H}$ of the apartment $\mathcal{A}_0$ of $\mathcal{B}_F$ which separates\ $\mathcal{A}^\Gamma$ from $\mathcal{A}'^\Gamma$ and contains neither of them; if we assume the converse, we reach the same contradiction as above. Let $s_{\mathcal{H}}$ be the orthogonal symmetry with respect to $\mathcal{H}$; we obviously have $s_{\mathcal{H}}(\mathcal{A})=\mathcal{A}'$. On the other hand, $\mathcal{H}$ being a wall in the building $\mathcal{B}_F$, the element of the affine Weyl group of $T_0$ corresponding to $s_{\mathcal{H}}$ admits representatives in $G_{F,der}$; the result follows. $\Box$

Note that the above proof does not work in the wildly ramified case because $\mathcal{A}^\Gamma$ and $\mathcal{A}'^\Gamma$ are then not contained in $\mathcal{B}_F$ in general. The author conjectures that proposition \ref{fpconj} still holds in that case, though.

Now we want to divide $Ch_E$ into $F$-anisotropy classes as well. Of course the $\Gamma$-stable apartment containing $C$, hence also the $E$-split maximal $F$-torus associated to it, are not unique, but we can still prove the following result:

\begin{prop}\label{anisconj}
Assume $E/F$ is tamely ramified. Let $C$ be any chamber of $X_E$ and let $A$ and $A'$ be two $\Gamma$-stable apartments of $X_E$ containing $C$. Then the pairs $(C,A)$ and $(C,A')$ are $G_{F,der}$-conjugates.
\end{prop}

When $R(C)$ is contained in $\mathcal{B}_F$, $R(A)$ and $R(A')$ must also be contained in $\mathcal{B}_F$, and \cite[I. proposition 2.3.8]{bt} implies that they are then always $G_{F,der}$-conjugates. Assume now $R(C)$ is not contained in $\mathcal{B}_F$ and let $g$ be an element of $G_{E,der}$ such that $gC=C$, $g\gamma(C)=\gamma(C)$ and $gA=A'$; such an element exists by \cite[I. proposition 2.3.8]{bt} again. Moreover, we also have $\gamma(g)C=C$, hence $g\in K_{C\cap\gamma(C),E}$, and $\gamma(g)\mathcal{A}=\mathcal{A}'$; if we set $h=\gamma^{-1}(g)g$, we then have $hC=C$ and $hA=A$, which implies, if $T$ is the $E$-split maximal torus of $G$ attached to $A$, that $h\in K_{T,E}$.

Since $C\cup\gamma(C)$ contains a chamber of $A$, $K_{C\cup\gamma(C),E}$ is contained in an Iwahori subgroup of $G_E$, and since it contains $K_{T,E}$, we have $K_{C\cup\gamma(C),E}=K^0K_{T,E}$, where $K^0$ is the pro-unipotent radical of $K_{C\cup\gamma(C),E}$. By multiplying $g$ by a suitable element of $K_{T,E}$ on the right, we may assume $g\in K^0$, which implies $h\in K^0\cap K_{T,E}$. On the other hand, by \cite[corollary 1]{cou1}, the cohomology group $H^1(\Gamma,K^0\cap K_{T,E})$ is trivial, which implies that since $h=\gamma^{-1}(g)g$ satisfies $h\gamma(h)=1$, and thus defines a $1$-cocycle of $\Gamma=\{1,\gamma\}$, it also defines a $1$-coboundary of that same group, hence must admit a decomposition of the form $h=\gamma(h')h'^{-1}$, with $h'$ being an element of $K_{T,E}\cap K^0$; hence $gh'=\gamma(gh')$, which implies that $gh'$ is an element of $G_{F,der}$ such that $gh'\mathcal{C}=\mathcal{C}$ and $gh'\mathcal{A}=\mathcal{A}'$, as desired. $\Box$

Note that the tame ramification hypothesis is needed for the above proof because it is used by \cite[corollary 1]{cou1}, but the author believes that in the wildly ramified case, a similar result should hold for chambers of $X_E$ contained in at least one $\Gamma$-stable apartment.

\begin{cor}\label{anisdim}
Assume $E/F$ is tamely ramified. Let $C$ and $A$ be defined as in proposition \ref{anisconj}, let $T$ be the maximal $E$-split $F$-torus of $G$ associated to $A$ and let $\Sigma_T$ be a subset of $\Phi$ attached to $T$ as in proposition \ref{acltorus}. Then up to conjugation, $\Sigma_T$ does not depend on the choice of $A$.
\end{cor}

This is an obvious consequence of proposition \ref{anisconj}. $\Box$

In other words, the $F$-anisotropy class of the torus $T$ associated to a $\Gamma$-stable apartment $A$ of $X_E$ containing $C$ does not depend on the choice of $A$. We can now state the following definition:

\begin{defi}
Let $C$ be a chamber of $X_E$. The $F$-anisotropy class of $C$ is the $F$-anisotropy class of the $E$-split maximal torus of $G$ associated to any $\Gamma$-stable apartment of $X_E$ containing $C$. 
\end{defi}

\section{The support of the $G_{F,der}$-invariant harmonic cochains}

In this section, we start the proof of theorem \ref{th2}. In the unramified case (see \cite[section 6]{bc}), in order to prove a similar result, we fix a chamber $C_0$ of $X_F\subset X_E$ and then, for every $C\in Ch_E$, we prove by induction on the combinatorial distance between $C$ and $X_F$ that for every $f\in\mathcal{H}(X_E)^{G_{F,der}}$, $f(C)$ depends only on $f(C_0)$. In the ramified case, a similar approach would be to start from a chamber of $X_E$ whose geometric realization is contained in $\mathcal{B}_F$; it turns out that although that kind of approach works in the case of a group of type $A_{2n}$, in the other cases, $f$ is identically zero on the set of such chambers and we have to find another starting point for our induction. For that reason, we start by determining the support of the elements of $\mathcal{H}(X_E)^{G_{F,der}}$.  In particular, when $\Phi$ is not of type $A_{2n}$, we prove that their support coincides with some given anisotropy class of $Ch_E$, namely the one given by proposition \ref{noncond}.

\subsection{The class $Ch_\emptyset$}

First we consider the trivial $F$-anisotropy class $Ch_\emptyset$ of $Ch_E$, or in other words the $F$-anisotropy class corresponding to $\Sigma=\emptyset$. When $E/F$ is tamely ramified, a chamber $C$ belongs to the trivial anisotropy class if and only if its geometric realization is contained in an apartment $\mathcal{A}$ of $\mathcal{B}_E$ whose associated torus is $F$-split, which is true if and only if $\mathcal{A}\subset\mathcal{B}_F$. When $E/F$ is wildly ramified, we also define $Ch_\emptyset$ as the set of chambers of $X_E$ satisfying that property.

Contrary to the unramified case, the action of $G_{F,der}$ on $Ch_\emptyset$ is not transitive, and we thus have to check that the space of the restrictions of elements of $\mathcal{H}(X_E)^{G_{F,der}}$ to $Ch_\emptyset$ is of dimension at most $1$. We start by the following lemma:

\begin{lemme}\label{fwall}
Let $f$ be an element of $\mathcal{H}(X_E)^{G_{F,der}}$, and let $C$ be a chamber of $X_E$ such that $R(C)$ is contained in $\mathcal{B}_F$ and that the geometric realization of at least one of its walls is contained in a codimension $1$ facet of $\mathcal{B}_F$. Then $f(C)=0$.
\end{lemme}

Let $C_F$ (resp. $D_F$) be a chamber (resp. a codimesion $1$ facet) of $X_F$ such that $R(C_F)$ contains $R(C)$ (resp. $R(D_F)$ contains some wall $R(D)$ of $R(C)$), and let $S$ be a set of representatives in $G_{F,der}$ of the quotient group $K_{D_F,F}/K_{C_F,F}$. Since $C$ (resp. $D$) and $C_F$ (resp. $D_F$) have the same closure in $\mathcal{B}_F$, we have $K_{C_F,F}=K_{C,F}$ (resp. $K_{D_F,F}=K_{D,F}$); moreover, since $E/F$ is totally ramified, $K_{D_F,F}/K_{C_F,F}=K_{D,F}/K_{C,F}$ is isomorphic to $K_{D,E}/K_{C,E}$, hence the chambers $gC$, $g\in S$, are precisely the chambers of $X_E$ containing $D$; by the harmonicity condition, we then have $\sum_{g\in S}f(gC)=0$. On the other hand, since $f$ is $G_{F,der}$-invariant, we have $f(gC)=f(C)$ for every $g\in S$, hence the result. $\Box$

Now we determine which chambers of $Ch_\emptyset$ do or do not satisfy the condition of the previous lemma.

\begin{prop}\label{chsol}
The following conditions are equivalent:
\begin{itemize}
\item There exists a chamber $C$ in $Ch_\emptyset$ such that none of the walls of $R(C)$ is contained in a codimension 1 facet of $\mathcal{B}_F$. Moreover, every chamber of $\mathcal{B}_F$ contains a unique chamber of $\mathcal{B}_E$ satisfying that property;
\item The root system $\Phi$ is of type $A_{2n}$, with $n$ being a positive integer.
\end{itemize}
\end{prop}

Let $C$ be any element of $Ch_\emptyset$, and set $\mathcal{C}=R(C)$. Assume $C$ satisfies the condition of the proposition; since for every $g\in G_F$, $gC$ satisfies it too, we can assume that $\mathcal{C}$ is contained in $\mathcal{C}_{0,F}$. Let $f_C$ be the concave function on $\Phi$ associated to $C$ and let $\Delta'_C$ be the extended set of simple roots of $\Phi$ associated to $C$, which is the set of elements of $\Phi$ corresponding to the $d+1$ half-apartments of $\mathcal{A}_0$ whose intersection is $\mathcal{C}$. Since the walls of $\mathcal{C}$ are not contained in any codimension 1 facet of $\mathcal{B}_F$, we must have $f_C(\alpha)\in{\mth{Z}}+\frac 12$ for every $\alpha\in\Delta'_C$. On the other hand, let $\Delta=\{\alpha'_1,\dots,\alpha'_d\}$ be a set of simple roots of $\Phi$ contained in $\Delta'_C$ and let $\alpha'_0=-\sum_{i=1}^d\lambda_i\alpha_i$ be the remaining element of $\Delta'_C$; we have, with an obvious induction:
\[f_C(\alpha'_0)+\sum_{i=1}^d\lambda_if_C(\alpha'_i)=f_C(\alpha'_0)+f_C(\sum_{i=1}^d\lambda_i\alpha'_i)=f_C(\alpha'_0)+f_C(-\alpha'_0)=\frac 12.\]
For every $i\in\{0,\dots,d\}$, we have $f_C(\alpha_i)\in{\mth{Z}}+\frac 12$, which implies:
\[\frac 12\in{\mth{Z}}+(1+\sum_{i=1}^d\lambda_i)\frac 12,\]
hence the integer $1+\sum_{i=1}^d\lambda_i$ must be odd. By \cite[\S 1, proposition 31]{bou}, this integer is the Coxeter number of $\Phi$, and by \cite[plates I to IX, (III)]{bou}, it is odd if and only if $\Phi$ is of type $A_{2n}$ for some $n$; the first implication of the proposition is then proved.

Now assume $G$ is of type $A_{2n}$ for some $n$. We prove that $\mathcal{C}_{0,F}$ contains exactly one chamber of $\mathcal{B}_E$ satisfying the required condition; since that property translates by the action of $G_F$, every chamber of $\mathcal{B}_F$ satisfies it as well.

Let $\Delta'$ be an extended set of simple roots of $\Phi$ and let $\mathcal{C}'$ be the geometric realization of the chamber $C'$ of $A_0$ defined by the concave function $f$ such that:
\begin{itemize}
\item $f(\alpha)=\frac 12$ for every element $\alpha$ of $\Delta'$ different from some given one $\alpha_0$, and $f(\alpha_0)=\frac 12-n$;
\item for every $\beta\in\Phi$, writing $\beta=\sum_{\alpha\in S}\alpha$ for a suitable proper subset $S$ of $\Delta'$ (since $\Phi$ is of type $A_d$, such a subset exists, and it is unique), we have $f(\beta)=\sum_{\alpha\in S}f(\alpha)$.
\end{itemize}
Since $f(\alpha)$ is not an integer for any $\alpha\in\Delta'$, none of the walls of $\mathcal{C}'$ are contained in codimension $1$ facets of $\mathcal{B}_F$. The chamber $\mathcal{C}'$ is generally not contained in $\mathcal{C}_{0,F}$, but is always conjugated by an element of $G_F$ to some chamber $\mathcal{C}$ contained in $\mathcal{C}_{0,F}$ which satisfies the same property.

Now we prove the unicity of $\mathcal{C}$. We use the notations of \cite[plate I]{bou} (see also \cite[\S 4.4]{bou}): $\Phi$ is a subset of a free abelian group $X_0$ of rank $2n+1$ generated by elements $\varepsilon_1,\dots,\varepsilon_{2n+1}$ (this is the group denoted by $L_0$ in \cite{bou}; we rename it here to avoid confusion with the group $L_0$ of proposition \ref{acltorus}), the elements of $\Phi$ are the ones of the form $\alpha_{ij}=\varepsilon_i-\varepsilon_j$ with $i\neq j\in\{1,\dots,2n+1\}$, the elements of $\Phi^+$ being the ones such that $i<j$, and $W$ acts on $X_0$ by permutation of the $\varepsilon_i$. (The group $X_0$ is isomorphic to the character group of a maximal torus of $GL_{2n+1}$, and $W$ is isomorphic to the symmetric group $S_{2n+1}$.)

Let $\mathcal{C}=R(C)$ be a chamber of $\mathcal{B}_E$ contained in $\mathcal{C}_{0,F}$ and satisfying the required condition, and let $f_C$ be the concave function associated to $\mathcal{C}$. Since $\mathcal{C}$ is contained in $\mathcal{C}_{0,F}$, for every $\alpha\in\Phi^+$, we have $f_C(\alpha)\leq 0$ and $f_C(-\alpha)\leq 1$. On the other hand, we have $f_C(\alpha)+f_C(-\alpha)=\frac 12$, which implies $f_C(\alpha)\in\{-\frac 12,0\}$ and $f_C(-\alpha)\in\{\frac 12,1\}$.

Let $\Delta'$ be the extended set of simple roots associated to $C$; since for every $\alpha\in\Delta'$, we have $f_C(\alpha)\in{\mth{Z}}+\frac 12$, we must have $f_C(\alpha)=-\frac 12$ if $\alpha>0$ and $f_C(\alpha)=\frac 12$ if $\alpha<0$. On the other hand, the sum of the $f_C(\alpha)$, $\alpha\in\Delta'$, is $\frac 12$; $\Delta'$ must then contain exactly $n$ positive roots and $n+1$ negative roots.

Now we examine more closely the elements of $\Delta'$. Since $W$ acts transitively on the set of all extended sets of simple roots of $G$, there exists an element $w$ of $W$ such that $\Delta'$ is the conjugate by $w$ of the standard extended set of simple roots $\{\alpha_{12},\alpha_{23},\dots,\alpha_{2n+1,1}\}$, or in other words there exists a permutation $\sigma$ of $\{1,\dots,2n+1\}$ such that $\Delta'=\{\alpha_{\sigma(1)\sigma(2)},\alpha_{\sigma(2)\sigma(3)},\dots,\alpha_{\sigma(2n+1)\sigma(1)}\}$.

Assume that for some $i$ (with cycling indices), $\alpha_{\sigma(i)\sigma(i+1)}$ and $\alpha_{\sigma(i+1)\sigma(i+2)}$ are both positive. Then $f_C(\alpha_{\sigma(i)\sigma(i+2)})=-\frac 12-\frac 12=-1$, which is impossible by the previous remarks. Hence there must always be at least one negative root between two positive ones in the extended Dynkin diagram attached to $\Delta'$, which is a cycle of length $2n+1$. Since $\Delta'$ contains $n+1$ negative roots and $n$ positive roots, positive and negative roots must alternate on the diagram, except for two consecutive negative roots at some point. We can always choose $\sigma$ in such a way that the consecutive negative roots are $\alpha_{\sigma(2n+1)\sigma(1)}$ and $\alpha_{\sigma(1)\sigma(2)}$; in that case, $\alpha_{\sigma(i)\sigma(i+1)}$ is positive if and only if $i$ is even. We then easily obtain, for every $i<j$:
\begin{itemize}
\item if $i$ and $j$ are either both even or both odd, $f_C(\alpha_{\sigma(i)\sigma(j)})=0$, hence $\alpha_{\sigma(i)\sigma(j)}$ is positive, which implies $\sigma(i)<\sigma(j)$;
\item if $i$ is even and $j$ is odd, $f_C(\alpha_{\sigma(i)\sigma(j)})=-\frac 12$, hence $\alpha_{\sigma(i)\sigma(j)}$ is positive, which implies $\sigma(i)<\sigma(j)$;
\item if $i$ is odd and $j$ is even, $f_C(\alpha_{\sigma(i)\sigma(j)})=\frac 12$, hence $\alpha_{\sigma(i)\sigma(j)}$ is negative, which implies $\sigma(i)>\sigma(j)$.
\end{itemize}
In other words, the restriction of $\sigma$ to the subset of even (resp. odd) elements of $\{1,\dots,2n+1\}$ is an increasing function, and for every $i,j$ such that $i$ is even and $j$ odd, $\sigma(i)<\sigma(j)$. This is only possible if, for every $i$, $\sigma(2i)=i$ and $\sigma(2i+1)=n+i+1$, and $\Delta'$ is uniquely determined by these conditions. Since $\Delta'$ and the $f_C(\alpha)$, $\alpha\in\Delta'$, determine $C$, the unicity of $C$ is proved. $\Box$

\begin{cor}\label{chsol1}
When $\Phi$ is not of type $A_{2n}$ for any $n$, for every $f\in\mathcal{H}(X_E)^{G_{F,der}}$ and for every chamber $C$ of $Ch_\emptyset$, $f(C)=0$.

When $\Phi$ is of type $A_{2n}$ for some $n$, there exists a unique $G_F$-orbit $Ch_c$ of chambers of $X_E$ contained in $Ch_\emptyset$ and such that the elements of $\mathcal{H}(X_E)^{G_{F,der}}$ are identically zero on $Ch_\emptyset-Ch_c$.
\end{cor}

This is an immediate consequence of lemma \ref{fwall} and proposition \ref{chsol}. In the case $A_{2n}$, the orbit $Ch_c$ is the one described in the proof of proposition \ref{chsol}. $\Box$

Let $C_F$ be a chamber of $X_F$, and let $C$ be the unique element of $Ch_c$ whose geometric realization is contained in $R(C_F)$. We will call $C$ the {\em central chamber} of $C_F$.

\begin{cor}
The space of the restrictions to $Ch_\emptyset$ of the elements of $\mathcal{H}(X_E)^{G_{F,der}}$ is of dimension at most $1$.
\end{cor}

This is an immediate consequence of the previous corollary. $\Box$

\subsection{The other anisotropy classes}

Now we deal with the remaining anisotropy classes. First we prove that for every $C\in Ch_E$ which does not belong to $Ch_\emptyset$, $f(C)$ is entirely determined by the values of $f$ on some finite set of chambers in a given $\Gamma$-stable apartment containing $C$. We start with the following result:

\begin{prop}\label{gtransa1}
Assume $E/F$ is tamely ramified. Let $C$ be a chamber of $X_E$ whose geometric realization is not contained in $\mathcal{B}_F$. Let $A$ be a $\Gamma$-stable apartment of $X_E$ containing $C$, let $D$ be a wall of $C$ and let $C'$ be the other chamber of $A$ admitting $D$ as a wall. Assume that $C'$ is not contained in the closure $cl(C\cup\gamma(C))$ and that $D$ and $\gamma(D)$ are not contained in the same hyperplane of $A$. Let $Ch_D$ be the set of chambers of $X_E$ admitting $D$ as a wall and distinct from $C$; then $G_{F,der}\cap K_{C\cup\gamma(C)}$ acts transitively on $Ch_D$.
\end{prop}

First we observe that since $C$ and $\gamma(C)$ are both contained in the same half-plane delimited by the wall of $A$ containing $D$ (resp $\gamma(D)$), thay are then both contained in the closure of $C'\cup\gamma(C')$. In particular, we have $K_{C'\cup\gamma(C'),E}\subset K_{C\cup\gamma{C},E}$.

Let $T$ be the $E$-split maximal $F$-torus of $G$ corresponding to $A$; since $R(C)$ is not contained in $\mathcal{B}_F$, $T$ is not $F$-split. Let $g$ be an element of $G_E$ such that $gTg^{-1}=T_0$; $\Gamma$ then acts on the root system of $G$ relative to $T$, which is $Ad(g)^{-1}\Phi$, and its action is nontrivial. For every $\alpha\in\Phi$, let $U_{Ad(g)^{-1}\alpha}$ be the root subgroup of $G$ (relative to $T$) corresponding to $Ad(g)^{-1}\alpha$.

Let $H$ be the hyperplane of $A$ containing $D$, and let $\alpha$ be the element of $\Phi$ such that the root $Ad(g)^{-1}\alpha$ corresponds to the half-space $\mathcal{S}$ of $A$ delimited by $H$ and containing $C$; the group $U_{Ad(g)^{-1}\alpha,C}=U_{Ad(g)^{-1}\alpha}\cap K_{C,E}$ then acts transitively on $Ch_D$; moreover, since $\gamma(H)\neq H$, $\mathcal{S}$ contains both $\gamma(C)$ and $\gamma(C')$, hence $U_{Ad(g)^{-1}\alpha,C}$ fixes every element of $\gamma(Ch_D)$; we deduce from this that $\gamma(U_{Ad(g)^{-1}\alpha,C})$ fixes every element of $Ch_D$. Let now $C''$ be any element of $Ch_D$ and let $u$ be an element of the group $U_{Ad(g)^{-1}\alpha,C}$  such that $uC'=C''$; $u$ (resp. $\gamma(u))$ then fixes both $\gamma(C')$ and $\gamma(C'')$ (resp. both $C'$ and $C''$) and we obtain:
\[\gamma(u)uC'=u\gamma(u)C'=C''\]
and:
\[\gamma(u)u\gamma(C')=u\gamma(u)\gamma(C')=\gamma(C''),\]
We deduce from the above equalities that $h=u^{-1}\gamma(u^{-1})u\gamma(u)$ fixes both $C'$ and $\gamma(C')$, hence belongs to $K_{C'\cup\gamma(C'),E}$. Moreover, since $C'$ is a chamber, $K_{C'\cup\gamma(C'),E}$ is contained in an Iwahori subgroup of $G_E$, hence is pro-solvable, and since $h$ is a product of unipotent elements of $K_{C'\cap\gamma(C'),E}$, it then belongs to the pro-unipotent radical $K^0_{C'\cup\gamma(C'),E}$ of $K_{C'\cup\gamma(C'),E}\subset K_{C\cup\gamma(C),E}$.

Moreover, we have $h\gamma(h)=1$, hence $h$ defines once again a $1$-cocycle of $\Gamma=\{1,\gamma\}$ in $K^0_{C'\cup\gamma(C'),E}$. On the other hand, since $E/F$ is tamely ramified, by \cite[corollary 1]{cou1}, the cohomology set $H^1(\Gamma,K^0_{C'\cup\gamma(C'),E})$ is trivial, hence there exists $h'\in K^0_{C'\cup\gamma(C'),E}$ such that $h=h'^{-1}\gamma(h')$, which implies:
\[u^{-1}\gamma(u^{-1})u\gamma(u)\gamma(h')^{-1)}h'=1.\]
 Set $g'=u\gamma(u)\gamma(h')^{-1}$; we obtain $g'=\gamma(u)uh'^{-1}=\gamma(g')$, hence $g'\in G_{F,der}$, and $g'C'=C''$. Since this is true for every $C''$, $G_{F,der}\cap K_{C\cup\gamma(C),F}$ acts transitively on $Ch_D$, as required. $\Box$

\begin{cor}\label{gtransa2}
Assume $E/F$ is tamely ramified. Let $A$ be a $\Gamma$-stable apartment of $X_E$, let $Ch_A$ be the set of chambers of $X_E$ contained in $A$ and let $f$ be an element of $\mathcal{H}(X_E)^{G_{F,der}}$. The restriction of $f$ to $Ch_A$ is entirely determined by the values of $f$ on the chambers of $Ch_A$ containing a facet of maximal dimension of the set $A^\Gamma$ of $\Gamma$-stable elements of $A$. More precisely, if $C$ is any chamber of $Ch_A$ and $C'$ is a chamber of $Ch_A$ containing a facet of maximal dimension of $A^\Gamma$ and whose combinatorial distance to $C$ is the smallest possible, then $f(C)$ depends only on $f(C')$ and conversely.
\end{cor}

Let $C,C'$ be two elements of $Ch_A$; assume $C'$ contains a facet of maximal dimension of $A^\Gamma$. Let $(C_0=C',C_1,\dots,C_r=C)$ be a minimal gallery between $C'$ and $C$; assume also that $C'$ has been chosen in such a way that $r$ is the smallest possible. For every $i$, let $D_i=C_{i-1}\cap C_i$; if $D_i$ and $\gamma(D_i)$ are not contained in the same wall of $A$, by proposition \ref{gtransa1} (applied to the chambers containing $D_i$) and the harmonicity condition, we have either $f(C_i)+qf(C_{i-1})=0$ or $qf(C_i)+f(C_{i-1})=0$, hence $f(C_i)$ is determined by $f(C_{i-1})$ and conversely. Hence if for every $i$, $D_i$ satisfies that condition, by an obvious induction, we obtain that $f(C)$ is determined by $f(C')$ and conversely.

Assume now there exists some $i$ such that $D_i$ and $\gamma(D_i)$ are both contained in some wall $H$ of $A$; $H$ is then $\Gamma$-stable. Let $s_H$ be the reflection of $A$ relative to $H$, or in other words the only simplicial automorphism of $A$ fixing $H$ pointwise; since $H$ is $\Gamma$-stable, $\gamma\circ s_H\circ\gamma^{-1}$ is also such an automorphism, and must then be equal to $s_H$; in other words, the action of $\gamma$ on $A$ commutes with $s_H$, from which we deduce that $s_H(C')$ contains a facet of maximal dimension of $A^\Gamma$. On the other hand, we have $s_H(C_{i-1})=C_i$, hence $(s_H(C'),s_H(C_1),\dots,s_H(C_{i-2}),C_i,\dots,C_r=C)$ is a gallery (not necessarily minimal) between $s_H(C')$ and $C$ of length $r-1$; there must then exist a minimal gallery between them of length strictly smaller than $r$, which contradicts the minimality of $r$. Hence $D_i$ and $\gamma(D_i)$ are never contained in the same hyperplane of $A$ and the corollary is proved. $\Box$

Now we prove that when $G$ is not of type $A_{2n}$, the elements of $\mathcal{H}(X_E)^{G,{F,der}}$ are  identically zero on most of the $F$-anisotropy classes of $X_E$ (actually all but one, as we will see later with the help of proposition \ref{anismax}):

\begin{prop}\label{nonmax}
Assume $E/F$ is tamely ramified, and $G$ is not of type $A_{2n}$ for any $n$. Let $C$ be an element of $Ch_E$ such that $\Sigma_C$ is of cardinality $d-1$ and not maximal as a set of strongly orthogonal elements of $\Phi^+$. Then for every $f\in\mathcal{H}(X_E)^{G_{F,der}}$, $f(C)=0$.
\end{prop}

Note first that, by lemma \ref{sigma1} and the following remark, the condition on $\Sigma_C$ in fact already implies that $G$ is not of type $A_{2n}$. This is also true for the second assertion of proposition \ref{nonlstr}.

Let $f$ be any element of $\mathcal{H}(X_E)^{G_{F,der}}$, let $A$ be a $\Gamma$-stable apartment of $X_E$ containing $C$ and let $T$ be the $E$-split $F$-torus  of $G$ associated to $A$; by eventually conjugating $C$ by some element of $G_F$ we can assume that the split component $T_s$ of $T$ is contained in $T_0$, and even that $T$ is contained in the $F$-split reductive subgroup $L_0$ of $G$ defined as in proposition \ref{acltorus}. Moreover, since $\Sigma_C$ is of cardinality $d-1$ and not maximal, there exists a unique $\alpha\in\Phi^+$ which is strongly orthogonal to every element of $\Sigma_C$. The root subgroups $U_\alpha$ and $U_{-\alpha}$ are then normalized by $T_0$ and by every $U_{\pm\beta}$, $\beta\in\Sigma_C$, hence by $L_0$.

Let $h\in L_0$ be such that $hT_0h^{-1}=T$; since $\alpha$ is orthogonal to every element of $\Sigma_C$, the root $Ad(h)\alpha$ of $T$ does not depend on the choice of $h$. Let $H_\alpha$ be a wall of $A$ corresponding to $Ad(h)\alpha$ and containing some facet of $C$, and let $H'_\alpha$ be the wall of $A$ corresponding to the same root $\alpha$, neighboring $H_\alpha$ and such that $C$ is contained in the slice between them. Let $D$ be a facet of maximal dimension, hence of dimension $1$, of $A^\Gamma\subset A_0$, whose combinatorial distance to $C$ is the smallest possible; $D$ is then the unique edge of $A^\Gamma$ whose vertices are contained respectively in $H_\alpha$ and $H'_\alpha$. By corollary \ref{gtransa2}, $f(C)$ depends only on $f(C')$ for some chamber $C'$ of $A$ containing $D$, and conversely. 

Let $f_D$ be the concave function on $\Phi$ associated to $D$; since $\alpha$ is not a linear combination of the elements of $\Sigma_C$, we must have $f_D(\alpha)+f_D(-\alpha)=\frac 12$, hence either $f_D(\alpha)$ or $f_D(-\alpha)$, say for example $f_D(\alpha)$, is an integer. Let $D'$ be a facet of maximal dimension of $H_\alpha$ and let $C''$ be the unique chamber of $A$ containing $D'$ and whose remaining vertex is on the same side of $H_\alpha$ as $H'_\alpha$; we have $K_{D',F}/K_{C'',F}\subset K_{D',E}/K_{C'',E}$. If we prove that these two quotients are equal, then we obtain that $K_{D',F}$ acts transitively on the set of chambers containing $D'$; if in addition we prove that every class of $K_{D,F}/K_{C,F}$ contains elements of $G_{F,der}$, we then obtain by $G_{F,der}$-invariance and the harmonicity condition that the value of $f$ on every such chamber is zero, and in particular that $f(C'')=0$.

We thus prove that $K_{D',F}/K_{D',F}^0=K_{D',E}/K_{D',E}^0$, from which the first part of our claim follows immediately. Since $L_0$ normalizes the root subgroup $U_\alpha$ of $G$ associated to $\alpha$ and $K_D\cap L$ is a compact subgroup of $L$, we must have $hU_{\alpha,f_D(\alpha)}h^{-1}=U_{\alpha,f_D(\alpha)}$, and since $f_D(\alpha)$ is an integer, the quotient $U_{\alpha,f_D(\alpha)}/U_{\alpha,f_D(\alpha)+\frac 12}$ admits a system of representatives contained in $G_{F,der}$. Hence $U_{\alpha,f_D(\alpha)}$ is included in $K_{D',E}$, and $U_{\alpha,f_D(\alpha)+\frac 12}=U_{\alpha,f_D(\alpha)}\cap K_{D',E}^0$. On the other hand, by the same reasoning, we have $U_{-\alpha,-f_D(\alpha)}\subset K_{D',E}$ and $U_{-\alpha,-f_D(\alpha)+\frac 12}=U_{-\alpha,-f_D(\alpha)}\cap K_{D',E}^0$; hence the root subgroups of $K_{D',E}/K_{D',E}^0$ associated to both $\alpha$ and $-\alpha$ are contained in $K_{D',F}/K_{D',F}^0$, which is enough to prove that these two groups are equal. Moreover, at least $q$ classes of $K_{D',F}/K_{C'',F}$ out of $q+1$ contain elements of $U_\alpha\subset G_{F,der}$, hence the quotient $K_{D',F}\cap G_{F,der}/K_{C'',F}\cap G_{F,der}$, whose cardinality divides $q+1$, must be isomorphic to $K_{D',F}/K_{C'',F}$ and the second part of the claim is proved.

Now if we choose $D'$ in such a way that $C'$ is at minimal combinatorial distance from $C''$ among the chambers containing a facet of dimension $1$ of $\mathcal{A}^\Gamma$, by corollary \ref{gtransa2}, we then have $f(C')=0$, and then, also by the same corollary, $f(C)=0$, which proves the proposition. $\Box$

More generally, we have:

\begin{prop}\label{nonlstr}
\begin{itemize}
\item Assume $E/F$ is tamely ramified, and $G$ is not of type $A_{2n}$. Let $C$ be an element of $Ch_E$ which does not belong to $Ch_\emptyset$ and let $\Sigma_C$ be a subset of strongly orthogonal roots of $\Phi$ corresponding to the $F$-anisotropy class of $C$. Let $\Sigma_C^\perp$ be the set of elements of $\Phi$ which are strongly orthogonal to every element of $\Sigma_C$. Then $\Sigma_C^\perp$ is a closed root subsystem of $\Phi$.

\item Assume in addition that $\Sigma_C$ and $\Sigma_C^\perp$ are both nonempty and that $\Sigma_C^\perp$ is of rank $d-\#(\Sigma_C)$. Then for every $f\in\mathcal{H}(X_E)^{G_{F,der}}$, $f(C)=0$.
\end{itemize}
\end{prop}

To prove that $\Sigma_C^\perp$ is a closed root subsystem of $\Phi$, we only need to prove that:
\begin{itemize}
\item for every $\alpha,\alpha'\in\Sigma_C^\perp$ such that $\alpha+\alpha'\in\Phi$, $\alpha+\alpha'\in\Sigma_C^\perp$;
\item for every $\alpha\in\Sigma_C^\perp$, $-\alpha\in\Sigma_C^\perp$.
\end{itemize}

For every $\alpha\in\Sigma_C^\perp$, consider the reflection $s_\alpha$ associated to $\alpha$. Since $\alpha$ is orthogonal to every element of $\Sigma_C$, $s_\alpha$ fixes $\Sigma_C$ pointwise, which implies that $\Sigma_C^\perp$ is stable by $s_\alpha$, and in particular that it contains $s_\alpha(\alpha)=-\alpha$. Now let $\alpha,\alpha'$ be two elements of $\Sigma_C^\perp$ such that $\alpha+\alpha'$ is a root; since both of them are orthogonal to every element of $\Sigma_C$, then so is $\alpha+\alpha'$. Assume there exists $\beta\in\Sigma_C$ such that $\alpha+\alpha'+\beta$ is a root. Then $\beta$ is orthogonal to both $\alpha$ and $\alpha'$, which implies that $\alpha$, $\alpha'$ and $\beta$ are linearly independent; on the other hand, we deduce from lemma \ref{stoshort} that $\Phi$ is not simply-laced and $\alpha+\alpha'$ and $\beta$ are both short, which also implies, since $\alpha+\alpha'$ and $\beta$ are orthogonal, that $\alpha+\alpha'+\beta$ is long; the roots $\alpha$, $\alpha'$ and $\beta$ then generate a subsystem $\Phi'$ of $\Phi$ which is irreducible, not simply-laced and of rank $3$, hence of type either $B_3$ or $C_3$. Moreover, since $\alpha+\alpha'$ is short, either $\alpha$ or $\alpha'$, say $\alpha$, must be short.

In both cases below, the characters $\varepsilon_i$, $1\leq i\leq d$, are respectively defined as in plates II and III of \cite{bou}.

\begin{itemize}
\item Assume $\Phi'$ is of type $B_3$. In a system of type $B_d$, the sum of two nonproportional short roots $\pm\varepsilon_i$ and $\pm\varepsilon_j$ is always a long root $\pm\varepsilon_i\pm\varepsilon_j$. Hence $\alpha+\beta$ is a root, which contradicts the fact that $\alpha\in\Sigma_C^\perp$.
\item Assume $\Phi'$ is of type $C_3$. In a system of type $C_d$, two strongly orthogonal short roots are of the form $\pm\varepsilon_i\pm\varepsilon_j$ and $\pm\varepsilon_k\pm\varepsilon_l$, with $i,j,k,l$ being all distinct, which is obviously possible only if $d\geq 4$; hence $\alpha$ and $\beta$ cannot be strongly orthogonal, which once again leads to a contradiction.
\end{itemize}
Hence such a $\beta$ does not exist and $\alpha+\alpha'\in\Sigma_C^\perp$, which proves the first assertion of proposition \ref{nonlstr}.

Now we prove the second one. Assume first $\Sigma_C^\perp$ is irreducible. Let $A$, $D$ and $f_D$ be defined as in the proof of proposition \ref{nonmax} and let $D_1,\dots,D_{r+1}$ be the facets of $D$ of dimension $r-1$, with $r=d-\#(\Sigma_C)$ being the dimension of $D$. Let $H_1,\dots,H_{r+1}$ be the hyperplanes of $A$ respectively associated to the roots $\pm\alpha_1,\dots,\pm\alpha_{r+1}$ of $\Sigma_C^\perp$ which respectively contain $D_1,\dots,D_{r+1}$; the $H_i$ are then actually walls of $A$, and if for every $i$, $\alpha_i$ is the one among $\pm\alpha_i$ which is oriented towards $C$, the set $\{\alpha_1,\dots,\alpha_{r+1}\}$ is an extended set of simple roots of $\Sigma_C^\perp$. If $\lambda_1,\dots,\lambda_{r+1}$ are the smallest positive integers such that $\lambda_1\alpha_1+\dots+\lambda_{r+1}\alpha_{r+1}=0$, we must have $\lambda_1f_D(\alpha_1)+\dots+\lambda_{r+1}f_D(\alpha_{r+1})=\frac 12$; on the other hand, if $\Sigma_C^\perp$ is irreducible and not of type $A_{2n}$ for any $n$, by \cite[\S I, proposition 31]{bou}, the sum of the $\lambda_i$ is even, which implies that at least one of the $f_D(\alpha_i)$ must be an integer, and we finish the proof in a similar way as in proposition \ref{nonmax}. When $\Sigma_C^\perp$ is reducible and has no irreducible component of type $A_{2n}$ for any $n$, considering each irreducible component of $\Sigma_C^\perp$ separately, the proof is similar.

Now we check that $\Sigma_C^\perp$ cannot possibly have any irreducible component of type $A_{2n}$. Assume it admits such a component. Then the set $\Phi_C=\Sigma_C\cup -\Sigma_C\cup\Sigma_C^\perp$ is a proper closed root subsystem of $\Phi$ of rank $d$ admitting at least one component of type $A_1$ since $\Sigma_C$ is nonempty, and at least one component of type $A_{2n}$ for some $n$, which implies in particular that $d\geq 3$. Assume first that $\Phi_C$ is a parahoric subsystem of $\Phi$, or in other words the subsystem generated by $\Delta'-\{\alpha\}$, where $\Delta'$ is an extended set of simple roots of $\Phi$ and $\alpha$ is a nonspecial element of $\Delta'$; its Dynkin diagram is then the extended Dynkin diagram of $\Phi$ with the vertex corresponding to $\alpha$ removed. By examining the diagrams of the various possible parahoric subsystems of root systems of every type, we see that $\Phi_C$ can possibly have the required irreducible components only when $\Phi$ is of type $E_8$, $r=7$ and $\Sigma_C^\perp$ is of type $A_2\times A_5$, which implies that $\Phi_C$ is of type $A_1\times A_2\times A_5$. On the other hand, if $\Phi$ is of type $E_8$ and $\Sigma_C$ is a singleton, it is easy to check that $\Sigma_C^\perp$ must be of type $E_7$; we thus obtain a contradiction.

Now we look at the general case. By \cite[theorem 5.5]{dyn} and an obvious induction, we obtain a tower of root systems $\Phi=\Phi_0\supset\Phi_1\supset\dots\Phi_s=\Phi_C$ such that $\Phi_i$ is a parahoric subsystem of $\Phi_{i-1}$ for every $i$ and that $\Phi_s$ admits the required irreducible components. We deduce from the above discussion that if $\Phi_{s-1}$ is irreducible, it must be of type $E_8$, which, since $E_8$ is not contained in any other root system of rank $8$ ($A_8$ and $D_8$ are both strictly contained in $E_8$, as well as the systems of long roots of $B_8$ and $C_8$, which are respectively $D_8$ and $A_1^8$), implies $s=1$, we are then reduced to the previous case. Assume now $\Phi_{s-1}$ is reducible. Then in order for $\Phi_C$ to admit any component of type $A_{2n}$, there must exist an $i$ such that $\Phi_i$ admits such an irreducible component and $\Phi_{i-1}$ does not. The possible cases are, apart from the one which is already ruled out:
\begin{itemize}
\item $\Phi_{i-1}$ admits an irreducible component of type $E_6$ and that component breaks into three components of type $A_2$ in $\Phi_i$. According to the table on page $29$ of \cite{bt}, every vertex of the Dynkin diagram of a root system of type $A_n$ is special; we deduce from this that such a root system has no proper subsystems of the same rank. This implies that no component of type $A_1$ can arise in $\Phi_s$ from these three components, hence the components forming $\Sigma_C\cup-\Sigma_C$ must come from the other components of $\Phi_{i-1}$. But then the whole component of type $E_6$ of $\Phi_{i-1}$ is contained in $\Sigma_C^\perp\subset\Phi_s$, which contradicts the fact that it is already not contained in $\Phi_i$.
\item $\Phi_{i-1}$ admits a component of type $E_7$ which breaks into a system of type $A_2\times A_5$ in $\Phi_i$. For the same reason as above, $\Phi_{i-1}$ must admit components distinct from that component of type $E_7$ and containing the whole $\Sigma_C\cup-\Sigma_C$, and we reach the same contradiction.
\item $\Phi_{i-1}$ is of type $E_8$ and $\Phi_i$ is of type $E_6\times A_2$. Since the only possible way for $\Phi_s$ to have any component of type $A_1$ is that the component of type $E_6$ breaks into $A_1\times A_5$, we are reduced to a previous case.
\item $\Phi_{i-1}$ is of type $E_8$ and $\Phi_i$ is of type $A_4\times A_4$. There is no way that $\Phi_s$ can ever have any component of type $A_1$, since such a component should come from one of these two components of type $A_4$ and we already know that it is impossible.
\item $\Phi_{i-1}$ is of type $E_8$ and $\Phi_i$ is of type $A_8$. Same as above.
\item $\Phi_{i-1}$ is of type $F_4$ and $\Phi_i$ is of type $A_2\times A_2$. Same as above.
\item $\Phi_{i-1}$ is of type $G_2$ and $\Phi_i$ is of type $A_2$. This case is ruled out by the fact that we must have $d\geq 3$.
\end{itemize}
Since we always reach a contradiction, $\Sigma_C^\perp$ cannot have any irreducible component of type $A_{2n}$ and the proposition is proved. $\Box$

Note that in the course of the above proof, we have proved the following lemma which will be useful later:

\begin{lemme}\label{twoshort}
Let $\Sigma$ be a subset of strongly orthogonal elements of $\Phi$. Assume at least two elements of $\Sigma$ are short. Then $G$ is  of type $C_d$, with $d\geq 4$, and these two short elements of $\Sigma$ are of the form $\pm\varepsilon_i\pm\varepsilon_j,\pm\varepsilon_k\pm\varepsilon_l$, with $i,j,k,l$ being all distinct.
\end{lemme}

The following proposition allows us to eliminate more $F$-anisotropy classes from the support of the harmonic cochains:

\begin{prop}\label{noncond}
Assume $E/F$ is tamely ramified. Let $C,C'$ be two adjacent chambers of $X_E$, and let $D$ be the wall separating them. Let $A$ (resp. $A'$) be a $\Gamma$-stable apartment of $X_E$ containing $C$ (resp. $C'$) and let $T$ (resp. $T'$) be the corresponding $E$-split maximal $F$-torus of $G$. Let $\Sigma$ (resp. $\Sigma'$) be a subset of strongly orthogonal roots of $\Phi$ corresponding to the $F$-anisotropy class of $T$ (resp. $T'$); assume that:
\begin{itemize}
\item {\bf (C1)} there exist $\alpha\in\Sigma'$ and $\beta\in\Phi$ such that $\beta$ is orthogonal to every element of $\Sigma'$ except $\alpha$ and that $<\alpha,\beta^\vee>$ is odd;
\item $\Sigma'=\Sigma\cup\{\alpha\}$.
\end{itemize}
Let $Ch_D$ be the set of chambers of $X_E$ containing $D$ and distinct from both $C$ and the other chamber $C''$ containing $D$ and contained in $A$. Then $G_{F,der}\cap K_{C\cup\gamma(C)}$ acts transitively on $Ch_D$.
\end{prop}

By eventually conjugating $C$ and $C'$ by some element of $G_F$ we can assume that the $F$-split component of $T'$ is contained in $T_0$. Let $g$ (resp. $g'$) be an element of $G_E$ such that $gT_0g^{-1}=T$ (resp. $g'T_0g'^{-1}=T'$); define $\Sigma_T$ and $L_0$ as in proposition \ref{acltorus} and $\Sigma_{T'}$ and $L'_0$ in a similar way (relative to $T'$ instead of $T$), and set $L=gL_0g^{-1}$ and $L'=g'L'_0g'^{-1}$; since, by lemma \ref{lsplit}, $L$ and $L'$ are both $F$-split,  we obtain that $L$ is a $G_F$-conjugate of some subgroup of $L'$, and by multiplying $g'$ by a suitable element of the normalizer of $T_0$ in $G_F$, we actually obtain $L\subset L'$.. The roots corresponding to the hyperplane of $A'$ containing $D$ are then $\pm\alpha$; for every one-parameter subgroup $\xi$ of $T_0$ orthogonal to every element of $\Sigma$, $\xi(\mathcal{O}_F)$ then stabilizes $Ch_D$ globally. Moreover, by {\bf (C1)}, there exists a one-parameter subgroup $\xi$ in $X^\vee$ which is orthogonal to every element of $\Sigma$ and such that $<\alpha,\xi>$ is odd, and by adding to $\xi$ a suitable multiple of $\alpha^\vee$ we can assume that $<\alpha,\xi>=1$. Hence $\alpha\circ\xi$ is the identity on $F^*$, and in particular its restriction to $Ch_F$ is surjective, which implies that $\xi(\mathcal{O}_F^*)$, which is contained in $G_{F,der}\cup K_{C'\cup\gamma(C')}$, acts transitively on $Ch_D$. $\Box$

Now we consider the $F$-anisotropy classes which are not covered by the previous induction. Actually we prove that there is no such class when $G$ is of type $A_{2n}$, and exactly one when $G$ is of any other type:

\begin{prop}\label{anismax}
\begin{enumerate}
\item Assume $\Phi$ is not of type $A_{2n}$ for any $n$. There exists a subset $\Sigma_a$ of $\Phi$, unique up to conjugation by an element of the Weyl group of $\Phi$, satisfying the following properties:
\begin{itemize}
\item for every $\alpha,\beta\in\Sigma_a$, $\alpha$ and $\beta$ are strongly orthogonal, and $\Sigma_a$ is maximal for that property;
\item $\Sigma_a$ does not satisfy {\bf (C1)}.
\end{itemize}
\item With the same hypothese, $\Sigma_a$ is also maximal as a set of orthogonal roots of $\Phi$.
\item With the same hypothese once again,  every subset of strongly orthogonal elements of $\Phi$ which does not satisfy {\bf (C1)} is a conjugate of some subset of $\Sigma_a$.
\item Assume now $\Phi$ is of type $A_{2n}$ for some $n$. Then every nonempty subset of strongly orthogonal roots of $\Phi$ satisfies {\bf (C1)}. In particular, a subset $\Sigma_a$ of $\Sigma$ defined as above cannot exist.
\end{enumerate}
\end{prop}

First consider the case $A_{2n}$; we prove $(4)$ by induction on $n$. When $n=1$, no two roots of $\Phi$ are orthogonal to each other, which implies that every nonempty subset of orthogonal roots of $\Phi$ is a singleton; on the other hand, if $\alpha,\beta$ are any two nonproportional roots of $\Phi$, we have $<\alpha,\beta^\vee>=\pm 1$, hence $\{\alpha\}$ satisfies {\bf (C1)}. Assume now $n>1$, and let $\Sigma$ be any subset of strongly orthogonal elements of $\Phi$. Let $\alpha$ be any element of $\Sigma$; the subsystem $\Phi'$ of the elements of $\Phi$ which are orthogonal to $\alpha$ is then of type $A_{2n-2}$, and admits $\Sigma-\{\alpha\}$ as a subset of strongly orthogonal elements. If $\Sigma-\{\alpha\}$ is empty, then taking as $\beta$ any element of $\Phi$ which is neither proportional nor orthogonal to $\alpha$, we see that $\Sigma=\{\alpha\}$ satisfies {\bf (C1)}. Now assume $\Sigma-\{\alpha\}$ is nonempty. By the induction hypothesis, $\Sigma-\{\alpha\}$ must satisfy {\bf (C1)} as a subset of $\Phi'$. Let then $\alpha'\in\Sigma-\{\alpha\}$ and $\beta\in\Phi'$ be such that $\beta$ s orthogonal to every element of $\Sigma-\{\alpha,\alpha'\}$ and $<\alpha,\beta^\vee>$ is odd; by definition of $\Phi'$, $\beta$ is also orthogonal to $\alpha$. Hence $\Sigma$ satisfies {\bf (C1)} and $(4)$ is proved.

Assume now $\Phi$ is not of type $A_{2n}$; we first prove the existence of $\Sigma_a$. First consider the cases covered by lemma \ref{sigma1}, or in other words assume that there exists $w\in W$ such that $w(\alpha)=-\alpha$ for every $\alpha\in\Phi$; by lemma \ref{sigma1}, there exists then a subset $\Sigma_a$ of $d$ strongly orthogonal elements of $\Phi$; such a subset is necessarily maximal, and for every $\alpha\in\Sigma_a$, the only elements of $\Phi$ which are strongly orthogonal to every element of $\Sigma_a-\{\alpha\}$ are $\pm\alpha$, and since $<\alpha,\alpha^\vee>=2$ is even, $\Sigma_a$ does not satisfy {\bf (C1)}, as required.

Now we consider the remaining cases. Using the same algorithm as for lemma \ref{sigma1} (taking the highest root $\alpha_0$ of $\Phi^+$ and then considering the subsystem of the elements of $\Phi$ which are strongly orthogonal to $\alpha_0$), we also obtain a maximal subset $\Sigma_a$ of strongly orthogonal roots of $\Phi$, but this time, since $w=\prod_{\beta\in\Sigma_a}s_\beta$ cannot satisfy $w(\alpha)=-\alpha$ for every $\alpha\in\Phi$, by lemma \ref{sigma1}, $\Sigma_a$ contains strictly less than $d$ elements; we claim that for every $\alpha\in\Sigma_a$, the only elements of $\Phi$ which are strongly orthogonal to every element of $\Sigma_a-\{\alpha\}$ are $\pm\alpha$ once again.

Remember that the root systems we are considering here are the types $A_{2n-1}$ for some $n>1$ ($A_{2n}$ being already ruled out), $D_{2n+1}$ for some $n$ and $E_6$: since all these systems are simply-laced, by \cite[\S 1, 10, proposition 1]{bou}, two elements of $\Sigma_a$ are always conjugates of each other, which implies that we only have to prove the claim for one given $\alpha\in\Sigma_a$. By eventually conjugating $\Sigma_a$, we can always assume it contains $\alpha_0$. In the sequel, the simple roots $\alpha_1,\dots,\alpha_d$ of $\Phi^+$ are numbered as in \cite[plates I to IX]{bou}.
\begin{itemize}
\item Assume first $\Phi$ is of type $A_{2n-1}$, $n\geq 2$. The subsystem $\Phi'$ of the elements of $\Phi$ which are strongly orthogonal to $\alpha_0$ is then generated by the $\alpha_i$, $2\leq i\leq 2n-2$, hence of type $A_{2n-3}$. On the other hand, if $\alpha'$ is an element of $\Sigma_a$ distinct from $\alpha$, it is contained in $\Phi'$, and if the assertion is true for $\Phi'$, $\Sigma_a-\{\alpha_0\}$ and $\alpha'$, then it is also true fot $\Phi$, $\Sigma_a$ and $\alpha'$; we are then reduced to the case of type $A_{2n-3}$. By an obvious induction, after a finite number of such reductions we reach the case of a system of type $A_1$, and in that case, $\Sigma_a=\{\alpha_0\}$ obviously satisfies the required condition.
\item Assume now $\Phi$ is of type $D_{2n+1}$, $n\geq 2$. The subsystem of the elements of $\Phi$ which are strongly orthogonal to $\alpha_0$ is then generated by the $\alpha_i$, $i\neq 2$, hence of type $A_1\times D_{2n-1}$, the component of type $A_1$ being $\{\pm\alpha_1\}$. By eventually conjugating $\Sigma_a$ by the reflection $s_{\alpha_1}$, we may assume it contains $\alpha_1$ as well as $\alpha_0$, and by a similar reasoning as above (considering $\Sigma_a-\{\alpha_0,\alpha_1\}$ instead of $\Sigma_a-\{\alpha_0\}$), we are reduced to the case of type $D_{2n-1}$; after a finite number of such reductions we reach the case of a system of type $D_3=A_3$, which is an already known case.
\item Assume finally $\Phi$ is of type $E_6$. The subsystem of the elements of $\Phi$ which are strongly orthogonal to $\alpha_0$ is then generated by the $\alpha_i$, $i\neq 2$, hence of type $A_5$, and we are once again reduced to an already known case.
\end{itemize}
Now we prove the unicity of $\Sigma_a$ (up to conjugation) by induction on $d$, the case $d=1$ being obvious. Let $\Sigma$ be any subset of $\Phi$ satisfying the conditions of the proposition. Assume $\Sigma$ contains at least one long root (recall that by convention every root of a simply-laced system is considered long); by eventually conjugating $\Sigma$, we can assume that root is $\alpha_0$, and if $\Psi$ is the subsystem of the elements of $\Phi$ which are strongly orthogonal to $\alpha_0$, $\Sigma-\{\alpha_0\}$ satisfies the conditions of the proposition as a subset of $\Psi$, hence by induction hypothesis $\Sigma-\{\alpha_0\}$ and $\Sigma_a-\{\alpha_0\}$ are conjugated by an element $w$ of the Weyl group $W_\Psi$ of $\Psi$. Since $\alpha_0$ is orthogonal to every element of $\Psi$, it is fixed by $W_\Psi$, hence $\Sigma$ and $\Sigma_a$ are conjugated by $w$.

Assume now $\Phi$ is not simply-laced and $\Sigma$ contains only short roots. We now examine the different cases:
\begin{itemize}
\item Assume first $\Phi$ is of type $G_2$. Then no two roots of $\Phi$ are orthogonal, hence $\Sigma$ must be a singleton $\{\alpha\}$. Since there are long roots in $\Phi$ which are orthogonal to $\alpha$, hence strongly orthogonal by lemma \ref{stoshort}, $\Sigma$ cannot be maximal.
\item Assume now $\Phi$ is of type $C_d$. Let $\varepsilon_1,\dots,\varepsilon_d$ be defined as in \cite[plate III]{bou}. We have already seen (lemma \ref{twoshort}) that when $\Phi$ is of type $C_d$ and $\Sigma$ contains only short roots,  these roots must be of the form $\pm\varepsilon_i\pm\varepsilon_j$ with no two indices being identical; on the other hand, every possible index must show up in some $\pm\varepsilon_i\pm\varepsilon_j$, since if some index $k$ does not, the long root $2\varepsilon_k$ is strongly orthogonal to every element of $\Sigma$, which contradicts the maximality of $\Sigma$. Hence $d=2n$ is even and the only possible $\Sigma$ is, up to conjugation: $\Sigma=\{\varepsilon_1+\varepsilon_2,\dots,\varepsilon_{2n-1}+\varepsilon_{2n}\}$. On the other hand, the long root $\beta=2\varepsilon_1$ is orthogonal to every element of $\Sigma$ but $\alpha=\varepsilon_1+\varepsilon_2$, and we have $<\alpha,\beta^\vee>=1$, which contradicts the fact that $\Sigma$ must not satisfy {\bf (C1)}.
\item Assume now $\Phi$ is of type either $B_d$ or $F_4$. In both these cases, it is easy to check that no two orthogonal short roots are strongly orthogonal, hence $\Sigma$ must be a singleton. Let $\Phi'$ be a subsystem of type $B_2=C_2$ of $\Phi$ containing $\Sigma$; according to the previous case, $\Sigma$ satisfies {\bf (C1)} as a subset of $\Phi'$, hence also as a subset of $\Phi$ and we reach a contradiction once again.
\end{itemize}

In all the above cases, either $\Sigma$ is a conjugate of $\Sigma_a$ or we have reached a contradiction. Hence $(1)$ is proved.

Now we prove $(2)$. Assume there exists $\alpha\in\Phi$ which is orthogonal to every element of $\Sigma_a$. Then at least one element of $\Sigma_a$ is orthogonal but not strongly orthogonal to $\alpha$, which implies by lemma \ref{stoshort} that $\Phi$ is not simply-laced. On the other hand, $\Sigma_a$ is then of cardinality strictly smaller than $d$, which by lemma \ref{sigma1} and the following remark is possible only if $\Phi$ is of type $A_d$, with $d>1$ odd, $D_d$, with $d$ odd, or $E_6$, hence simply-laced.  We thus reach a contradiction, hence $\alpha$ cannot exist and $(2)$ is proved. 

Now we prove $(3)$. When $\Phi$ is simply-laced, we deduce from $(1)$ that every maximal subset of strongly orthogonal roots of $\Phi$ which does not satisfy $(C1)$ is conjugated to $\Sigma_a$, and $(3)$ follows immediately. When $\Phi$ is of type $G_2$, it is easy to check that every maximal subset of strongly orthogonal roots of $\Phi$ must always contain one long root and one short root, which also implies $(3)$. When $\Phi$ is of type $B_{2n+1}$ for some $n$, every subset of strongly orthogonal elements of $\Phi$ contains at most one short root (since in a system of type $B_d$, the sum of two nonproportional short roots is always a long root), and at most $2n$ long roots (since all of these long roots must be contained in the subsystem of the long roots of $\Phi$, which is of type $D_{2n+1}$ and, as we have already seen, does not contain any subset of strongly orthogonal elements of cardinality $2n+1$); using the induction of lemma \ref{sigma1} once again, we easily see that such a subset must also be contained in a conjugate of $\Sigma_a$; the assertion $(3)$ follows immediately in that case too.

It remains to consider the cases $B_{2n}$, $C_d$ and $F_4$. In all these cases, $\Sigma_a$ contains only long roots: this is easy to check by examining the subsystem $\Phi_l$ of the long roots of $\Phi$, which is of type respectively $D_{2n}$, $A_1^d$ and $D_4$; in all three cases, $\Phi_l$ contains a subset of $d$ strongly orthogonal roots which does not satisfy {\bf (C1)}, and such a subset must then be a conjugate of $\Sigma_a$ in $\Phi_l$, hence also in $\Phi$. If $\Sigma$ contains only long roots, by replacing $\Phi$ by $\Phi_l$, we are reduced to the simply-laced case. Assume now $\Sigma$ contains at least one short root $\alpha$; we prove by induction on the number of short roots it contains that it must satisfy {\bf (C1)}. By induction hypothesis, if $\Phi'$ is the subsystem of the elements of $\Phi$ which are strongly orthogonal to $\alpha$, $\Sigma-\{\alpha\}$ either satisfies {\bf (C1)} as a subset of $\Phi'$  or is contained in some conjugate of $\Sigma_a$ that by conjugating $\Sigma$ we may assume to be $\Sigma_a$ itself. In the first case, by the same argument as in the case $A_{2n}$, $\Sigma$ must satisfy {\bf (C1)} as a subset of $\Phi$. In the second case, since $\Sigma_a$ is of cardinality $d$, $\alpha$ is a linear combination of the elements of $\Sigma_a$, which is possible only if there exist $\beta_1,\beta_2\in\Sigma_a$ such that $\alpha=\frac 12(\pm\beta_1\pm\beta_2)$. We then have:
\[<\alpha,\beta_1^\vee>=\pm<\alpha,\beta_2^\vee> =\pm 1,\]
which proves at the same time that $\beta_1$ and $\beta_2$ do not belong to $\Sigma$ and that $\Sigma$ satisfies {\bf (C1)}. Hence $(3)$ is proved. $\Box$

\begin{cor}\label{clzero}
Assume $\Phi$ is not of type $A_{2n}$ for any $n$, $E/F$ is tamely ramified and $\Sigma'$ is a subset of strongly orthogonal roots of $\Phi$ which either is not maximal or satisfies {\bf (C1)}. Then for every $f\in\mathcal{H}(X_E)^{G_{F,der}}$ and every $C\in Ch_E$ of anisotropy class $\Sigma'$, $f(C)=0$.
\end{cor}

Assume first $\Sigma_a$ is of cardinality $d$. If $\Sigma'$ is a conjugate of some subset of $\Sigma_a$, then the set $\Sigma'^\perp$ of elements of $\Phi$ which are strongly orthogonal to $\Sigma'$ contains some conjugate of $\Sigma_a-\Sigma'$, hence is of dimension $d-\#(\Sigma')$  and we can just apply proposition \ref{nonlstr} if $\Sigma'$ is nonempty, and corollary \ref{chsol1} if $\Sigma'$ is empty. Assume now $\Sigma'$ is not a conjugate of any subset of $\Sigma_a$. By proposition \ref{anismax}(3), $\Sigma'$ satisfies {\bf (C1)} and we can proceed by induction. Let $\Sigma,C,C'$ be defined as in proposition \ref{noncond}; if we assume $f(C)=0$, by proposition \ref{noncond} we have $f(C')=0$ as well. By proposition \ref{anismax}(3), either $\Sigma$ is a conjugate of some subset of $\Sigma_a$, in which case we have $f(C')=0$ by the previous case, or $\Sigma$ satisfies {\bf (C1)}, in which case we can just iterate the process; since $\Sigma'$ is finite, after a finite number of steps we reach a situation where $\Sigma$ is conjugated to a subset (eventually empty) of $\Sigma_a$, hence by the previous case $f(C)=0$, and by proposition \ref{noncond} and an obvious induction, we must have $f(C')=0$. The fact that $f$ is then zero on the whole anisotropy class $\Sigma'$ of $Ch_E$ follows from corollary \ref{gtransa2}.

Assume now $\Sigma_a$ is of cardinality smaller than $d$. Then $\Phi$ is simply-laced, hence, as we have seen during the proof of proposition \ref{anismax}, $\Sigma'$ is always a conjugate of a subset of $\Sigma_a$. Now we examine the different cases:
\begin{itemize}
\item Assume $\Phi$ is of type $A_{2n-1}$. It is easy to check that for every $\alpha\in\Phi$, the subsystem of elements of $\Phi$ which are orthogonal to $\alpha$ is of type $A_{2n-3}$; we deduce from this that every proper subset of $\Sigma_a$, and more generally every nonmaximal subset of strongly orthogonal roots of $\Phi$, is contained in a subsystem of $\Phi$ of type $A_{2n-3}$, hence also in a subsystem of $\Phi$ of type $A_{2n-2}$; by proposition \ref{anismax}(4), $\Sigma'$ then satisfes {\bf (C1)}. We thus can apply proposition \ref{noncond} and an easy induction to get the desired result.
\item Assume $\Phi$ is of type $D_{2n+1}$, and, the $\varepsilon_i$ being defined as in \cite[plate IV]{bou}, set $\Sigma_a=\{\varepsilon_1\pm\varepsilon_2,\dots,\varepsilon_{2n-1}\pm\varepsilon_{2n}\}$. It is easy to check (details are left to the reader) that the set of $W$-conjugacy classes of sets of strongly orthogonal elements of $\Phi$ admits as a set of representatives the set of subsets $\{\Sigma_{i,j}|0\leq i\leq j\leq n\}$, with $\Sigma_{i,j}=\{\varepsilon_1\pm\varepsilon_2,\dots,\varepsilon_{2i-1}\pm\varepsilon_{2i},\varepsilon_{2i+1}+\varepsilon_{2i+2},\dots,\varepsilon_{2j-1}+\varepsilon_{2j}\}$; in particular, $\Sigma_{n,n}=\Sigma_a$. When $i<j$, setting for example $\alpha=\varepsilon_{2j-1}+\varepsilon_{2j}$ and $\beta=\varepsilon_{2j}+\varepsilon_{2j+1}$, we see that $\Sigma_{i,j}$ satisfies {\bf (C1)}. However, this is not true for the $\Sigma_{i,i}$, $0\leq i\leq n-1$, and we must then deal with them first. For every $i<n$, $\Sigma_{i,i}^\perp$ is a subsystem of type $D_{2(n-i)+1}$ of $\Phi$,  namely the set of roots which are linear combinations of the $\varepsilon_j$, $2i+1\leq j\leq 2n+1$; its rank is then equal to $d-\#(\Sigma_{i,i})$, and we can then apply proposition \ref{nonlstr} (or corollary \ref{chsol1} if $i=0$) to obtain that $f(C)=0$ in these cases. The cases $\Sigma_{i,j}$, $i<j$, then follow from the cases $\Sigma_{i,i}$ by proposition \ref{noncond} and an easy induction.
\item Assume $\Phi$ is of type $E_6$; $\Sigma_a$ is then contained in a Levi subsystem $\Phi'$ of type $D_4$ of $\Phi$, hence also in a Levi subsystem  $\Phi''$ of type $D_5$; we can thus define subsets $\Sigma_{i,j}$, $0\leq i\leq j\leq 2$, of that last subsystem in a similar way as in the previous proposition, and we can even assume they are contained in $\Phi'$. Mreover, if we assume that $\Phi'$ (resp. $\Phi''$) is generated by the elements $\alpha_2,\dots,\alpha_5$ (resp. $\alpha_1,\dots,\alpha_5$) of $\Delta$ (the simple roots being numbered as in \cite[plate V]{bou}), the elements of $W$ corresponding to the order $3$ automorphisms of the extended Dynkin diagram of $\Phi$ act on $\Phi'$ by automorphisms of order $3$ of its Dynkin diagram, and in particular permute the subsets $\{\alpha_2,\alpha_3\}$, $\{\alpha_2,\alpha_5\}$ and $\{\alpha_3,\alpha_5\}$ of $\Phi$; we deduce from this that the sets $\Sigma_{1,1}=\{\alpha_2,\alpha_5\}$ and $\Sigma_{0,2}=\{\alpha_3,\alpha_5\}$ belong to the same conjugacy class of sets of strongly orthogonal elements of $\Phi$. By proposition \ref{noncond} and the previous induction applied to $\Sigma_{0,0}\rightarrow\Sigma_{0,1}\rightarrow\Sigma_{0,2}$ and then to $\Sigma_{1,1}\mapsto\Sigma_{1,2}$, we obtain the desired result.
\end{itemize}
The corollary is then proved. $\Box$

\begin{prop}\label{sigmaa}
In the various cases, the sets $\Sigma_a$ are, up to conjugation, the following ones:
\begin{itemize}
\item when $G$ is of type $A_{2n-1}$, $\Sigma_a=\{-\varepsilon_1+\varepsilon_{2n},-\varepsilon_2+\varepsilon_{2n-1},\dots,-\varepsilon_n+\varepsilon_{n+1}\}$;
\item when $G$ is of type $B_{2n}$, $\Sigma_a=\{-\varepsilon_1\pm\varepsilon_2,-\varepsilon_3\pm\varepsilon_4,\dots,-\varepsilon_{2n-1}\pm\varepsilon_{2n}\}$;
\item when $G$ is of type $B_{2n+1}$, $\Sigma_a=\{-\varepsilon_1\pm\varepsilon_2,-\varepsilon_3\pm\varepsilon_4,\dots,-\varepsilon_{2n-1}\pm\varepsilon_{2n},-\varepsilon_{2n+1}\}$;
\item when $G$ is of type $C_d$, $\Sigma_a=\{-2\varepsilon_1,\dots,-2\varepsilon_d\}$;
\item when $G$ is of type $D_d$, with $d$ being either $2n$ or $2n+1$, $\Sigma_a=\{-\varepsilon_1\pm\varepsilon_2,-\varepsilon_3\pm\varepsilon_4,\dots,-\varepsilon_{2n-1}\pm\varepsilon_{2n}\}$;
\item when $G$ is of type $E_6$, $\Sigma_a=\{-\alpha_0,-\alpha_1-\alpha_3-\alpha_4-\alpha_5-\alpha_6,-\alpha_3-\alpha_4-\alpha_5,-\alpha_4\}$;
\item when $G$ is of type $E_7$, $\Sigma_a=\{-\alpha_0,-\alpha_2-\alpha_3-2\alpha_4-2\alpha_5-2\alpha_6-\alpha_7,-\alpha_2-\alpha_3-2\alpha_4-\alpha_5,-\alpha_2,-\alpha_3,-\alpha_5,-\alpha_7\}$;
\item when $G$ is of type $E_8$, $\Sigma_a=\{-\alpha_0,-2\alpha_1-2\alpha_2-3\alpha_3-4\alpha_4-3\alpha_5-2\alpha_6-\alpha_7,-\alpha_2-\alpha_3-2\alpha_4-2\alpha_5-2\alpha_6-\alpha_7,-\alpha_2-\alpha_3-2\alpha_4-\alpha_5,-\alpha_2,-\alpha_3,-\alpha_5,-\alpha_7\}$;
\item when $G$ is of type $F_4$, $\Sigma_a=\{-\alpha_0,-\alpha_2-2\alpha_3-2\alpha_4,-\alpha_2-2\alpha_3,-\alpha_2\}$;
\item when $G$ is of type $G_2$, $\Sigma_a=\{-\alpha_0,-\alpha_1\}$.
\end{itemize}
\end{prop}

The above sets $\Sigma_a$ are simply the ones we obtain by applying the algorithm of lemma \ref{sigma1}. Details are left to the reader. $\Box$

Note that for convenience (to be able to make the best possible use of lemma \ref{sroot}), we may want in the sequel to use representatives of $\Sigma_a$ which contain as many negatives of simple roots as possible, and we thus obtain: 

\begin{prop}\label{sigmaa2}
In the following cases, these alternative $\Sigma_a$ are also valid choices:
\begin{itemize}
\item when $G$ is of type $A_{2n-1}$, $\Sigma_a=\{-\alpha_1,-\alpha_3,\dots,-\alpha_{2n-1}\}$;
\item when $G$ is of type $D_{2n+1}$,  $\Sigma_a=\{-\varepsilon_2\pm\varepsilon_3,-\varepsilon_4\pm\varepsilon_5,\dots;-\varepsilon_{2n}\pm\varepsilon_{2n+1}\}$;
\item when $G$ is of type $E_6$, $\Sigma_a=\{\alpha_2+\alpha_3+2\alpha_4+\alpha_5,-\alpha_2,-\alpha_3,-\alpha_5\}$.
\end{itemize}
\end{prop}

Checking that these sets are also valid representatives of $\Sigma_a$ in their respective cases is straightforward, details are left to the reader. In the other cases, the representative of $\Sigma_a$ we pick up is still the one given by proposition \ref{sigmaa}.

In particular, we have proved the following result:

\begin{prop}\label{odd}
It is possible to choose $\Sigma_a$ in such a way that it is contained in a standard Levi subsystem $\Phi'$ of rank $\#(\Sigma_a)$ of $\Phi$ and that every one of its elements is the negative of the sum of an odd number of simple roots of $\Phi^+$ (counted with multiplicities).
\end{prop}

Checking that the condition of proposioin \ref{odd} is satisfied by the sets $\Sigma_a$ given by proposition \ref{sigmaa2} in the cases covered by that proposition and by proposition \ref{sigmaa} in the other cases is straightforward. $\Box$

The last three results of this section are three more corollaries to proposition \ref{anismax}.

 Let $Ch_a$ be the subset of chambers of $X_E$ of anisotropy class $\Sigma_a$, and let $Ch_a^0$ be the subset of the elements of $Ch_a$ containing a $\Gamma$-fixed facet of maximal dimension of any $\Gamma$-stable apartment containing them.

\begin{cor}\label{suppaut}
Assume $\Phi$ is not of type $A_{2n}$ for any $n$ and $E/F$ is tamely ramified. Let $f$ be an element of $\mathcal{H}(X_E)^{G_{F,der}}$; the support of $f$ is contained in $Ch_a$, and $f$ is entirely determined by its values on $Ch_a^0$.
\end{cor}

By proposition \ref{anismax}(3), every subset of strongly orthogonal roots of $\Phi$ which is not a conjugate of $\Sigma_a$ either is not maximal or satisfies {\bf (C1)}; the corollary then follows from corollaries \ref{gtransa2} and \ref{clzero}. $\Box$

In the case of groups of type $A_{2n}$, our induction actually works on the whole set $Ch_E$ and we obtain:

\begin{cor}\label{th2a2n}
Assume $\Phi$ is of type $A_{2n}$ for some $n$ and $E/F$ is tamely ramified. Let $f$ be an element of $\mathcal{H}(X_E)^{G_{F,der}}$; $f$ is then entirely determined by its value on some given element of $Ch_c$. In particular, theorem \ref{th2} holds for groups of type $A_{2n}$.
\end{cor}

By proposition \ref{anismax}(4), every subset of strongly orthogonal roots of $\Phi$ satisfies {\bf (C1)}; by corollary \ref{gtransa2}, proposition \ref{noncond} and an easy induction, $f$ is then entirely determined by its values on the set $Ch_\emptyset$ of chambers of $X_E$ whose geometric realization is contained in $\mathcal{B}_F$. On the other hand, by proposition \ref{chsol}, $Ch_c$ is the only $G_{F,der}$-orbit of chambers satisfying that condition and on which $f$ can be nonzero, hence $f$ is entirely determined by its values on $Ch_c$. In particular, $\mathcal{H}(X_E)^{G_{F,der}}$ is of dimension at most $1$. As in \cite{bc}, theorem \ref{th2} follows. $\Box$

Note that we did not need to determine precisely the support of the elements of $\mathcal{H}(X_E)^{G_{F,der}}$ to prove the above corollary, so we do it now.

For every apartment $A$ of $X_E$, we denote by $Ch_A$ the set of chambers of $A$.

\begin{cor}\label{suppa2n}
Assume $\Phi$ is of type $A_{2n}$ for some $n$ and $E/F$ is tamely ramified. Then assuming $\mathcal{H}(X_E)^{G_{F,der}}$ contains nonzero elements, their support is the union of $Ch_c$ and of the $Ch_A$, with $A$ being a  $\Gamma$-stable apartment of $X_E$ whose geometric realization is not contained in $\mathcal{B}_F$ and such that every facet of maximal dimension of $A^\Gamma$ is a facet of some element of $Ch_c$.
\end{cor}

Let $A'$ be any $\Gamma$-stable apartment of $X_E$, and let $\Sigma'$ be a set of strongly orthogonal roots of $\Phi$ corresponding to the $F$-anisotropy class of the $E$-split maximal torus associated to $A'$. By proposition \ref{anismax}(4), every nonempty subset of $\Sigma'$ satisfies {\bf (C1)}. Let $D$ be a facet of maximal dimension of $A'^\Gamma$; we prove by induction on the cardinality of $\Sigma'$ that for every nonzero $f\in\mathcal{H}(X_E)^{G_{F,der}}$, assuming such an $f$ actually exists, $f$ is nonzero on the set of chambers of $A'$ containing $D$ only if $D$ is contained in some chamber of $Ch_c$, and $f$ is then constant on that set. Let $\Sigma$, $A$, $C$, $C'$ and $C''$ be defined as in proposition \ref{noncond} relatively to $\Sigma'$ and $A'$; by that proposition and the harmonicity condition, $f(C')=0$ if and only if $f(C)+f(C'')=0$, and we have:
\[f(C')=\frac{f(C)+f(C'')}{1-q}.\]

On the other hand, if $\Sigma'$ is a singleton, then $C$ and $C''$ are two adjacent chambers in $Ch_\emptyset$, hence by definition of $Ch_c$, one of them can belong to $Ch_c$ only if the wall separating them, which is $D$, is such that $R(D)$ is not contained in any wall of $\mathcal{B}_F$, and in such a case, $R(C)$ and $R(C'')$ are contained in the same chamber of $\mathcal{B}_F$; by proposition \ref{chsol}, if, say, $C$ belongs to $Ch_c$, then $C''\not\in Ch_c$. Hence by corollary \ref{chsol1}, we have $f(C)+f(C'')\neq 0$, which implies $f(C')\neq 0$, if and only if $D$ is contained in some chamber of $Ch_c$. On the other hand, the value of $f$ on $C'$ is then always equal to the constant value of $f$ on $Ch_c$ multiplied by $\frac 1{1-q}$, hence nonzero. The fact that $f$ is then nonzero on the whole set of chambers of $A'$ is a consequence of corollary \ref{gtransa2}.

Assume now $\Sigma'$ contains at least two elements. By induction hypothesis, we have $f(C)=f(C'')$, and they are nonzero if and only if both $C$ and $C''$ contain a facet $D'$ of maximal dimension of $A^\Gamma$ contained in some chamber of $Ch_c$. Hence if $f(C')\neq 0$, $D$ must be contained in $Ch_c$ and we have $f(C')=\frac 2{1-q}f(C)$. Conversely, if $D$ is contained in some chamber of $Ch_c$, then it is contained in some $D'$ satisfying the same condition and $f(C')$ is then nonzero. As in the previous case, we use corollary \ref{gtransa2} to obtain that $f$ is then nonzero on the whole set of chambers of $A'$. $\Box$

\section{The spherical part}

In this section, we prove theorem \ref{th2} when $\Phi$ is not of type $A_{2n}$ for any $n$. From now on until the end of the paper, we assume that $E/F$ is tamely ramified.

Let $\Sigma_a$ be a subset of strongly orthogonal roots of $\Phi$ satisfying the conditions of proposition \ref{anismax}, let $\mathcal{A}$ be a $\Gamma$-stable apartment of $\mathcal{B}_E$ whose associated torus is of $F$-anisotropy class $\Sigma_a$, let $T$ be the $E$-split maximal torus associated to $\mathcal{A}$, and let $D$ be a facet of $X_E$ whose geometric realization is a facet of maximal dimension of $\mathcal{A}^\Gamma$; we denote by $Ch_D$ the set of chambers of $X_E$ containing $D$. First we prove that the elements of $\mathcal{H}(X_E)^{G_{F,der}}$ are entirely determined by their restrictions to $Ch_D$ for a suitably chosen $D$, then we prove that the space of restrictions of the elements of $\mathcal{H}(X_E)^{G_{der,F}}$ to $Ch_D$ is of dimension at most $1$.

To achieve that, we will continue to restrict our harmonic cochains to smaller sets. The general strategy is the following one: starting with the whole set $Ch_E$, we successively prove that we only have to consider the following subsets:
\begin{itemize}
\item the subset $Ch_D$ of the elements of $Ch_E$ which contain $D$;
\item the subset $Ch_{D,a}$ of the elements of $Ch_D$ whose $F$-anisotropy class is (up to conjugation) $\Sigma_a$;
\item the subset $Ch_{D,a,L}$ of the elements of $Ch_{D,a}$ contained in some $\Gamma$-stable apartment $A$ of $X_E$ whose associated torus is contained in some given reductive subgroup $L$ of $G$ (namely, the one of proposition \ref{acltorus});
\item the subset $Ch_{D,a,L,C_0}$ of the elements of $Ch_{D,a,L}$ of the form $uC_0$ for some given $\Gamma$-fixed chamber $C_0$ of $X_E$ containing $D$, where $u$ is a product of elements of the root subgroups of $L$ which correspond to elements of $\Sigma_a$.
\end{itemize}
Finally, we compute explicitely the restrictions to $Ch_{D,a,L,C_0}$ of our harmonic cochains; by proposition \ref{clcoh}, that set happens to be in 1-1 correspondence with some cohomology group which is easier to study.

\subsection{Some preliminary results}

We choose $D$ arbitrarily for the moment.
First we prove the following results:

\begin{prop}\label{spv}
Assume $D$ is a single vertex $x$; $x$ is then a special vertex of $X_E$.
\end{prop}

(See section $2$ for the definition of a special vertex.)

By eventually conjugating it by some element of $G_{F,der}$ we can always assume that $x\in\mathcal{A}_0$. The above statement can then be rewritten in terms of concave functions the following way: let $\Sigma$ be a set of strongly orthogonal roots of $\Phi$ conjugated to $\Sigma_a$. Assume the cardinality of $\Sigma$ is equal to the rank $d$ of $\Phi$ and let $f$ be a concave function from $\Phi$ to $\frac 12{\mth{Z}}$ such that $f(\alpha)\in{\mth{Z}}+\frac 12$ for every $\alpha\in\pm\Sigma$ (this property corresponds to the fact that $D$ is a $\Gamma$-fixed facet of maximal dimension of an apartment of $F$-anisotropy class $\Sigma_a$); we then have $f(\alpha)+f(-\alpha)=0$ for every $\alpha\in\Sigma$.

Let $f'$ be the element of $Hom(X^*(T_0),{\mth{Q}})$ which coincides with $f$ on $\pm\Sigma$; for every $\alpha\in\Phi$, we have $f'(\alpha)=\alpha(x)$ (remember that $\mathcal{A}_0=X_*(T_0)\otimes{\mth{R}}$). We then have $f(\alpha)+f(-\alpha)=0$ for every $\alpha\in\Phi$ if and only if $f$ coincides with $f'$ on $\Phi$, and by definition of $f$, this is the case if and only if the image of $f'$ is contained in $\frac 12{\mth{Z}}$. Proposition \ref{spv} is then an immediate consequence of the following proposition:

\begin{prop}\label{f12z}
The function $f$ being defined as above, the image of $f'$ is actually contained in $\frac 12{\mth{Z}}$.
\end{prop}

Let $\beta_1,\dots,\beta_d$ be the elements of $\Sigma$, and let $\alpha$ be any element of $\Phi$, which we can assume to be different from the $\pm\beta_i$, $i\in\{1,\dots,d\}$. Write $\alpha=\sum_{i=1}^d\lambda_i\beta_i$, the $\lambda_i$ being elements of ${\mth{Q}}$; we then have $f'(\alpha)=\Sigma_{i=1}^d\lambda_if'(\beta_i)$. On the other hand, for every $i$, we have:
\[<\alpha,\beta_i^\vee>=\lambda_i<\beta_i,\beta_i^\vee>=2\lambda_i,\]
hence $\lambda_i\in\frac 12{\mth{Z}}$. Let $(.,.)$ be a nontrivial $W$-invariant scalar product on $X^*(T)\otimes{\mth{Q}}$; we also have:
\begin{equation}\label{winv}(\alpha,\alpha)=\sum_{i=1}^d\lambda_i^2(\beta_i,\beta_i).\end{equation}
We now consider the possible cases. To simplify notations, we can assume that the nonzero $\lambda_i$ are the ones with the lowest indices, and are positive (because we can always replace some of the $\beta_i$ with their opposites by simply conjugating $\Sigma$ by a product of reflections $s_{\beta_i}$).
\begin{itemize}
\item Assume first $\Phi$ is simply-laced. Then $(\alpha,\alpha)$ and the $(\beta_i,\beta_i)$ are all equal to each other, and there is only one possibility: $\lambda_i=\frac 12$ for $1\leq i\leq 4$ and $\lambda_i=0$ for $i>4$; we then obtain:
\[f'(\alpha)=\frac 12(f'(\beta_1)+f'(\beta_2)+f'(\beta_3)+f'(\beta_4))\in\frac 12(2+{\mth{Z}})=\frac 12{\mth{Z}}.\]
as desired.
\item Assume now $\Phi$ is not simply-laced and every $\beta_i$ such that $\lambda_i\neq 0$ is long; since there are then at least two long $\beta_i$ orthogonal to each other, we cannot be in the case $G_2$ here. If $\alpha$ is long as well, we are reduced to the previous case. If $\alpha$ is short, then for every $i$, $(\alpha,\alpha)=\frac 12(\beta_i,\beta_i)$ and there is again only one possibility: $\lambda_1=\lambda_2=\frac 12$ and $\lambda_i=0$ for $i>2$; we then obtain:
\[f'(\alpha)=\frac 12(f'(\beta_1)+f'(\beta_2))\in\frac 12(1+{\mth{Z}})=\frac 12{\mth{Z}},\]
\item Assume now $\Phi$ is not simply-laced and some of the $\beta_i$ such that $\lambda_i\neq 0$ are short. We first treat the case $G_2$; in this case, assuming $\beta_1$ is short and $\beta_2$ is long, we have $(\beta_2,\beta_2)=3(\beta_1,\beta_1)$, and $3\lambda_2^2+\lambda_1^2$ is either $1$ (if $\alpha$ is short) or $3$ (if $\alpha$ is long). In the first (resp. second) case, it implies $\lambda_1=\frac 12$ and $\lambda_2=\frac 12$ (resp. $\lambda_1=\frac 12$ and $\lambda_2=\frac 32$), and in both cases, we obtain $f'(\alpha)=\lambda_1f'(\beta_1)+\lambda_2f'(\beta_2)\in\frac 12{\mth{Z}}$.
\item Assume now $\Phi$ is not of type $G_2$, not simply-laced and $\beta_1$ is short. First assume $\beta_1$ is the only short $\beta_i$ such that $\lambda_i\neq 0$. If $\alpha$ is long, this is only possible if there are three nonzero $\lambda_i$, $\lambda_1=1$ and $\lambda_2=\lambda_3=\frac 12$, and we then have:
\[f'(\alpha)=f'(\beta_1)+\frac 12(f'(\beta_2)+f'(\beta_3))\in\frac 12{\mth{Z}}+\frac 12(1+{\mth{Z}})=\frac 12{\mth{Z}}.\]

Assume now $\alpha$ is short, still with only one of the $\beta_i$ such that $\lambda_i\neq 0$ being short. We deduce from the relation (\ref{winv}) that we must have:
\[1=\lambda_1^2+2\sum_{i=2}^d\lambda_i^2.\]
Since $\lambda_2$ is nonzero, we must have $\lambda_1=\frac 12$. But then the right-hand side of the above equality belongs to $\frac 14+\frac 12{\mth{Z}}$, hence cannot be equal to $1$. We are then in an impossible case.
\item Assume finally that at least two of the $\beta_i$ such that the $\lambda_i$ are nonzero are short. By lemma \ref{twoshort}, this is possible only if $\Phi$ is of type $C_d$ with $d\geq 4$. On the other hand, by proposition \ref{sigmaa}, if $\Phi$ is of type $C_d$, $\Sigma_a$ contains only long roots. Hence this case is impossible too.
\end{itemize}
$\Box$

Now we consider the cases where $D$ is of positive dimension or in other words the ones where $\Sigma_a$ contains less than $d$ elements; as we have already seen, these cases are $A_d$, $d>1$ odd (remember that we rule out the case $A_{2n}$ in this whole section), $D_d$, $d=2n+1$ odd, and $E_6$.

\begin{prop}\label{spv1}
Assume $D$ is of positive dimension. Then its vertices are all special.
\end{prop}

When $G$ is of type $A_d$, we see on the table of page $29$ of \cite{bt} that every vertex of $X_E$ is special and the result of the proposition is trivial; we then only have to consider the cases $D_{2n+1}$ and $E_6$.

Assume $\Phi$ is of type $D_{2n+1}$. The facet $D$ is then  of dimension $1$, and we have, for example, $\Sigma_a=\{\varepsilon_2\pm\varepsilon_3,\dots,\varepsilon_{2n}\pm\varepsilon_{2n+1}\}$; $\Sigma_a$ is then contained in the Levi subsystem $\Phi'$ of type $D_{2n}$ of $\Phi$ generated by $\alpha_2,\dots,\alpha_{2n+1}$. Let $Y$ be the subgroup of $X^*(T_0)$ generated by $\Sigma_a$, let $f$ be a concave function from $\Phi\cap Y$ to $\frac 12{\mth{Z}}$ such that $f(\alpha)\in{\mth{Z}}+\frac 12$ for every $\alpha\in\pm\Sigma$, and let $f'$ be the element of $Hom(Y,{\mth{Q}})$ which coincides with $f$ on $\pm\Sigma$; if we extend $f'$ linearly to $X^*(T_0)\otimes{\mth{Q}}$ by choosing $f'(\varepsilon_1)$ arbitrarily in $\frac 12{\mth{Z}}$, we obtain on $\Phi$ a concave function satisfying $f'(\alpha)+f'(-\alpha)=0$ for every $\alpha\in\Phi$ and associated to some vertex of $\mathcal{A}^\Gamma$, and it is easy to check that every vertex of $\mathcal{A}^\Gamma$ is associated to such a concave function, hence special.

Assume now $\Phi$ is of type $E_6$. The facet $D$ is then of dimension $2$, and, up to conjugation, we have $\Sigma_a=\{\alpha_2,\alpha_3,\alpha_5,\alpha_2+\alpha_3+2\alpha_4+\alpha_5\}$; $\Sigma_a$ is then contained in the Levi subsystem $\Phi'$ of type $D_4$ of $\Phi$ generated by $\alpha_2,\dots,\alpha_5$. Once again, $Y$ and $f'$ being defined as in the previous case, we can extend $f'$ linearly to $X^*(T_0)\otimes{\mth{Q}}$ by choosing $f'(\alpha_1)$ and $f'(\alpha_6)$ arbitrarily in $\frac 12{\mth{Z}}$, and we conclude similarly as above. $\Box$

Now assume the geometric realization of $D$ is contained in $\mathcal{A}_0$; let $f_{D,E}$ be the concave function associated to $D$ (as a facet of $X_E$; we have to specify here since $D$ may be a vertex of both $X_E$ and $X_F$). The following corollary follows immediately from propositions \ref{spv} and \ref{spv1}:

\begin{cor}\label{spv2}
Let $\alpha$ be an element of $\Phi$ which is a linear combination of elements of $\Sigma$. Then $f_{D,E}(\alpha)+f_{D,E}(-\alpha)=0$.
\end{cor}

Let now $A$ be any $\Gamma$-stable apartment of $X_E$ of anisotropy class $\Sigma_a$ containing $D$.

\begin{prop}\label{appar}
The subvomplex $A^\Gamma$ is isomorphic to any apartment of a building of type $A_r$, where $r$ is its dimension.
\end{prop}

Let $\Phi_D$ be the root system of $K_{D,E}/K_{D,E}^0$, viewed as a root subsystem of $\Phi$, and let $S_0\subset T_0$ be the intersection of the $Ker(\alpha)$, $\alpha\in\Phi$. Assume there exists a field $F_1$ (not necessarily related to $F$ in any way) on which $G$ is defined and a $F_1$-inner form $G'$ of $G$ such that $S_0$ is precisely the maximal $F_1$-split subtorus of $T_0$ in $G'$. Then by \cite[corollary 5.8]{bot}, the set of nonzero restrictions to $S_0$ of the elements of $\Phi$ is a root system, which implies that $\mathcal{A}^\Gamma$ is isomorphic as an affine simplicial complex to an apartment of a building of the same type as that root system.

Now we check that the above assumption is true. If $r=0$, then $\mathcal{A}^F$ consists of a single vertex and the result is trivial; assume $r>0$. We then obtain, with the help of \cite[section 17]{spr}:
\begin{itemize}
\item when $G$ is of type $A_{2n-1}$, we can take $F_1=F$, and $G'$ is then, up to isogeny, the group $GL_n(\mathcal{\mathcal{D}})$, where $\mathcal{D}$ is a quaternionic division algebra over $F$. The group $G'$ is then of relative type $A_{n-1}$, hence $A^\Gamma$ is isomorphic to an apartment of type $A_{n-1}$;
\item when $G$ is of type $D_{2n+1}$, $G'$ is isogeneous to a special orthogonal group in $4n+2$ variables defined by a quadratic form of index $1$, and $F_1$ is any field on which such a quadratic form exists (for example ${\mth{R}}$, in which case $G_{\mth{R}}=SO_{4n+1,1}({\mth{R}})$, but not $F$ this time); $G'$ is then of relative type $A_1$;
\item when $G$ is of type $E_6$, $G'$ is the case $(1,6)$ of \cite[proposition 17.7.2]{spr}, and $F_1$ is any field on which such an inner form of $G$ exists (again, $F_1=F$ does not work, but $F_1={\mth{R}}$ does according to the classification of \cite{tits}); $G'$ is then of relative type $A_2$.
\end{itemize}
The proposition is now proved. $\Box$

\subsection{Restriction to $Ch_D$}

Now we go back to the proof of theorem \ref{th2}. Let $A$ and $A'$ be two $\Gamma$-stable apartments of $X_E$ corresponding to tori of $F$-anisotropy class $\Sigma_a$. By proposition \ref{fpconj}, there exists $g\in G_{F,der}$ such that $gA^\Gamma=A'^\Gamma$. If $A^\Gamma$ (resp. $A'^\Gamma$) is a single vertex $x$ (resp. $x'$), we have $Ch_{x',a}=gCh_{x,a}$, and the $G_{F,der}$-invariance of the elements of $\mathcal{H}(X_E)^{G_{F,der}}$ implies that their restrictions to $Ch_{x,a}$ and $Ch_{x',a}$ are related. Assume now $A^\Gamma$  and $A'^\Gamma$ are of dimension at least $1$. Then for every facet $D$ of $A^\Gamma$ of maximal dimension, $gD$ is a facet of $A'^\Gamma$ of maximal dimension, and the restriction to $Ch_{D,a}$ of every element of $\mathcal{H}(X_E)^{G_{F,der}}$ depends only on its restriction to $Ch_{gD,a}$ and conversely. To prove that $f\in\mathcal{H}(X_E)^{G_{F,der}}$ only depends on its restriction to $Ch_{D,a}$, we thus only need to prove that its restrictions to respectively $Ch_{D,a}$ and $Ch_{D',a}$, where $D'$ is any other facet of maximal dimension of $A^\Gamma$, are related as well.

\begin{prop}\label{cclink}
Let $D,D'$ be two facets of maximal dimension of $A^\Gamma$. Then $D$ and $D'$ are $G_{F,der}$-conjugates.
\end{prop}

The result is trivial when $A^\Gamma$ consists of a single vertex; assume it is not the case. By an obvious induction we only have to prove the proposition wuen $D$ and $D'$ are neighboring each other. Let $D''$ be their common facet of codimension $1$, and set ${\mth{G}}_D=K_{D,E}/K_{D,E}^1$; define similarly ${\mth{G}}_{D'}$ and ${\mth{G}}_{D''}$. It is easy to check (details are left to the reader) that:
\begin{itemize}
\item when $G$ is of type $A_{2n-1}$, ${\mth{G}}_D$ and ${\mth{G}}_{D'}$ are of type $A_1\times A_1$ and ${\mth{G}}_{D''}$ of type $A_3$;
\item when $G$ is of type $D_{2n+1}$, ${\mth{G}}_D$ and ${\mth{G}}_{D'}$ are of type $D_{2n}$ and ${\mth{G}}_{D''}$ of type $D_{2n+1}$;
\item when $G$ is of type $E_6$, ${\mth{G}}_D$ and ${\mth{G}}_{D'}$ are of type $D_4$ and ${\mth{G}}_{D''}$ of type $D_5$.
\end{itemize}
Hence in every case (including the first one, remember that $D_3=A_3$ and $D_2=A_1\times A_1$), ${\mth{G}}_D$ and ${\mth{G}}_{D'}$ are of type $D_{2r}$ and ${\mth{G}}_{D''}$ of type $D_{2r+1}$ for some $r$; we thus are reduced to the case $D_{2n+1}$, and we may assume $G$ is $SO'_{4n+2}$. It is then easy to check that, depending on the case, $D$ and $D'$ are conjugated by some $G_{F,der}$-conjugate of either:
\[\left(\begin{array}{cccccccccc}&&&&&&&&&1\\&1\\&&\ddots\\&&&1\\&&&&&1\\&&&&1\\&&&&&&1\\&&&&&&&\ddots\\&&&&&&&&1\\1\end{array}\right),\]
or:
\[\left(\begin{array}{cccccccccc}&&&&&&&&&\varpi_F^{-1}\\&1\\&&\ddots\\&&&1\\&&&&&1\\&&&&1\\&&&&&&1\\&&&&&&&\ddots\\&&&&&&&&1\\\varpi_F\end{array}\right),\]
The proposition is then proved. $\Box$

\begin{cor}
Let $f$ be an element of $\mathcal{H}(X_E)^{G_{F,der}}$ and let $D$ be the facet of $X_E$ defined as in proposition \ref{cclink}. Then the restriction of $f$ to $Ch_a^0$ is entirely determined by its restriction to $Ch_{D,a}$.
\end{cor}

Let $D'$ be another facet of maximal dimension of some $A'^\Gamma$. If $g\in G_{F,der}$ is such that $A'^\Gamma=gA^\Gamma$, the restrictions of $f$ to $Ch_{gD,a}$ and $Ch_{D',a}$ depend only on each other by the previous proposition and $G_{F,der}$-invariance, and by $G_{F,der}$-invariance again, its restrictions to $Ch_{D,a}$ and $Ch_{gD,a}$ are also linked. The result follows. $\Box$

\subsection{The harmonic cochains on $Ch_{D,a}$}

Now we prove that, for some convenient $D$, the dimension of the space of the restrictions to $Ch_{D,a}$ of the elements of $\mathcal{H}(X_E)^{G_{F,der}}$ is at most $1$.

We fix $D$ arbitrarily for the moment among the possible ones contained in $A_{0,E}$. Let $\Phi_D$ be the root system of $K_{D,E}/K_{D,E}^0$ relative to $K_{T_0,E}/K_{T_0,E}^0$, viewed as a root subsystem of $\Phi$.

Let $\beta_1,\dots,\beta_r$ be the elements of some fixed representative of $\Sigma_a$, and let $L$ be the subgroup of $G$ generated by $T_0$ and the $U_{\pm\beta_i}$, $i=1,\dots,r$; by proposition \ref{acltorus} we know that every $E$-split maximal $F$-torus of $G$ of $F$-anisotropy class $\Sigma_a$ is $G_{F,der}$-conjugated to some maximal torus of $L$. Hence we can replace $Ch_{D,a}$ by the subset $Ch_{D,a,L}$ of the elements $C\in Ch_{D,a}$ contained in a $\Gamma$-stable apartment of $X_E$ whose associated $E$-split maximal torus is also contained in $L$.

\begin{prop}\label{oppch}
Let $C$ be any chamber of $A_E$ containing $D$; there exist chambers $C_0,C'_0$ of $A_{0,E}$ containing $D$ and corresponding to opposite Borel subgroups of $K_{D,E}/K_{D,E}^0$ and an element $u\in L_{E,der}\cap K_{C'_0,E}$ such that $C=uC_0$.
\end{prop}

Since $T$ and $T_0$ are both contained in $L$, there exists $h\in L_{E,der}$ such that $hT_0h^{-1}=T$ and $hD=D$, hence $h\in K_D\cap L_{E,der}$. Therefore, we have $h^{-1}C=C_0$ for some chamber $C_0$ of $A_{0,E}$ containing $D$. Moreover, we have:

\begin{lemme}\label{uuprime}
Let $B,B'$ be two opposite Borel subgroups of $L_{E,der}$ containing $T_0$ and let $U,U'$ be their respective unipotent radicals. Then $T_0$ and $T'$ are conjugated by some element $h=uu'$ of $U_EU'_E$, where $U_E$ (resp. $U'_E$) is the group of $E$-points of $U$ (resp. $U'$). Moreover, if $h\in L_{E,der}\cap K_{D,E}$, then $u$ and $u'$ also belong to $L_{E,der}\cap K_{D,E}$.
\end{lemme}

Let $h'$ be any element of $L_{E,der}$ such that $h'T_0h'^{-1}=T'$. Using the Bruhat decomposition of $L_E$ (see for example \cite[16.1.3]{spr}) and the fact that both $B$ and $B'$ contain $T_0$, we can write $h'=unu''$, with $u,u''\in U_E$ and $n\in N_{L_{E,der}}(T_0)$, and we can even assume that $u''$ belongs to $n^{-1}U'_En$, hence $u'=nu''n^{-1}\in U'_E$; if we set $h=h'n^{-1}$, then $h=uu'$ satisfies $hT_0h^{-1}=h'T_0h'^{-1}=T'$, as required.

Assume now $h\in K_{D,E}\cap L_{E,der}$. Since the intersections of $K_{D,E}$ with respectively $U_E$ and $U'_E$ are products of the intersections with $K_{D,E}$ of the root subgroups respectively contained in these two subgroups, and since these two sets of root subgroups are disjoint, we obtain that $u$ and $u'$ belong to $K_{D,E}\cap L_{E,der}$ as well. $\Box$

Note that since $T_0$ is split and $T$ is of anisotropy class $\Sigma_a$, the element $n$ of $N_{L_{E,der}}(T_0)$ used in the above proof always corresponds to the element of the Weyl group of $L$ relative to $T_0$ which sends every root of $L$, hence also every root of $K_{D,E}$ by linearity, to its opposite (more precisely, $w$ is the product of $d$ copies of $w_0$, where $w_0$ is the unique nontrivial element of the Weyl group of $SL_2$). Since $h'$ has been chosen arbitrarily, we obtain that every $h\in L_{E,der}$ such that $hT_0h^{-1}=T'$ satisfies $h\in U_EU'_ET_{0,E}$, and that when $h$ belongs to $K_{D,E}$, its three components also belong to $L_{E,der}\cap K_{D,E}$.

Now we prove proposition \ref{oppch}. According to lemma \ref{uuprime} and the previous remark, for every choice of $U_E$ and $U'_E$, we have $C=uu'C_0$ for some $C_0$, some $u\in L_{E,der}\cap U_E$ and some $u'\in L_{E,der}\cap U'_E$. Hence for every $C_0$, if we choose $U_E$, $U'_E$ in such a way that the image of $u'$ in $K_{D,E}/K_{D,E}^0$ belongs to the Borel subgroup of $K_{D,E}/K_{D,E}^0$ corresponding to $C_0$, or in other words that $u'\in K_{C_0,E}$, we have in fact $C=uC_0$. Let then $C'_0$ be the unique chamber of $A_{0,E}$ containing $D$ and corresponding to a Borel subgroup of $K_{D,E}/K_{D,E}^0$ opposite to the previous one; by definition of $U_E$ and by lemma \ref{uuprime}, we must then have $u\in L_{E,der}\cap K_{C'_0,E}$, as required. $\Box$

For every $\alpha\in\Phi$, let $u_\alpha$ be a group isomorphism between the additive group $E$ and $U_\alpha$ compatible with the valued root datum $(G,T_0,(U_\alpha)_{\alpha\in\Phi},(\phi_\alpha)_{\alpha\in\Phi})$; for every one-parameter subgroup $\xi$ of $T_0$, we then have, for every $x,y\in E^*$:
\[\xi(x)u_\alpha(y)\xi(x)^{-1}=u_\alpha(x^{<\alpha,\xi>}y),\]
where $<.,.>$ denotes the usual pairing between $X^*(T_0)$ and $X_*(T_0)$.

\begin{cor}
There exist elements $\lambda_1,\dots,\lambda_r\in\mathcal{O}_E^*$ such that the element $u$ of lemma \ref{uuprime} is of the form $u=\prod_{i=1}^ru_{\beta_i}(\varpi_E^{2f_{D,E}(\beta_i)}\lambda_i)$ for some choice of $\Sigma_a=\{\beta_1,\dots,\beta_n\}$.
\end{cor}

(Note that since the elements of $\Sigma_a$ are strongly orthogonal, the root subgroups $U_{\beta_i}$ commute, hence the above product can be taken in any order.)

Assume $\Sigma_a$ has been chosen in such a way that for every $\beta\in\Sigma_a$, the root subgroup $U_\beta$ of $G$ is contained in the group $U_E$ defined as in proposition \ref{oppch}; since, using lemma \ref{negsorth}, we can always replace some of its elements with their opposites, this is always possible. 

Since $u$ is unipotent, it belongs to the derived group $L_{E,der}$ of $L_E$, and we can work componentwise. Write $u=u_1\dots u_r$, where for every $i$, $u_i$ belongs to $U_{\beta_i}$. For every $i$, $u_i$ then belongs to $K_{C'_0}$ but not to $K_{C_0}$, hence is of the form $u_{\beta_i}(\varpi_E^{2f_{D,E}(\beta_i)}\lambda_i)$ for some $\lambda_i\in\mathcal{O}_E^*$; the result follows. $\Box$

Note that for every $i$, since $f_{D,E}(\beta_i)\in{\mth{Z}}+\frac 12$, $\varpi_E^{2f_{D,E}}(\beta_i)\lambda_i$ cannot be an element of $F$.

For every chamber $C_0$ of $A_0$ containing $D$, let $Ch_{D,a,L,C_0}$ be the subset of the $C\in Ch_{D,a,L}$ such that, with $C'_0$ being defined as in proposition \ref{oppch}, $C=uC_0$ for some $u\in K_{C'_0,E}$. We deduce from the previous corollary that,$Ch_{D,a,L}$ is the union of the $Ch_{D,a,L,C_0}$, with $C_0$ being such that the corresponding Borel subgroup of $K_{D,E}/K_{D,E}^0$ contains every root subgroup associated to the elements of some fixed representative of $\Sigma_a$.

Now we fix arbitrarily such a chamber $C_0$. For every $\lambda_1,\dots,\lambda_r\in\mathcal{O}_E^*$, where $r$ is the cardinality of $\Sigma_a$, let $C(\lambda_1,\dots,\lambda_r)$ be the chamber $\prod_{i=1}^ru_{\beta_i}(\varpi_E^{2f_{D,E}(\beta_i)}\lambda_i)C_0$, where the $\beta_i$ are the elements of $\Sigma_a$. The chamber $C(\lambda_1,\dots,\lambda_r)$ only depends on the classes mod $\mathfrak{p}_E$ of the $\lambda_i$, hence by a slight abuse of notation we can consider them as elements of the residual field $k_E^*=k_F^*$.

\begin{prop}\label{chdalconj}
The subsets $Ch_{D,a,L,C_0}$ of $Ch_{D,a,L}$ are all $K_{D,E}\cap G_{F,der}$-conjugates.
\end{prop}

Let $C_0,C'_0$ be two chambers of $A_{0,E}$ containing $D$, and let $C$ be any element of $Ch_{D,a,L,C_0}$; there exists then $n\in N_{G_{E,der}}(T_0)\cap K_{D,E}$ such that $nC_0=C'_0$. Let $g\in G_{E,der}$ be such that $gT_0g^{-1}=T$ and $gC_0=C$, and set $n'=gng^{-1}$; the chamber $C'=n'C$ then belongs to $Ch_{D,a,L,C'_0}$. We thus only have to prove that $C'$ is a $G_{F,der}$-conjugate of some element of $Ch_{D,a,L,C_0}$.

By an obvious induction it is enough to prove the result when $C$ and $C'$ are neighboring each other. Assume first $n$ is the reflection associated to some element $\beta_i$ of $\Sigma_a$; then $C'$ is of the form $C'=C(\lambda_1,\dots,\lambda_{i-1},\mu,\lambda_{i+1},\dots,\lambda_r)$ for some $\mu\in k_F^*$ distinct from $\lambda_i$, and the result follows.

Next we prove the following lemma:

\begin{lemme}
Let $n'$ be any element of $N_{G_E}(T)\cap K_{D,E}$. Then $L'=nLn^{-1}$ is $F$-split.
\end{lemme}

If we assume that $L'$ is defined over $F$, then it is $F$-split by lemma \ref{lsplit}. Therefore, we only have to prove that $L'$ is defined over $F$. Let $w$ be the element of the Weyl group of $G/T$ corresponding to $n$: since $T$ is of $F$-anisotropy class $\Sigma_a$, for every $\alpha$ belonging to the root system $\Phi_{L',T}$ of $L'$ relative to $T$, we have $\gamma(\alpha)=-\alpha$, hence for every root $\beta$ of $L'/T$, $w(\beta)$ is a root of $L'/T$ and:
\[\gamma(w(\beta))=-w(\beta)=w(-\beta)\]
is also a root of $L'/T$. Hence $L'$ is $\Gamma$-stable, hence defined over $F$. $\Box$ 

According to this lemma, replacing $L$ by some $K_{D,F}$-conjugate if needed, we see that the result of proposition \ref{chdalconj} holds when $n$ is the reflection associated to any conjugate of any element of $\Sigma_a$. Since by \cite[\S I, proposition 11]{bou}, two roots of $\Phi$ of the same length are always conjugates, proposition \ref{chdalconj} holds when $\Sigma_a$ contains roots of every length.

Now assume $n$ is any element of $K_{D,E}\cap N_{G_{E,der}}(T)$, and let $g$ be an element of $G_{E,der}$ such that $gC_0=C$. We then have:
\[C'=ngC_0=g(g^{-1}ng)C_0=gn_0C_0,\]
where $n_0=g^{-1}ng$ is an element of $K_{D,E}\cap N_{G_{E,der}}(T_0)$, which we can assume to be in $G_{F,der}$ since $T_0$ is $F$-split. On the other hand, by lemma \ref{uuprime}, $g$ is of the form $(n_0un_0^{-1})(n_0u'n_0^{-1})$, with $u'\in K_{C_0}$ and $u$ being of the form $u=\prod_{i=1}^ru_{\beta_i}(\varpi_E^{2f_{D,E}(\beta_i)}\mu_i)$, with $\mu_1,\dots,\mu_r$ being elements of $k_F^*$. Hence $n_0^{-1}C'=uu'C_0$ belongs to $C(\mu_1,\dots,\mu_r)$ and $C'$ is then $G_{F,der}$-conjugated to some element of $Ch_{D,a,L C_0}$.

It remains to prove that every $n$ such that $C$ and $nC$ are neighboring each other is a $G_F$-conjugate of some element of $K_{D,E}\cap N_{G_{E,der}}(T)$. It is of course true when $n$ is a representative of the reflection associated to some conjugate of some element of $\Sigma_a$, hence since $\Sigma_a$ always contains some long roots, we only have to consider the case where $\Sigma_a$ contains only long roots and $n$ is a representative of the reflection associated to some short root of $\Phi$.
Since $\Phi$ is then not simply-laced, we deduce from the remark following lemma \ref{sigma1} that $\Phi$ satisfies the equivalent conditions of that lemma, which implies in particular that $\Phi_D$ is of the same rank as $\Phi$, and proposition \ref{spv} then implies that $\Phi_D=\Phi$. Hence $n$ is a representative of the reflection associated to some element of $\Phi_D$, and the result follows. $\Box$

By the above proposition, to prove theorem \ref{th2}, we only have to prove that the space of the restrictions of the elements of $\mathcal{H}(X_E)^{G_{F,der}}$ to $Ch_{D,a,L,C_0}$ is of dimension at most $1$ for some given $D,a,L,C_0$. We start by dividing that set into $L_{F,der}$-conjugacy classes, which happen to be easier to handle than the full $G_{F,der}$-conjugacy classes. 

\begin{prop}\label{clcoh}
The $L_{F,der}$-conjugacy classes of elements of $Ch_{D,a,L,C_0}$ are in 1-1 correspondence with the elements of the cohomology group $H^1(\Gamma,K_{T\cap L_{E,der}})$. Moreover, that group is isomorphic to $({\mth{Z}}/2{\mth{Z}})^r$, where $r$ is the cardinality of $\Sigma_a$.
\end{prop}

First we compute $H^1(\Gamma,K_{T\cap L_{E,der}})$. It is obvious from the definitions that the group $L_{E,der}$ is $F$-anisotropic, hence $T\cap L_{E,der}$ is nothing else than the $F$-anisotropic component of $T$. Let $\xi$ be any $1$-parameter subgroup of $T\cap L_{E,der}$; its intersection with $K_{T\cap L_{E,der}}$ is $\xi(\mathcal{O}_E^*)$. On the other hand, since $Im(\xi)$ is contained in $T_a$, for every $\lambda\in\mathcal{O}_E^*$, we have $\gamma(\xi(\lambda))=\xi(\gamma(\lambda)^{-1})$. Hence $\xi(\lambda)$ defines a $1$-cocycle of $\Gamma$ if and only if $\gamma(\lambda)^{-1}\lambda=1$, or in other words if and only if $\lambda\in\mathcal{O}_F^*$. (Note that it does not mean that $\xi(\lambda)\in G_F$.) Moreover, $\xi(\lambda)$ defines a $1$-coboundary if and only if $\lambda=\gamma(\mu)\mu$ for some $\mu\in\mathcal{O}_E^*$, or in other words if and only if $\lambda$ is the norm of some element of $\mathcal{O}_E^*$, which is true if and only if its image in $k_F^*$ is a square. Since $X_*(T\cap L_{E,der})$ is generated by the coroots $\beta_1^\vee,\dots,\beta_r^\vee$ associated to the elements $\beta_1,\dots,\beta_r$ of $\Sigma_a$, we obtain that $H^1(\Gamma,K_{T\cap L_{E,der}})$ is isomorphic to a product of $r$ copies of $k_F^*/(k_F^*)^2\simeq{\mth{Z}}/2{\mth{Z}}$.

Now we prove some lemmas.

\begin{lemme}\label{lfdc}
Let $F'$ be the unique quadratic unramified extension of $F$. Then the elements of $Ch_{D,a,L,C_0}$ are all $L_{F',der}$-conjugates.
\end{lemme}

Let $T$ be a maximal torus of $G$ satisfying the conditions of lemma \ref{idanis}. By simply replacing $F$ by $F'$ in the discussion above, we obtain that when $\lambda$ is an element of $k_{E'}^*$, $\xi(\lambda)$ defines a $1$-cocycle in $K_{T\cap L_{E',der}}$ if and only if $\lambda\in k_{F'}^*$ and a $1$-coboundary if and only if $\lambda$ is the norm of an element of $k_{E'}^*$, which is true if and only if it is a square in $k_{F'}^*$. On the other hand, $[k_{F'}^*:k_F^*]=q+1$ is even, hence every element of $k_F^*$ is a square in $k_{F'}^*$. Lemma \ref{lfdc} follows. $\Box$

Set $E'=EF'$; $E'/E$ is then a quadratic unramified extension, and $X_E$ is a simplicial subcomplex of the building $X_{E'}$ of $G_{E'}$; the set $Ch_{D,a,L,C_0}$ is then a subset of the set of chambers of $X_{E'}$ containing $D$. Moreover, the extension $E'/F'$ is quadratic and tamely ramified.

\begin{lemme}\label{idanis}
Assume $C$ is an element of $Ch_E$; let $A$ be a $\Gamma$-stable apartment of $X_E$ containing $C$, and let $T$ be the corresponding $E'$-split torus of $G$. Then we can choose $A$ in such a way that $T$ is defined over $F$, $E$-split and that its $F$-anisotropic and $F'$-anisotropic components are identical.
\end{lemme}

Since $C\in Ch_E$, it is possible to choose $A$ in  such a way that $A$ is contained in $X_E$, which, since it is $\Gamma$-stable, ensures that $T$ is defined over $F$ and $E$-split.

Moreover, since $E'/F'$ is tamely ramified, the geometrical building $\mathcal{B}_{F'}$ of $G_{F'}$ is the set of $\Gamma$-fixed points of $\mathcal{B}_{E'}$, and in particular we have $\mathcal{B}_{F'}\cap\mathcal{B}_E=\mathcal{B}_F$. Hence the affine subspaces of $R(A)$ contained in respectively $\mathcal{B}_F$ and $\mathcal{B}_{F'}$ are the same, which proves that the $F$-anisotropic and $F'$-anisotropic components of $T$ have the same dimension. Since the second one is obviously contained in the first one, the lemma follows. $\Box$

Now we prove the first assertion of proposition \ref{clcoh}. For every $i\in\{1,\dots,r\}$, every $\lambda_1,\dots,\lambda_r\in\mathcal{O}_E^*$ and every $\mu\in\mathcal{O}_{E'}^*$ whose square is an element of $\mathcal{O}_E^*$, we have:
\[\beta_i^\vee(\mu)C(\lambda_1,\dots,\lambda_r)=C(\lambda_1,\dots,\mu^2\lambda_i,\dots,\lambda_r).\]
The chamber $C(\lambda_1,\dots,\lambda_r)$ being stable by $\beta_i^\vee(1+\mathfrak{p}_{E'})\subset K_{C(\lambda_1,\dots,\lambda_r),E'}$, we can assume $\mu\in\mathcal{O}_{F'}^*$, which implies $\mu^2\in\mathcal{O}_E^*\cap\mathcal{O}_{F'}^*=\mathcal{O}_F^*$. Since every element of $k_F^*$ is a square in $k_{F'}^*$, the image of $\mu^2$ in $k_F^*$ can be any element of $k_F^*$; we thus obtain that the subgroup $L$ of the elements of $(T_0)_F$ such that $tC(\lambda_1,\dots,\lambda_r)$ belongs to $Ch_{D,a,L,C_0}$, contains representatives of every element of $(k_F^*/(k_F^*)^2)^r\simeq H^1(\Gamma,K_{T\cap L_{E,der}})$; this proves that the set of $L_{F,der}$-conjugacy classes of elements of $Ch_{D,a,L,C_0}$ is in 1-1 correspondence with $H^1(\Gamma,K_{T\cap L_{F,der}})$, and proposition \ref{clcoh} is now proved. $\Box$

For every $h=(\sigma_1,\dots,\sigma_r)\in H^1(\gamma,K_{T,E}\cap L_{E,der})$, the $\sigma_i$ being elements of ${\mth{Z}}/2{\mth{Z}}$ that we will denote by $+$ or $-$ signs in the sequel, let $Ch(h)=Ch(\sigma_1,\dots,\sigma_r)$ be the $L_{F,der}$-conjugacy class of chambers of $X_E$ containing the $C(\lambda_1,\dots,\lambda_r)$ such that for every $i$, $\lambda_i$ is a square (resp. not a square) if $\sigma_i=+$ (resp. $\sigma_i=-$). Of course the $Ch(h)$ depend on the choices we have made for $D$ and $\Sigma_a$.

We denote by $(e_1,\dots,e_r)$ the canonical basis of $H^1(\Gamma,K_{T,E}\cap L_{E,der})$ viewed as a ${\mth{Z}}/2{\mth{Z}}$-vector space. More precisely, for every $i\in\{1,\dots,r\}$, $e_i$ is the element $(+,\dots,+,-,+,\dots,+)$, where the minus sign is in $i$-th position, and corresponds by the above correspondence to the root $\beta_i\in\Sigma_a$ (or in other words, $(\sigma_1,\dots,\sigma_r)\in H^1(\Gamma,K_{T,E}\cap L_{E,der})$ corresponds to elements of $Ch_E$ of the form $u_{\beta_i}(\varpi_E^{2f_{D,E}(\beta_i)}\lambda_i)C_0$, where for every $i$, $\lambda_i$ is a square if and only if $\sigma_i=+$).

By a slight abuse of notation, for every $h=(\sigma_1,\dots,\sigma_r)\in H^1(\Gamma,K_{T,E}\cap L_{E,der})$ and every $f\in\mathcal{H}(X_E)^{G_{F,der}}$, we write $f(h)=f(\sigma_1,\dots,\sigma_r)$ for the constant value of $f$ on $Ch(\sigma_1,\dots,\sigma_r)$.

In the whole beginning of this section,  $D$ and $C_0$ have been chosen arbitrarily among the ones satisfying the required conditions. (We did not impose explicitely any particular conditions on $\Sigma_a$ either, but we of course still assume $\Sigma_a$  is the one given by either proposition \ref{sigmaa} or proposition \ref{sigmaa2} depending on the case.) Now it is time to make more precise choices. Let then $D$ be such that $\Phi_D$ is a standard Levi subsystem of $\Phi$; every element of $\Phi_D$ is then a sum of simple roots contained in $\Phi_D$. Let $C_0$ be the chamber of $A_{0,E}$ corresponding to the following concave function: for every $\alpha\in\Phi^+$, define $h(\alpha)$ the following way:
\begin{itemize}
\item if $\Sigma_a$ contains roots of every length, then $h(\alpha)$ is the number of simple roots (counted with multiplicities) $\alpha$ is the sum of;
\item if $\Phi$ is not simply-laced and $\Sigma_a$ contains only long roots, $h(\alpha)$ is the number of long roots (again, counted with multiplicities) among the simple roots $\alpha$ is the sum of.
\end{itemize}

Note that we see from proposition \ref{sigmaa} that the case where $\Phi$ is not simply-laced and $\Sigma_a$ contains only short roots cannot happen.

Set $f(\alpha)=-\frac{h(\alpha)}2$. Set also $f(-\alpha)=\frac{h(\alpha)+1}2$. It is easy to check that $f$ is concave; details are left to the reader. Moreover, since $f$ is concave and $f(\alpha)+f(-\alpha)=\frac 12$ for every $\alpha$, $f$ is the concave function $f_{C_0}$ associated to some chamber $C_0$ of $A_{0,E}$. We can also easily check that the extended set of simple roots associated to $C_0$ is $\Delta\cup\{-\alpha_0\}$.
(Note that $R(C_0)$ is not contained in $R(C_{0,F})$ in general.)

For every $\alpha<0$ which is the inverse of the sum of an odd number of simple roots in $\Phi^+$, $f_{C_0}(\alpha)$ is an integer, hence when $\Sigma_a=\{\beta_1,\dots,\beta_r\}$ contains roots of every length, we see with the help of proposition \ref{odd} that $f_{C_0}(\beta_i)$ is an integer for every $i$. Now we check that it is also true when $\Phi$ is not simply-laced and $\Sigma$ contains only long roots. In that case, the assertion is an immediate consequence of proposition \ref{odd} and the following lemma:

\begin{lemme}\label{shoev}
Assume $\Phi$ is of type $B_d$, $C_d$ or $F_4$. Let $\alpha$ be a positive long root, and write $\alpha=\sum_{i=1}^d\lambda_i\alpha_i$, with $\alpha_1,\dots,\alpha_d$ being the elements of $\Delta$. Then for every $i$ such that $\alpha_i$ is short, $\lambda_i$ is even.
\end{lemme}

We prove the result by induction on $h(\alpha)$. If $h(\alpha)=1$, then $\alpha$ is a long simple root and the result is trivial. Assume $h(\alpha)>1$ and let $i$ be such that $\alpha-\alpha_i$ is a root. If $\alpha_i$ is long, then $\alpha-\alpha_i$ is also long and positive and $h(\alpha-\alpha_i)=h(\alpha)-1$; the result then follows from the induction hypothesis. Assme now $\alpha_i$ is short. Then $\alpha$ and $\alpha_i$ generate a subsystem of type $B_2$ of $\Phi$, which implies in particular, since $\alpha_i$ is a simple root and $\alpha\neq\alpha_i$, that $\alpha-2\alpha_i$ is also a positive root and is long. The result then follows from the induction hypothesis applied to $\alpha-2\alpha_i$. $\Box$

Now we prove that for every $f\in\mathcal{H}(X_E)^{G_{F,der}}$, the $f(h)$, $h\in H^1(\Gamma,K_{T,E}\cap L_{E,der})$, are all determined by $f(1)$.  We then establish relations between the $f(h)$ using the $G_{F,der}$-invariance of $f$ and the following two lemmas:

\begin{lemme}\label{sroot}
Let $i$ be an element of $\{1,\dots,r\}$; assume $\beta_i$ is the negative of a simple root in $\Phi^+$. Then for every $h\in H^1(\Gamma,K_{T,E}\cap L_{E,der})$ and every $f\in\mathcal{H}(X_E)^{G_{F,der}}$, $f(e_ih)=-f(h)$.
\end{lemme}

Let $C=C(\lambda_1,\dots,\lambda_r)$ be an element of $Ch(h)$. Set ${\mth{G}}_D=K_{D,E}/K_{D,E}^0$ and let ${\mth{P}}_i$ be the parabolic subgroup of ${\mth{G}}_D$ generated by ${\mth{B}}_0$ and the root subgroup $U_{\beta_i}$. Let $K_i\subset K_D$ be the corresponding parahoric subgroup of $G_E$ and let $D'_i$ be the codimension $1$ facet of $X_E$ associated to $K_i$. The chambers of $X_E$ admitting $D'_i$ as a wall are precisely the ones corresponding to the Iwahori subgroups contained in $K_i$. Out of these $q+1$ chambers, two do not belong to $Ch_{D,a,L}$ (the chamber $C="C(\lambda_1,\dots,\lambda_{i-1},0,\lambda_{i+1},\dots,\lambda_r)$" (with a slight abuse of notation)
and $C'=s_{\beta_i}(C)$), which implies that every element of $\mathcal{H}(X_E)^{G_{F,der}}$ is zero on them, and the remaining $q-1$ are the ones of the form $C(\lambda_1,\dots,\lambda_{i-1},\mu,\lambda_{i+1},\dots,\lambda_r)$ with $\mu\in k_F^*$; since exactly half of the elements of $k_F^*$ are squares, the lemma follows immediately from the harmonicity condition. $\Box$

\begin{lemme}\label{cedeux}
Let $\beta_i,\beta_j$ be two elements of $\Sigma_a$ satisfying the following conditions:
\begin{itemize}
\item $\alpha=\frac{\beta_j-\beta_i}2$ is a root, and $\beta_j$ and $\alpha$ generate a subsystem of $\Phi$ of type $B_2$;
\item $\alpha$ is the negative of a simple root of $\Phi^+$, and $f_{D,E}(\alpha)$ is an integer.
\end{itemize}
Then for every $h\in H^1(\Gamma,K_{T,E}\cap L_{E,der})$ and every $f\in\mathcal{H}(X_E)^{G_{F,der}}$, $f(e_ie_jh)=-f(h)$.
\end{lemme}

We first remark that by corollary \ref{spv2}, we have $f_{D,E}(\alpha)+f_{D,E}(-\alpha)=0$, hence if $f_{D,E}(\alpha)$ is an integer, $f_{D,E}(-\alpha)$ is an integer as well.

Set ${\mth{G}}_D=K_{D,E}/K_{D,E}^0$, let ${\mth{T}}_0$ be the image of $K_{T_0}$ in ${\mth{G}}_D$ and let ${\mth{B}}_0$ be the Borel subgroup of ${\mth{G}}_D$ containing ${\mth{T}}_0$ associated to $\Phi^+$. The root $-\alpha$, being a simple root in $\Phi^+$, is also a simple root in ${\mth{G}}_D$ in the set of positive roots associated to ${\mth{B}}_0$. Let ${\mth{P}}'$ be the parabolic subgroup of ${\mth{G}}_D$ generated by ${\mth{B}}_0$ and the root subgroup associated with $-\alpha$, and let $K$ and $D'$ be defined as the $K_i$ and $D_i$ of lemma \ref{sroot}, relatively to ${\mth{P}}'$ this time. The chambers of $X_E$ admitting $D'$ as a wall are the ones of the form:
\[C_l=(\prod_{i=1}^ru_{-\beta_i}(\lambda_i))lC_0,\]
where $l$ is an element of the Levi component ${\mth{M}}'$ of ${\mth{P}}'$, which is the product of a subgroup ${\mth{M}}''$ of type $A_1$ by the image of $K_{T_0}$ in ${\mth{G}}_D$; since $K_{T_0}$ stabilizes $C_0$ we can assume that $l\in{\mth{M}}''$, and to simplify notations we can consider $l$ as an element of $GL_2(k_F)$. On the other hand, $l$ admits representatives in $G_F$, hence:
\[f(C_l)=f(l^{-1}C_l)=f(l^{-1}(\prod_{i=1}^ru_{-\beta_i}(\lambda_i))lC_0,).\]
Since conjugating $\prod_{i=1}^ru_{-\beta_i}(\lambda_i)$ by $l$ leaves every term of the product but the $i$-th and $j$-th unchanged, we are reduced to the case where $d=2$ and $\Phi$ itself is of type $B_2$, in which case $G$ is the group $SO_5=PGSp_4$.

It turns out to be more convenient to work with $G=GSp_4$. The harmonicity condition applied to the chambers containing $D'$ can then, up to $G_F$-conjugation of the involved chambers, be rewritten as follows, if $h=(\sigma_1,\sigma_2)$: 
\begin{equation}\label{bedeux}\sum_{l\in R}f\left(\left(\begin{array}{cc}Id&0\\{}^\tau l\left(\begin{array}{cc}0&\sigma_1\varpi_E\\\sigma_2\varpi_E&0\end{array}\right)l&Id\end{array}\right)C_0\right)=0,\end{equation}
where $R$ is a set of representatives of the right classes of $GL_2(k_F)$ modulo ${\mth{B}}_0$. We thus have to find a set $R$ such that for every $l\in R$, if $C'_l$ is the chamber defined in the above sum, either $C'_l$ belongs to $Ch(h')$ for some $h'\in H^1(\Gamma,K_{T,E}\cap L_{E,der})$ or $f(C'_l)=0$.

To simplify the notations, we only write down the proof of  the case $h=1$; the other cases can be treated in a similar way. If $l=\left(\begin{array}{cc}a&b\\c&d\end{array}\right)\in GL_2({\mth{F}}_q)$, then we have:
\[{}^\tau l\left(\begin{array}{cc}0&1\\1&0\end{array}\right)l=\left(\begin{array}{cc}ab+cd&a^2+c^2\\b^2+d^2&ab+cd\end{array}\right),\]
which means that we only have to consider the $C_l$ such that there exists $l'\in{\mth{B}}$ such that $ll'$ satisfies the condition $ab+cd=0$; since that conditon is obviously right ${\mth{T}}$-invariant we can even assume that $l'$ is unipotent. A simple computation shows that in this case, $a^2+c^2$ and $b^2+d^2$ are either both squares or both non-squares, which implies that $C_l$ belongs to either $Ch(1)$ or $Ch(e_1e_2)$.

Consider first the element $l_\infty=\left(\begin{array}{cc}0&1\\1&0\end{array}\right)$. This element satisfies $a^2+c^2=b^2+d^2=1$, hence we have $C_{l_\infty}\in Ch_x(1)$. Moreover, none of the $l_\infty u$, with $u\neq 0$ belonging to the unipotent radical ${\mth{U}}$ of ${\mth{B}}_0$, satisfies the condition $ab+cd=0$.

Consider now, for every $y\in k_F$, the element $l_y=\left(\begin{array}{cc}1&0\\y&1\end{array}\right)$ of $GL_2(k_F)$. Another simple computation shows that there exists an element of $l_y{\mth{U}}$ satisfying the condition $ab+cd=0$ if and only if $1+y^2\neq 0$, and that in that case, $\left(\begin{array}{cc}1&\frac{-y}{1+y^2}\\y&\frac 1{1+y^2}\end{array}\right)$ is the only such element. To prove the lemma, we now only have to compute the number of $y\in k$ such that $1+y^2=a^2+c^2$ is nonzero and a square (resp. not a square).

Assume there exists $e\in k_F^*$ such that $1+y^2=e^2$; we then have $(e+y)(e-y)=1$. Set $\lambda=e+y$; we then have $\lambda(\lambda-2y)=1$, hence $\lambda-\frac 1\lambda=2y$. Moreover, it is easy to check that $\lambda-\frac 1\lambda=\mu-\frac 1\mu$ if and only if either $\lambda=\mu$ or $\lambda=-\frac 1\mu$.

Assume first $-1$ is not a square in $k_F^*$. Then $1+y^2$ is always nonzero, and we only have to count the number of different $y$ such that $1+y^2$ is a square. On the other hand, we never have $\lambda=-\frac 1\lambda$, hence for every $y\in k_F$, there are always either $0$ or $2$ values of $\lambda$ such that $\lambda-\frac 1\lambda=2y$. Hence the number of possible values for $y$ is $\frac{q-1}2$, which proves that for a suitable choice of $R$, taking into account $l_\infty$, there are exactly $\frac{q+1}2$ terms in the sum such that $C_l\in Ch(1)$ (resp. $C_l\in Ch(e_1e_2)$). The lemma then follows immediately from the harmonicity condition.

Assume now $-1$ is a square in $k_F^*$. Then each one of its square roots $\lambda$ satisfies $\lambda=-\frac 1\lambda$ and is its own image by $\lambda\mapsto\frac 12(\lambda-\frac 1\lambda)$, hence by the previous remark is also its only inverse image by that same map. On the other hand, every $y$ such that $y^2\neq -1$ has either $0$ or $2$ inverse images, hence there are exactly $\frac{q+1}2$ elements $y$ such that $y^2+1$ is a square, $\frac{q-3}2$ of them not being roots of $-1$, and $\frac{q-1}2$ elements $y$ such that $y^2+1$ is not a square. Taking into account $l_\infty$ once again, we conclude as in the previous case. $\Box$

Now we use these lemmas to prove theorem \ref{th2}. We already know that every $f\in\mathcal{H}(X_E)^{G_{F,der}}$ is entirely determined by the $f(h)$, $h\in H^1(\Gamma,K_{T\cap L_{E,der}})$; it then only remains to prove the following proposition:

\begin{prop}\label{uniccha}
Let $f$ be any element of $\mathcal{H}(X_E)^{G_{F,der}}$, viewed as a function on $H^1(\Gamma,K_{T\cap L_{E,der}})$. Then $f$ is entirely determined by $f(1)$.
\end{prop}

If $\lambda_1,\dots,\lambda_r$ are elements of $k^*$ such that $C(\lambda_1,\dots,\lambda_r)\in G_F C(1,\dots,1)$, and if $h'$ is the element of $H^1(\Gamma,K_{T\cap L_{E,der}})$ corresponding to the elements $\lambda_1,\dots,\lambda_r$, then we have $f(h'h)=f(h)$ for every $h\in H^1(\Gamma,K_{T\cap L_{E,der}})$. Moreover, if $i$ is such that $\beta_i$ is the negative of a simple root, by lemma \ref{sroot}, setting $h'=e_i$, $f(h'h)=-f(h)$ for every $h\in H^1(\Gamma,K_{T\cap L_{E,der}})$. Finally, if $\beta_i,\beta_j$ are two elements of $\Sigma_a$ satisfying the conditions of lemma \ref{cedeux}, then by that lemma, setting $h'=e_ie_j$, we have $f(h'h)=-f(h)$ for every $h\in H^1(\Gamma,K_{T\cap L_{E,der}})$. We thus only have to prove that the set $S$ of all these various elements $h'$ always generates $H^1(\Gamma,K_{T\cap L_{E,der}})$ as a ${\mth{Z}}/2{\mth{Z}}$-vector space.

We proceed by a case-by-case analysis. In the rest of the proof, the $\alpha_i$ and the $\varepsilon_i$ are defined the same way as in \cite[plates I to IX]{bou}.
.
\begin{itemize}
\item Assume first $\Phi$ is of type $A_d$, with $d=2n-1$ being odd; by proposition \ref{sigmaa2}, $\Phi_D$ is then the Levi subsystem of $\Phi$ generated by the simple roots $\alpha_{2i-1}$, $i=1,\dots,n$, and we can set for every $i$ $\beta_i=-\alpha_{2i-1}$, which is always the negative of a simple root of $\Phi^+$; by lemma \ref{sroot}, for every $i\in\{1,\dots,n\}$, $e_i\in S$ and $f(e_i)=-f(1)$ for every $f\in\mathcal{H}(X_E)^{G_{F,der}}$. The result follows.

\item Assume now $\Phi$ is of type $B_d$, with $d=2n$ being even; we have $\Phi_D=\Phi$. By proposition \ref{sigmaa}, for every $i\in\{1,\dots,n\}$, we can set $\beta_{2i-1}=-\varepsilon_{2i-1}-\varepsilon_{2i}$ and $\beta_{2i}=-\varepsilon_{2i-1}+\varepsilon_{2i}$. The $\beta_{2i}$ are then negatives of simple roots of $\Phi^+$, hence by lemma \ref{sroot}, for every $i$, $e_{2i}\in S$ and $f(e_{2i})=-f(1)$ for every $f\in\mathcal{H}(X_E)^{G_{F,der}}$. Moreover, for every element of $\Phi$ of the form $\alpha=\varepsilon_{2i}+\varepsilon_{2i+1}$, it is easy to check that $<\beta_j,\alpha^\vee>$ is odd if and only if $j\in\{2i-1,2i,2i+1, 2i+2\}$, hence if $c$ is an element of $\mathcal{O}_F^*$ which is not a square, $\alpha^\vee(c)$ acts on $H^1(\Gamma,K_{T\cap L_{E,der}})$ by multiplication by $e_{2i+1}e_{2i}e_{2i+1}e_{2i+2}$, which implies that $e_{2i-1}e_{2i}e_{2i+1}e_{2i+2}\in S$ and $f(e_{2i-1}e_{2i}e_{2i+1}e_{2i+2})=f(1)$ for every $f\in\mathcal{H}(X_E)^{G_{F,der}}$. We thus have obtained $2n-1$ linearly independent elements of $S$; we still need one more.

We will now prove that $e_{2n+1}e_{2n}\in S$ and $f(e_{2n+1}e_{2n})=-f(1)$ for every $f\in\mathcal{H}(X_E)^{G_{F,der}}$. Let $\alpha=\varepsilon_n$ be the unique short simple root in $\Phi^+$; the roots $\beta_{2n-1}$, $\beta_{2n}$ and $\alpha$ then satisfy the conditions of lemma \ref{cedeux}, and the desired result follows.

\item Assume now $\Phi$ is of type $B_d$, with $d=2n+1$ being odd; we have $\Phi_D=\Phi$. By proposition \ref{sigmaa}, we can define the $\beta_i$, $i\leq 2n$, as in the previous case and set $\beta_d=-\varepsilon_d$. Then for $j$ being either an even integer or $d$, $\beta_j$ is the negative of a simple root, hence by lemma \ref{sroot} $e_j\in S$ and $f(e_j)=-f(1)$ for every $f\in\mathcal{H}(X_E)^{G_{F,der}}$; moreover, for every $i\in\{1,\dots,n-1\}$, we obtain $e_{2i-1}e_{2i}e_{2i+1}e_{2i+2}\in S$ and $f(e_{2i-1}e_{2i}e_{2i+1}e_{2i+2})=f(1)$ for every $f\in\mathcal{H}(X_E)^{G_{F,der}}$ by the same reasoning as in the previous case; we also similarly obtain $e_{d-2}e_{d-1}e_d\in S$ and $f(e_{d-2}e_{d-1}e_d)=f(1)$ for every $f\in\mathcal{H}(X_E)^{G_{F,der}}$. This makes $2n+1$ linearly independent elements of $S$, as desired.

\item Assume now $\Phi$ is of type $C_d$; we have $\Phi_D=\Phi$. By proposition \ref{sigmaa}, we can set $\beta_i=-2\varepsilon_i$ for every $i$. The root $\beta_d$ is the negative of a simple root, hence $e_d\in S$ and $f(e_d)=-f(1)$ for every $f\in\mathcal{H}(X_E)^{G_{F,der}}$; moreover, for every $i\in\{1,\dots,d-1\}$, $\alpha_i=\varepsilon_i-\varepsilon_{i+1}$ is a simple root and $\beta_i$, $\beta_{i+1}$ and $\alpha_i$ satisfy the conditions of lemma \ref{cedeux}, hence $e_ie_{i+1}\in S$ and $f(e_ie_{i+1})=-f(1)$ for every $f\in\mathcal{H}(X_E)^{G_{F,der}}$. We thus obtain $d$ linearly independent elements of $S$, as desired.

\item Assume now $\Phi$ is of type $D_d$, with $d=2n$ being even; we have $\Phi_D=\Phi$. By proposition \ref{sigmaa}, we can choose the $\beta_i$ the same way as in the case $B_{2n}$, and it is easy to check that the first $2n-1$ linearly independent elements of $S$ are the same, with the same relative values of $f\in\mathcal{H}(X_E)^{G_{F,der}}$; to get one more, we simply remark that $-\beta_{2n-1}$ is now also the negative of a simple root of $\Phi^+$, which implies that $e_{2n-1}\in S$ and $f(e_{2n-1})=-f(1)$.

\item Assume now $\Phi$ is of type $D_d$, with $d=2n+1$ being odd; we deduce from proposition \ref{sigmaa2} that $\Phi_D$ is then the Levi subsystem of $\Phi$ generated by the simple roots $\alpha_i=\varepsilon_i-\varepsilon_{i+1}$, $i=2,\dots,d-1$, and $\alpha_d=\varepsilon_{d-1}+\varepsilon_d$, and $\Sigma_a$ is of cardinality $2n$. By that same proposition, the $\beta_i$ are defined the same way as in the cases $B_{2n}$ and $D_{2n}$, except that we add $1$ to every index of the $\varepsilon_i$ (i.e. $\varepsilon_i$ becomes $\varepsilon_{i+1}$): more precisely, we now have $\beta_{2i-1}=-\varepsilon_{2i}-\varepsilon_{2i+1}$ and $\beta_{2i}=-\varepsilon_{2i}+\varepsilon_{2i+1}$ for every $i\in\{1,\dots,n\}$. The $2n$ linearly independent elements of $S$ and the relative values of $f\in\mathcal{H}(X_E)^{G_{F,der}}$ are obtained as in the case $D_{2n}$, taking into account the shift of indices.

\item Assume now $\Phi$ is of type $E_6$; by proposition \ref{sigmaa2}, $\Phi_D$ is then the Levi subsystem of $\Phi$ generated by $\alpha_2,\dots,\alpha_5$, and $\Sigma_a$ is of cardinality $4$. By that same proposition, we can set $\beta_1=-\alpha_2-\alpha_3-2\alpha_4-\alpha_5$, $\beta_2=-\alpha_2$, $\beta_3=-\alpha_3$ and $\beta_4=-\alpha_5$. Then $\beta_2$, $\beta_3$ and $\beta_4$ are negatives of simple roots, hence for every $i\in\{2,3,4\}$, $e_i\in S$ and $f(e_i)=-f(1)$ for every $f\in\mathcal{H}(X_E)^{G_{F,der}}$. Moreover, it is easy to check that $<\beta_i,\alpha_4^\vee>$ is odd for every $i$, hence if $c$ is an element of $\mathcal{O}_E^*$, which is not a square, we have the following coroot action on $H^1(\Gamma,K_{T\cap L_{E,der}})$:
\[\alpha_4^\vee(c)h=e_1e_2e_3e_4h.\]
Hence $e_1e_2e_3e_4\in S$ and $f(e_1e_2e_3e_4)=f(1)$ for every $f\in\mathcal{H}(X_E)^{G_{F,der}}$. This makes $4$ linearly independent elements of $S$, as desired.

\item Assume now $\Phi$ is of type $E_7$; we have $\Phi_D=\Phi$. By proposition \ref{sigmaa}, we can set $\beta_1=-\alpha_0$, $\beta_2=-\alpha_2$, $\beta_3=-\alpha_3$, $\beta_4=-\alpha_2-\alpha_3-2\alpha_4-2\alpha_5-2\alpha_6-\alpha_7$, $\beta_5=-\alpha_5$, $\beta_6=-\alpha_2-\alpha_3-2\alpha_4-\alpha_5$ and $\beta_7=-\alpha_7$. For every $i\in\{2,3,5,7\}$, $\beta_i$ is the negative of a simple root, hence by lemma \ref{sroot} $e_i\in S$ and $f(e_i)=-f(1)$ for every $f\in\mathcal{H}(X_E)^{G_{F,der}}$; on the other hand, we have:
\begin{itemize}
\item $<\beta_i,\alpha_1^\vee>$ is odd if and only if $i=1,3,4,6$;
\item $<\beta_i,\alpha_4^\vee>$ is odd if and only if $i=2,3,5,6$;
\item $<\beta_i,\alpha_6^\vee>$ is odd if and only if $i=4,5,7$;
\end{itemize}
hence if $c$ is an element of $\mathcal{O}_E^*$ which is not a square, we have the following coroot actions on $H^1(\Gamma,K_{T\cap L_{E,der}})$:
\begin{itemize}
\item $\alpha_1^\vee(c)h=e_1e_3e_4e_6h$;
\item $\alpha_4^\vee(c)h=e_2e_3e_5e_6h$;
\item $\alpha_6^\vee(c)h=e_4e_5e_7h$,
\end{itemize}
Hence $e_1e_3e_4e_6$, $e_2e_3e_5e_6$ and $e_4e_5e_7$ belong to $S$; for every $f\in\mathcal{H}(X_E)^{G_{F,der}}$, the value of $f$ on them is then equal to $f(1)$. We thus obtain $7$ linearly independent elements of $S$, as desired.

\item Assume now $\Phi$ is of type $E_8$; we have $\Phi_D=\Phi$. By proposition \ref{sigmaa}, we can set $\beta_1=-\alpha_0$, $\beta_2=-\alpha_2$, $\beta_3=-\alpha_3$, $\beta_4=-2\alpha_1-2\alpha_2-3\alpha_3-4\alpha_4-3\alpha_5-2\alpha_6-\alpha_7$, $\beta_5=-\alpha_5$, $\beta_6=-\alpha_2-\alpha_3-2\alpha_4-2\alpha_5-2\alpha_6-\alpha_7$, $\beta_7=-\alpha_7$ and $\beta_8=-\alpha_2-\alpha_3-2\alpha_4-\alpha_5$. For every $i\in\{2,3,5,7\}$, as in the case $E_7$, $\beta_i$ is the negative of a simple root, hence by lemma \ref{sroot} $e_i\in S$ and $f(e_i)=-f(1)$ for every $i$; on the other hand, we have:
\begin{itemize}
\item $<\beta_i,\alpha_1^\vee>$ is odd if and only if $i=3,4,6,8$;
\item $<\beta_i,\alpha_4^\vee>$ is odd if and only if $i=2,3,5,8$;
\item $<\beta_i,\alpha_6^\vee>$ is odd if and only if $i=5,6,7$;
\item $<\beta_i,\alpha_8^\vee>$ is odd if and only if $i=1,4,6,7$;
\end{itemize}
hence if $c$ is an element of $\mathcal{O}_E^*$ which is not a square, we have the following coroot actions on $H^1(\Gamma,K_{T\cap L_{E,der}})$:
\begin{itemize}
\item $\alpha_1^\vee(c)h=e_3e_4e_6e_8h$;
\item $\alpha_4^\vee(c)h=e_2e_3e_5e_8h$;
\item $\alpha_6^\vee(c)h=e_5e_6e_7h$;
\item $\alpha_8^\vee(c)h=e_1e_4e_6e_7$.
\end{itemize}
Hence $e_2e_3e_4e_6$, $e_4e_5e_6e_7$, $e_3e_7e_8$ and $e_1e_2e_3e_8$ belong to $S$ and for every $f\in\mathcal{H}(X_E)^{G_{F,der}}$, the value of $f$ on them is equal to $f(1)$. We thus obtain $8$ linearly independent elements of $S$, as desired.

\item Assume now $\Phi$ is of type $F_4$; we have $\Phi_D=\Phi$. By proposition \ref{sigmaa}, we can set $\beta_1=-\alpha_0$, $\beta_2=-\alpha_2$, $\beta_3=-\alpha_2-2\alpha_3$ and $\beta_4=-\alpha_2-2\alpha_3-2\alpha_4$. Since $(\alpha_2,\alpha_3,\alpha_4)$ is the set of simple roots of a standard Levi subsystem of type $C_3$ of $\Phi$, with the help of the case $C_d$ applied to that subsystem, we obtain that $e_2$, $e_3$ and $e_4$ belong to $S$ and $f(e_4)=-f(e_3)=f(e_2)=-f(1)$ for every $f\in\mathcal{H}(X_E)^{G_{F,der}}$; on the other hand, $<\beta_i,\alpha_1^\vee>$ is odd for every $i$, hence if $c$ is an element of $\mathcal{O}_E^*$ which is not a square, we have the following coroot action on $H^1(\Gamma,K_{T\cap L_{E,der}})$:
\[\alpha_1^\vee(c)h=e_1e_2e_3e_4h.\]
Hence $e_1e_2e_3e_4$ belongs to $S$ as well, and $f(e_1e_2e_3e_4)=f(1)$ for every $f\in\mathcal{H}(X_E)^{G_{F,der}}$. The result follows.

\item Assume finally $\Phi$ is of type $G_2$; we have $\Phi_D=\Phi$. By proposition \ref{sigmaa}, we can set $\beta_1=-\alpha_1$ and $\beta_2=-\alpha_0$; $\beta_1$ is then the negative of a simple root, hence by lemma \ref{sroot} $e_1\in S$ and $f(e_1)=-f(1)$ for every $f\in\mathcal{H}(X_E)^{G_{F,der}}$; on the other hand, $<-\alpha_0,\alpha_2^\vee>$ and $<-\alpha_1,\alpha_2^\vee>$ are both odd, hence if $c$ is an element of $\mathcal{O}_E^*$ which is not a square, we have the following coroot action on $H^1(\Gamma,K_{T\cap L_{E,der}})$:
\[\alpha_2^\vee(c)h=e_1e_2h.\]
Hence $e_1e_2$ belongs to $S$ and $f(e_1e_2)=f(1)$ for every $f\in\mathcal{H}(X_E)^{G_{F,der}}$. The result follows.
\end{itemize}
The proposition is now proved. $\Box$

\begin{cor}
Assume $\Phi$ is not of type $A_{2n}$ for any $n$. Then theorem \ref{th2} holds.
\end{cor}

\subsection{Action of some elements of $G_F$}

We finish this section by summarizing the action of the simple coroots of $\Phi^+$ on $H^1(\Gamma,K_{T\cap L_{E,der}})$ (proposition \ref{sract}) and the elements of the canonical basis of $H^1(\Gamma,K_{T\cap L_{E,der}})$ on $\mathcal{H}(X_E)^{G_{F,der}}$ (proposition \ref{eic}):

\begin{prop}\label{sract}
Assume $\Phi$ is not of type $A_{2n}$ for any $n$; let $h$ be an element of $H^1(\Gamma,K_{T\cap L_{E,der}})$, and let $c$ be an element of $k_F^*$ which is not a square. We have:
\begin{itemize}
\item if $\Phi$ is of type $A_{2n-1}$:
\begin{itemize}
\item $\alpha_{2i+1}^\vee(c)h=h$ for every $i$;
\item $\alpha_{2i}^\vee(c)h=e_ie_{i+1}h$ for every $i$.
\end{itemize}
\item if $\Phi$ is of type $B_d$:
\begin{itemize}
\item $\alpha_i^\vee(c)h=h$ if either $i$ is odd or $i=d$;
\item $\alpha_i^\vee(c)h=e_{i-1}e_ie_{i+1}e_{i+2}h$ if $i$ is even and $<d-1$;
\item $\alpha_{d-1}^\vee(c)h=e_{d-2}e_{d-1}e_d$ if $d$ is odd.
\end{itemize}
\item if $\Phi$ is of type $C_n$, $\alpha_i^\vee(c)h=h$ for every $i$;
\item if $\Phi$ is of type $D_n$:
\begin{itemize}
\item $\alpha_{d-i}^\vee(c)h=h$ for every odd $i$;
\item $\alpha_{d-i}^\vee(c)h=e_{d-i-1}e_{d-i}e_{d-i+1}e_{d-i+2}h$ for every $i$ even, positive and such that $d-i>1$;
\item $\alpha_d^\vee(h)=h$;
\item when $d$ is odd, $\alpha_1^\vee(c)h=e_1e_2h$.
\end{itemize}
\item if $\Phi$ is of type $E_6$:
\begin{itemize}
\item $\alpha_i^\vee(c)h=h$ for every $i<4$;
\item $\alpha_4^\vee(c)h=e_1e_2e_3e_4h$;
\end{itemize}
\item if $\Phi$ is of type $E_7$:
\begin{itemize}
\item $\alpha_i^\vee(c)h=h$ for $i=2,3,5,7$;
\item $\alpha_1^\vee(c)h=e_1e_3e_4e_6h$;
\item $\alpha_4^\vee(c)h=e_2e_3e_5e_6h$;
\item $\alpha_6^\vee(c)h=e_4e_5e_7h$,
\end{itemize}
\item if $\Phi$ is of type $E_8$:
\begin{itemize}
\item $\alpha_i^\vee(c)h=h$ for $i=2,3,5,7$;
\item $\alpha_1^\vee(c)h=e_3e_4e_6e_8h$;
\item $\alpha_4^\vee(c)h=e_2e_3e_5e_8h$;
\item $\alpha_6^\vee(c)h=e_5e_6e_7h$;
\item $\alpha_8^\vee(c)h=e_1e_4e_6e_7$.
\end{itemize}
\item if $\Phi$ is of type $F_4$:
\begin{itemize}
\item $\alpha_1^\vee(c)h=e_1e_2e_3e_4h;$
\item $\alpha_i^\vee(c)h=h$ for every $i\geq 2$;
\end{itemize}
\item if $\Phi$ is of type $G_2$, $\alpha_1^\vee(c)h=h$ and $\alpha_2^\vee(c)h=e_1e_2h$.
\end{itemize}
\end{prop}

\begin{prop}\label{eic}
Assume $\Phi$ is not of type $A_{2n}$ for any $n$. Let $f$ be a nonzero element of $\mathcal{H}(X_E)^{G_{F,der}}$; we have:
\begin{itemize}
\item if $\Phi$ is of type $A_{2n-1}$, $f(e_i)=-f(1)$ for every $i$; 
\item if $\Phi$ is of type $B_d$, $f(e_i)=-f(1)$ if either $i=d$ or $i$ is even, and $f(e_i)=f(1)$ if $i$ is odd and $<d$; 
\item if $\Phi$ is of type $C_d$, $f(e_i)=(-1)^{d+1-i}f(1)$ for every $i$;
\item if $\Phi$ is of type $D_d$ (either odd or even), $f(e_i)=-f(1)$ for every $i$;
\item if $\Phi$ is of type $E_6$, $f(e_i)=-f(1)$ for every $i$;
\item if $\Phi$ is of type $E_7$, $f(e_i)=f(1)$ if $i$ is either $1$ or $2$ and $-f(1)$ if $i\geq 3$;
\item if $\Phi$ is of type $E_8$, $f(e_i)=f(1)$ if $i$ is either $1$ or $3$ and $-f(1)$ in the other cases;
\item if $\Phi$ is of type $F_4$, $f(e_i)=-f(1)$ for every $i$;
\item if $\Phi$ is of type $G_2$, $f(e_i)=-f(1)$ for every $i$.
\end{itemize}
\end{prop}

These relations either are already contained in the proof of proposition \ref{uniccha} or can be deduced from the relations established during that proof by easy computations. Details are left to the reader. $\Box$

\section{Proof of the $\chi$-distinction}

\subsection{A convergence result}

Now we go to the proof of theorem \ref{th1}. Before defining our linear form $\lambda$, we have to prove a preliminary result, which plays here the same role as \cite[lemma 4.5]{bc} for the unramified case, except that it now works for any value of $q$ thanks to the use of the Poincar\'e series.

To make notations clearer, we denote by $d_E(.,.)$ (resp. $d_F(.,.)$) the combinatorial distance between two chambers of $X_E$ (resp. $X_F$).
 
\begin{prop}\label{elun}
Let $f$ be an element of $\mathcal{H}(X_E)^\infty$, and let $O$ be any $G_F$-orbit of chambers of $X_E$. Then we have:
\[\sum_{C\in O}|f(C)|<+\infty.\]
\end{prop}

Fix an element $C$ of $O$. Let $C_0$ be an element of $Ch_E$ whose geometric realization is contained in $\mathcal{B}_F$ and such that $d_E(C,C_0)$ is minimal, and let $C_F$ be the chamber of $X_F$ whose geometric realization contains $R(C_0)$. We first prove the following lemmas:

\begin{lemme}\label{dbdist}
Let $A_F$ be an apartment of $X_F$ containing $C_F$ and let $T$ be the associated $F$-split torus of $G$. For every $t\in T_F$, we have $d_E(C_0,tC_0)=2d_F(C_F,tC_F)$.
\end{lemme}

By eventually conjugating $C_0$ and $A_F$ by the same element of $G_F$ we may assume that $A_F=A_{0,F}$. Let $f_0,f_t,f_F,f_{t,F}$ be the concave functions associated respectively to $C_0$, $tC_0$, $C_F$ and $tC_F$; we have:
\[d_E(C_0,tC_0)=2\sum_{\alpha\in\Phi^+}|f_t(\alpha)-f_0(\alpha)|;\]
\[d_F(C_F,tC_F)=\sum_{\alpha\in\Phi^+}|f_{t,F}(\alpha)-f_F(\alpha)|.\]
On the other hand, since $t\in T_F$, for every $\alpha$, $f_t(\alpha)-f_0(\alpha)$ is an integer, and we deduce from this that $f_{t,F}(\alpha)-f_F(\alpha)=f_t(\alpha)-f_0(\alpha)$. The result follows immediately. $\Box$

\begin{lemme}\label{dbdist2}
There exists an integer $N_0$ such that for every $g\in G_F$, we have $d_E(C,gC)\geq 2d_F(C_F,gC_F)-N_0$.
\end{lemme}

Let $g$ be an element of $G_F$, let $\mathcal{A}$ be an apartment of $\mathcal{B}_F$ containing both $R(C_F)$ and $R(gC_F)$, let $T$ be the corresponding maximal $F$-split torus of $G_F$ and let $N_G(T)$ be the normalizer of $T$ in $G$; we have $gC_F=nC_F$ for some element $n$ of $N_G(T)_F$, hence $g$ is of the form $nh$, with $h\in K_{C_F,F}$.

Let $x$ be a special vertex of $C_F$; we can write $n=tn_0$, with $t\in T$ and $n_0\in K_{x,F}$. Set $C'=gC$, $C''=n_0hC$ and $C''_F=n_0C_F$; $C''_F$ also admits $x$ as a vertex. Since $n_0h$ always belongs to the open compact subgroup $K_{x,F}$ of $G_F$, the union of the $n_0hC$ (resp. of the $n_0C_F$) is bounded, which implies that there exists an integer $N''$ (resp. $N''_F$) such that we always have $d_E(C,C'')\leq N''$ (resp. $d_F(C_F,C''_F)\leq N''_F$). Moreover, according to lemma \ref{dbdist}, setting $C''_0=n_0hC_0$ and $C'_0=tC''_0=gC_0$, we have $d_E(C''_0,C'_0)=2d_F(C''_F,C'_F)$. On the other hand, since $C''=n_0hC$ and $C''_0=n_0hC_0$, we have $d_E(C'',C''_0)=d_E(C,C_0)$; similarly, $d_E(C',C'_0)=d_E(C,C_0)$. We finally obtain:
\[d_E(C,C')\geq d_E(C',C'')-d_E(C'',C)\]
\[\geq d_E(C'_0,C''_0)-d_E(C',C'_0)-d_E(C'',C''_0)-d_E(C'',C)\]
\[\geq 2d_F(C'_F,C''_F)-2d_E(C,C_0)-N''\]
\[\geq 2d_F(C_F,C'_F)-2d_F(C_F,C''_F)-2d_E(C,C_0)-N''\]
\[\geq 2d_F(C_F,C'_F)-2N''_F-2d_E(C,C_0)-N''.\]
We thus can set $N_0=2N''_F+2d_E(C,C_0)+N''$; the lemma is now proved. $\Box$

We can now prove the proposition.
We can write:
\[\sum_{C\in O}|f(C)|=\frac 1{[K_{C,F}:K_{C_F,F}\cap K_{C,F}]}\sum_{g\in G_F/K_{C,F}}|f(gC)|.\]
It is easy to check by induction that the number of chambers $C''_F$ of $X_F$ whose retraction on $A_{0,F}$ relatively to $C_F$ is some given chamber $C_{1,F}$ is $q^{d_F(C_F,C_{1,F})}$. By lemma \ref{dbdist2}, we obtain, $W'$ being the affine Weyl group of $G$ relative to $T_0$:
\[\sum_{g\in G_F/K_{C,F}}|f(gC)|\leq\sum_{g\in G_F/K_{C_F,F}}\frac 1{q^{2d_F(C_F,gC_F)-N_0}}\]
\[=\sum_{w\in W'}\frac {q^{d_F(C_F,wC_F)}}{q^{2d_F(C_F,wC_F)-N_0}}.\]
\[=\sum_{w\in W'}\frac 1{q^{l(w)-N_0}}.\]
By \cite[section 3]{mcdo}, the above sum converges for every $q>1$. The result follows immediately. $\Box$

\subsection{The case $A_d$, $d$ even}

Now we prove theorem \ref{th1} when $\Phi$ is of type $A_d$, with $d=2n$ being even, and $q$ is large enough. First we have:

\begin{prop}
Assume $\Phi$ is of type $A_d$, $d$ even. Then the Prasad character $\chi$ of $F$ is trivial.
\end{prop}

Let $\rho$ be tha half-sum of the elements of $\Phi^+$. Write $\rho=\sum_{i=1}^d\lambda_i\alpha_i$, the $\alpha_i$ being the elements of $\Delta$; by \cite[plate I, (VII)]{bou}, we have:
\[\rho=\sum_{i=1}^d\frac{i(d+1-i)}2\alpha_i.\]
Hence $\lambda_i$ is an integer for every $i$, and the proposition follows immediately from \cite[lemma 3.1]{cou2}. $\Box$

Now define the set $Ch_c$ of chambers of $X_E$ as in corollary \ref{chsol1}; by proposition \ref{chsol}, $Ch_c$ is $G_F$=stable and $G_F$ acts transitively on it. Set:
\[\lambda:f\in\mathcal{H}(X_E)^\infty\longmapsto\sum_{C\in Ch_c}f(C).\]
The linear form $\lambda$ is well-defined by proposition \ref{elun}, and obviously $G_F$-invariant. We want to prove that it is not identically zero on $\mathcal{H}(X_E)^\infty$.

By a slight abuse of notation, for every $C,C'\in Ch_c$, we write $d_F(C,C')$ for the combinatorial distance between the chambers of $X_F$ whose geometric realizations contain respectively $R(C)$ and $R(C')$.

Let $C$ be any element of $Ch_E$, and let $I$ be the Iwahori subgroup of $G_E$ fixing $C$. A well-known result about the Steinberg representation (see \cite{sha} for example) says that there exists a unique (up to a multiplicative constant) $I$-invariant element in the space of $St_E$, hence also in $\mathcal{H}(X_E)^\infty$. More precisely, set:
\[\phi_C:C'\in Ch_E\longmapsto (-q)^{-d_E(C,C')}.\]
It is easy to check that $\phi_C$ is $I$-invariant and satisfies the harmonicity condition. Hence every $I$-invariant element of $\mathcal{H}(X_E)^\infty$ is proportional to $\phi_C$; $\phi_C$ is called the (normalized) {\em Iwahori-spherical vector} of $\mathcal{H}(X_E)^\infty$ attached to $C$. Of course $\phi_C$ depends on $C$.

Now we prove the following proposition, from which theorem \ref{th1} follows immediately when $G$ is of type $A_{2n}$:

\begin{prop}\label{tva2n}
Let $C_0$ be any element of $Ch_c$. Then $\phi_{C_0}$ is a test vector for $\lambda$. More precisely, we have $\lambda(\phi_{C_0})=1$.
\end{prop}

Let $C_{0,F}$ be the chamber of $X_F$ whose geometric realization contains $R(C_0)$, let $C'_{0,F}$ be any chamber of $X_F$ adjacent to $C_{0,F}$ and let $C'_0$ be the unique element of $Ch_c$ whose geometric realization is contained in $R(C'_{0,F})$.

Let $\mathcal{A}$ be an apartment of $\mathcal{B}_F$ containing both $R(C_{0,F})$ and $R(C'_{0,F})$. Then $\mathcal{A}$ also contains both $R(C)$ and $R(C')$, hence also every minimal gallery between them.

First we prove the following lemmas:

\begin{lemme}\label{dis3}
The combinatorial distance between $C_0$ and $C'_0$ is $3$.
\end{lemme}

Let $A_E$ be the apartment of $X_E$ whose geometric realization is $\mathcal{A}$ and let $C$ be a chamber of $A_E$ adjacent to $C_0$; since by definition of $Ch_c$ none of the walls of $R(C_0)$ is contained in a codimension $1$ facet of $\mathcal{B}_F$, $R(C)$  is also contained in $R(C_{0,F})$, and since $C'_0$ is not an element of $Ch_c$, at least one of its walls has its geometric realization contained in $R(D)$, where $D$ is a wall of $C_{0,F}$. On the other hand, since $G_F$ is of type $A_{2n}$, the group of isomorphisms of $\mathcal{B}_F$ which stabilize $C_{0,F}$ is of order $2n+1$ by \cite[plate I]{bou}, hence acts transitively on the set of its walls; we can then assume without loss of generality that $D_F$ is the wall between $C_{0,F}$ and $C'_{0,F}$. Let $C'$ be the chamber of $A_E$ which is separated from $C$ by some wall whose geometric realization is contained in $R(D_F)$; by symmetry, $C'$ is adjacent to $C'_0$, Hence $(C_0,C,C',C'_0)$ is a gallery of length $3$ between $C_0$ and $C'_0$. On the other hand, every gallery $(C_0=C,C_1,\dots,C_s=C')$ between $C$ and $C'$ contained in $A$ must contain two chambers $C_{i+1}$ and $C_i$ separated by  the hyperplane of $A_E$ whose geometric realization contains $R(D_F)$, and the geometric realization of their common wall is then contained in $\mathcal{B}_F$; since $C_0$ and $C'_0$ are both elements of $Ch_c$, they are both distinct from both $C_i$ and $C_{i+1}$, and the length of any gallery beween them is then at least $3$. The result follows. $\Box$

We deduce immediately from the lemma the following corollary:

\begin{cor}
Let $\mathcal{H}$ be the hyperplane of $\mathcal{A}$ containing $R(D_F)$. For every $C\in Ch_c$ whose geometric realization is contained in $\mathcal{A}$, $d_E(C'_0,C)-d_E(C_0,C)$ is contained in $\{-3,-1,1,3\}$, and is positive (resp. negative) if $R(C)$ is contained in the same half-apartment with respect to $\mathcal{H}$ as $R(C_0)$ (resp. $R(C'_0)$).
\end{cor}

Now we examine more closely the structure of the subcomplex $Ch_\emptyset$.

\begin{lemme}\label{tcc0}
There are exactly two chambers of $A_E$ adjacent to $C_0$ and such that the geometric realization of one of their walls is contained in $\mathcal{H}$.
\end{lemme}

Let $H_E$ be the hyperplane of $A_E$ whose geometric realization is $\mathcal{H}$. We already know that there exists at least one chamber satisfying these conditions, namely the chamber $C$ of the gallery of length $3$ between $C_0$ and $C'_0$ defined during the proof of lemma \ref{dis3}. Since every such chamber contains a wall of $C_0$, its intersection with $C_0$ contains a facet $D$ of $H_E$ of codimension at most $2$, and in fact of codimension exactly $2$ since by hypothesis $H_E$ does not contain any wall of $C_0$. Since exactly two walls of $C_0$ contain $D$, there are also two chambers of $A_E$ adjacent to $C_0$ and containing $D$.

Let $C'$ be the unique chamber distinct from $C$ satisfying these conditions; we now only have to prove that one of the walls of $C'$ is contained in $H_E$. Let $K_D$ be the connected fixator of $D$, and let ${\mth{G}}_D$ be the quotient of $K_D$ by its pro-unipotent radical; ${\mth{G}}_D$ is then the group of $k_E$-points of a reductive group defined over $k_E$ whose root system is of rank $2$ and contained in a system of type $A_{2n}$, hence of type either $A_1^2$ or $A_2$, and the combinatorial distance between two chambers containing $D$ is equal to the combinatorial distance between the corresponding chambers in the spherical building of ${\mth{G}}_D$. If ${\mth{G}}_D$ is of type $A_1^2$, the combinatorial distance between $C_0$ and $C'_0$ can be at most $2$, which contradicts lemma \ref{dis3}; hence ${\mth{G}}_D$ must be of type $A_2$. Since the order of its Weyl group is then $6$, $K_D$ contains exactly $6$ Iwahori subgroups of $G_E$ containing the maximal compact subgroup $K_{T,E}$ of $T_E$, where $T$ is the maximal torus of $G$ associated to $A_E$, or equivalently, $D$ is contained in exactly $6$ chambers of $A_E$. Out of these six chambers, exactly four admit as a wall some facet of maximal dimension of any given hyperplane of $A_E$ containing $D$; this is in particular true for $H_E$. On the other hand, $C_0$ is one of these six chambers, and by symmetry $C'_0$ must be another one. Since none of these two admit any facet of maximal dimension of $H_E$ as a wall, then $C'$ must admit one and the lemma is proved. $\Box$

\begin{lemme}
Let $C$ be a chamber of $A_E$ adjacent to $C_0$. There are exactly two walls of $C$ whose geometric realizations are contained in walls of $R(C_F)$.
\end{lemme}

Let $C$, $C'$, $D$ and $H_E$ be defined as in the previous lemma. Since $C$ and $C'$ are both adjacent to $C_0$ and all three of them belong to $A_E$, $C$ and $C'$ cannot be adjacent to each other, hence their intersection is $D$, which proves that the walls of $C$ and $C'$ contained in $H_E$ are distinct. Hence by the previous lemma, the total number of walls of chambers of $A_E$ adjacent to $C_0$ whose geometric realizations are contained in the walls of $R(C_F)$ is $2(2n+1)$. On the other hand, as we have already seen, the group of automorphisms of $A_E$ stabilizing $C_0$ acts transitively on the set of its walls, hence also on the set of chambers of $X_E$ adjacent to $C_0$; since its action obviously preserves the number of walls of $C$ whose geometric realization is contained in walls of $R(C_F)$, that number must be two. $\Box$

Let $I_0$ be the Iwahori subgroup of $G_E$ fixing $C_0$; we have the following lemma:

\begin{lemme}\label{numch}
The number of elements of $Ch_c$ which are conjugated to $C$ by some element of $I_0$ is $q^{d_F(C_0,C)}$.
\end{lemme}

By \cite[lemma 4.2]{bc} and an obvious induction, it is enough to prove that two elements of $Ch_c$ are conjugated by an element of $I_0$ if and only if they are conjugated by an element of $I_{0,F}=I_0\cap G_F$. Let $C''$ be an element of $Ch_c$ conjugated to $C$ by some element of $I_0$, and let $C_F$ (resp. $C''_F$)) be the chamber of $X_F$ whose geometric realization contains $R(C)$ (resp. $R(C'')$). There exists then an element of $I_{0,F}=I_0\cap G_F$ sending $C_F$ to $C''_F$, and by unicity of the central chamber in the geometric realization of $C_F$ (resp. $C''_F$), that element must send $C$ on $C''$. The other implication being obvious, the lemma is proved. $\Box$

Now we prove proposition \ref{tva2n}. Let $C_0$ be the only element of $Ch_c$ whose geometric realization is contained in $R(C_{0,F})$, let $C$ be any element of $Ch_c$, set $d=d_E(C_0,C)$, and let $C_1$ be a chamber of $A_{0,E}$ adjacent to $C$ and such that $H_E$ contains a wall $D_1$ of $C_1$, First we assume that $C$ satisfies the following property:

{\bf (P1)}: There exists a minimal gallery of the form $(C_0,C_1,\dots,C_\delta=C)$,

and that $C_0$ and $C$ are in the same half-space of $A_{0,E}$ with respect to $H_E$. Let $C'_1$ be the other chamber of $\mathcal{B}_E$ admitting $D_1$ as a wall. Then $(C'_0,C'_1,C_1,\dots,C_\delta)$ is a minimal gallery of length $\delta+1$, from which we deduce by symmetry that if $C'$ is the image of $C$ by the orthogonal reflection with respect to $H_E$, $d(C_0,C')=\delta+1$. Hence we have $\phi_{C_0}(C)=(-q)^{-\delta}$ and $\phi_{C_0}(C')=(-q)^{-\delta-1}$.

On the other hand, by the same reasoning, if we set $\delta'=d_F(C_0,C')$, we have $d_F(C_0,C')=\delta'+1$.
From lemma \ref{numch}, we deduce that the sum of the $f(C'')$, when $C''$ runs through the set of conjugates of $C$ (resp. $C'$) by elements of $I_{0,F}$ is $\frac{q^{d_F(C_0,C)}}{(-q)^\delta}$ (resp. $\frac{q^{d_F(C_0,C')}}{(-q)^{\delta+1}}=\frac{q^{d_F(C_0,C)+1}}{(-q)^{\delta+1}}$). Since these two values are opposite to each other, their sum is zero. Since this is true for every $C$ satisfying {\bf (P1)} and on the same side of $H_E$ as $C_0$, we obtain the following lemma:

\begin{lemme}\label{zerosing}
The sum of the $\phi_{C_0}(C)$, when $C$ runs through the set of all conjugates by elements of $I_{0,F}$ of all elements of $Ch_c$ satisfying {\bf (P1)} and on the same side of $H_E$ as $C_0$ and of their images by the reflection with respect to $H_E$, is zero.
\end{lemme}

From now on, we denote by $Ch_{c,C_1}$ the set of such $C$.

Now let $C''_1$ be the other chamber adjacent to $C$ and such that $H$ contains a wall $D''_1$ of $C''_1$; we have:

\begin{lemme}
Let $C$ be any element of $Ch_c$ contained in $A_{0,E}$. The following conditions are equivalent:
\begin{itemize}
\item $C$ is either a chamber satisfying {\bf (P1)} and on the same side of $H_E$ as $C_0$ or the image by the reflection with respect to $H_E$ of such a chamber;
\item there exist minimal galleries between $C_0$ and $C$ containing $C_1$ but none containing $C''_1$.
\end{itemize}
\end{lemme}

Let $D_0$ (resp. $D'_0$) be the wall separating $C_0$ from $C_1$ (resp. $C'_1$), and let $H_0$ (resp. $H'_0$) be the hyperplane of $A_E$ containing it.  A chamber $C$ of $A_E$ satisfies the second condition if and only if it is separated from $C_0$ by $H_0$ but not by $H'_0$. On the other hand, since $H_0$, $H'_0$ and $H$ are the only three hyperplanes of $A_E$ containing $D_0\cap D'_0$, $H'_0$ must be the image of $H_0$ by the orthogonal reflection with respect to $H$. Both conditions are then equivalent to: $R(C)$ is contained either in the connected component of $R(A_E)-(R(H_E)\cup R(H_0)\cup R(H'_0))$ containing $C_1$ or in its image by the orthogonal reflection with respect to $R(H_E)$. The lemma follows immediately. $\Box$

On the other hand, since $H_0$ and $H'_0$ both contain walls of $C_0$ and are not perpendicular to each other, they correspond to consecutive roots in the extended Dynkin diagram of $\Phi$. Since, $\Phi$ being of type $A_{2n}$, its extended Dynkin diagram  is a cycle, we can label the hyperplanes $H_{0,1}, \dots,H_{0,2n+1}$ containing walls of $C_0$ in such a way that for every $i$, with $H'_{0,i}$ being defined relatively to $H_{0,i}$ the same way as $H'_0$  is defined relatively to $H_0$, we have $H'_{0,i}=H_{0,i+1}$ (the indices being taken modulo $2n+1$). More precisely, for every $i$, let $C_{1,i}$ be the chamber of $A_E$ separated from $C_0$ by $H_{0,i}$, let $D_i$ be their common wall and let $D_{F,i}$ be the wall of $C_F$ whose geometric realization contains $D_i$. Let $C'_{1,i}$ be the unique chamber of $A_E$ neighboring $C_0$, containing a wall whose geometric realization is contained in $D_{F,i}$ and distinct from $C_{1,i}$; such a chamber exists and is unique by lemma \ref{tcc0}. Let $H'_{0,i}$ be the hyperplane of $A_E$ separating $C_0$ from $C'_{1,i}$, we then have $H'_{0,i}=H_{0,i+1}$.

Let also $A_F$ be the apartment of $X_F$ whose geometric realization is $\mathcal{A}$, and for every $i$, let $H_i$ be the hyperplane of $A_E$ whose geometric realization contains $R(D_{F,i})$, let $C'_{F,i}$ be the chamber of $A_F$ separated from $C_F$ by $D_{F,i}$ and let $C'_{0,i}$ be the unique element of $Ch_c$ whose geometric realization is contained in $C'_{F,i}$.

Let now $C$ be any element of $Ch_c$ contained in $A_E$ and different from $C_0$. Let $I_C$ be the subset of the elements $i\in{\mth{Z}}/(2n+1){\mth{Z}}$ such that $C$ is separated from $C_0$ by $H_{0,i}$; since $C\neq C_0$; $I_C$ is nonempty, and since the closure of $C\cup C_0$ must contain at least one wall of $C_0$, $I_C$ is not the whole set ${\mth{Z}}/(2n+1){\mth{Z}}$ either. Hence the set $I'_C$ of elements $i$ of ${\mth{Z}}/(2n+1){\mth{Z}}$ such that $i\in I_C$ and $i+1\not\in I_C$ is nonempty.

For every $i$, set $Ch_{c,i}=Ch_{c,C_{1,i}}$, and for every $I'\subset{\mth{Z}}/(2n+1){\mth{Z}}$, set $Ch_{c,I'}=\bigcap_{i\in I'}Ch_{c,i}$; for every $I'$ and every $C\in Ch_c$ contained in $A_E$, we have $C\in Ch_{c,I'}$ if and only if $I'\subset I'_C$, and we thus obtain:
\[\sum_{C\in Ch_c}\phi_{C_0}(C)=\phi_{C_0}(C_0)+\sum_{I'\subset{\mth{Z}}/(2n+1){\mth{Z}}, I'\neq\emptyset}(-1)^{\#(I')+1}\sum_{C\in Ch_{c,I'}}\phi_{C_0}(C).\]
Since $\phi_{C_0}(C_0)=1$, to prove proposition \ref{tva2n}, it is now enough to prove the following result:

\begin{prop}
For every nonempty subset $I'$ of $\{1,\dots,2n+1\}$, we have $\sum_{C\in Ch_{c,I'}}\phi_{C_0}(C)=0$.
\end{prop}

We already know by lemma \ref{zerosing} that the assertion of the proposition holds when $I'$ is a singleton; we now have to prove it in the other cases.

First we remark that since for every $C$ and for every $i\in I'_C$, $i$ belongs to $I_C$ but $i+1$ does not, a necessary condition for $Ch_{c,I'}$ to be nonempty is that $I'$ does not contain two consecutive elements of ${\mth{Z}}/(2n+1){\mth{Z}}$. In the sequel, we assume that $I'$ satisfies that condition.

For every $i\in{\mth{Z}}/(2n+1){\mth{Z}}$, let $|i|$ be the distance between $i$ and $0$ in the cyclic group: for example, $|1|$ is $1$, and $|2n|$ is also $1$. We have:

\begin{lemme}\label{orth3}
Let $i,j\in I'$ be such that $|i-j|\geq 3$. Then all three of $H_{0,i}$, $H'_{0,i}$, $H_i$ are orthogonal to all three of $H_{0,j}$, $H'_{0,j}$, $H_j$.
\end{lemme}

Let $\varepsilon_1,\dots,\varepsilon_d$ be elements of $X^*(T)\otimes{\mth{Q}}$ defined as in \cite[plate I]{bou}. Assume the $\varepsilon_i$ are numbered in such a way that for every $i$, $H_{0,i}$ corresponds to the roots $\pm(\varepsilon_i-\varepsilon_{i+1})$. Then $H'_{0,i}$ (resp. $H_i$) corresponds to the roots $\pm(\varepsilon_{i+1}-\varepsilon_{i+2})$ (resp. $\pm(\varepsilon_i-\varepsilon_{i+2})$). The lemma follows immediately. $\Box$

This lemma proves that the union of the elements of the intersection $Ch_{c,\{i,j\}}=Ch_{c,i}\cap Ch_{c,j}$ whose geometric realization is contained in $\mathcal{A}$ is symmetrical with respect to $H_i$ (or $H_j$, for that matter); we deduce from this, using the same reasoning as for $Ch_{c,i}$ in lemma \ref{zerosing}, that $\sum_{C\in Ch_{c,\{i,j\}}}f(C)=0$. More generally, we divide $I'$ into segments the following way: $I'=I'_1\cup\dots\cup I'_r$, where every $I'_k$ is of the form $\{i,i+2,\dots,i+2(l_k-1)\}$, $l_k$ being the length of the segment, and if $i\in I'_k$ and $j\in I'_l$ with $k\neq l$, then $|i-j|\geq 3$; such a partition of $I'$ into segments exists since $I'$ cannot contain two consecutive elements of ${\mth{Z}}/(2n+1){\mth{Z}}$, and is obviously unique up to permutation of the segments. We then prove in a similar manner as for $I'=\{i,j\}$ that we have $\sum_{C\in Ch_{c,I'}}f(C)=0$ as soon as one of the $I'_k$ is a singleton.

Consider now the case where $I'$ is a single segment of length $l>1$, say for example $I'=\{1,3,\dots,2l-1\}$. Then if $C$ is an element of $Ch_{c,I'}$ contained in $\mathcal{A}$, the concave function $f_C$ associated to $C$ (normalized by taking $C_0$ as the standard Iwahori) must satisfy the following conditions:
\begin{itemize}
\item for every $i\in\{0,\dots,l-1\}$, $f_C(\varepsilon_{1+2i}-\varepsilon_{2+2i})\geq \frac 12$;
\item for every $i\in\{0,\dots,l-1\}$, $f_C(\varepsilon_{2+2i}-\varepsilon_{3+2i})\leq 0$.
\end{itemize}
Since $f_C(\alpha)+f_C(-\alpha)=1$ for every $\alpha\in\Phi$, we obtain:
\begin{itemize}
\item for every $i\in\{0,\dots,l-1\}$, $f_C(\varepsilon_{2+2i}-\varepsilon_{1+2i})\leq 0$;.
\item for every $i\in\{0,\dots,l-1\}$, $f_C(\varepsilon_{3+2i}-\varepsilon_{2+2i})\geq \frac 12$.
\end{itemize}
We can associate to $C$ the $(l+1)\times l$ matrix $M=(m_{ij})$ defined the following way: for every $i\in\{0,\dots,l\}$ and every $j\in\{1,\dots,l\}$, $m_{ij}=1$ (resp. $m_{ij}=0$) if $f_C(\varepsilon_{2j}-\varepsilon_{1+2i})\geq \frac 12$ (resp. $\leq 0$). For every $M$, let $Ch_{c,I',M}$ be the set of $C'\in Ch_{c,I'}$ which are conjugated by an element of $I_F$ to some chamber contained in $A_E$ whose associated matrix is $M$; we now prove that for every $M$, we have $\sum_{C\in Ch_{c,I',M}}f(C)=0$.

We first investigate the conditions for $Ch_{c,I',M}$ to be nonempty. From the above conditions we see that we must have $m_{i-1,i}=m_{ii}=0$ for every $i$. We now prove the following lemma:

\begin{lemme}
Assume there exist $i,i',j,j'$ such that $m_{ij}=m_{i'j'}=1$ and $m_{ij'}=m_{i'j}=0$. Then $Ch_{c,I',M}$ is empty.
\end{lemme}

Let $C$ be an element of $Ch_{c,I',M}$ contained in $A_E$. In terms of concave functions, the assertion of the lemma translates into: $f_C(\varepsilon_{2j}-\varepsilon_{1+2i}),f_C(\varepsilon_{2j'}-\varepsilon_{1+2i'})\geq \frac 12$ and $f_C(\varepsilon_{2j}-\varepsilon_{1+2i'}),f_C(\varepsilon_{2j'}-\varepsilon_{1+2i})\leq 0$. We deduce from this that we have $f_C(\varepsilon_{1+2i}-\varepsilon_{2j}))\leq 0$ and $f_C(\varepsilon_{1+2i'}-\varepsilon_{2j'})\leq 0$, hence by concavity:
\[f_C(\varepsilon_{1+2i}-\varepsilon_{1+2i'})\leq f_C(\varepsilon_{1+2i}-\varepsilon_{2j})+f_C(\varepsilon_{2j}-\varepsilon_{1+2i'})\leq 0,\]
\[f_C(\varepsilon_{1+2i'}-\varepsilon_{1+2i})\leq f_C(\varepsilon_{1+2i'}-\varepsilon_{2j'})+f_C(\varepsilon_{2j'}-\varepsilon_{1+2i})\leq 0.\]
On the other hand, since $C$ is a chamber, we must have $f_C(\varepsilon_{1+2i}-\varepsilon_{1+2i'})+f_C(\varepsilon_{1+2i'}-\varepsilon_{1+2i})=\frac 12$, which is impossible given the above inequalities. Hence $Ch_{c,I',M}$ must be empty and the lemma is proved. $\Box$

From now on we assume that $M$ is such that $Ch_{c,I',M}$ is nonempty.

\begin{cor}
For every $i$, let $Z_i$ be the set of indices $j$ such that $m_{ij}=0$. Then for every $i,i'$, we have either $Z_i\subset Z_{i'}$ or $Z_{i'}\subset Z_i$.
\end{cor}

Assume there exist $j,j'$ such that $j\in Z_{i'}-Z_i$ and $j'\in Z_i-Z_{i'}$. Then $i,i',j,j'$ satisfy the conditions of the previous lemma, and $M$ cannot then be nonempty. $\Box$

Using this corollary, we define a total preorder on $\{0,\dots,l\}$ by $i\leq_M i'$ if and only if $Z_i\subset Z_{i'}$.

\begin{lemme}
Let $i$ be a maximal element for that preorder. Then $Z_i$ is the full set $\{1,\dots,l\}$.
\end{lemme}

As we have already seen, for every $j\in\{1,\dots,n\}$, $m_{jj}=0$, hence $j\in Z_j\subset Z_i$. $\Box$

\begin{lemme}\label{immax}
There exists an $i\in\{0,\dots,l\}$ such that both $i$ and $i-1$ are maximal for the order $\leq_M$.
\end{lemme}

Let $i_0$ be any maximal element of $\{0,\dots,l\}$ for $\leq_M$. If either $i_0-1$ or $i_0+1$ is maximal, there is nothing to prove; assume that none of them is maximal. Let $j$ be an element of $Z_{i_0}-Z_{i_0+1}$; since $m_{j-1,j}=m_{jj}=0$, $j$ belongs to both $Z_j$ and $Z_{j-1}$, and we then have $i_0+1<_Mj$ and $i_0+1<_Mj-1$. If both $j$ and $j-1$ are maximal, the lemma is proved, if either $j$ or $j-1$ is not maximal, assuming for example $j$ is not, we now consider an index $k$ not belonging to $Z_j$ and we use the same reasoning as above to obtain that $j<_Mk$ and $j<_Mk-1$; since our set of indices is finite, after a finite number of iterations we must reach an $i$ such that both $i$ and $i-1$ are maximal, as desired. $\Box$

\begin{cor}\label{msymm}
Assume $i$ is such that both $i$ and $i-1$ are maximal for $<_M$. Then the set of chambers in $Ch_{c,I',M}$ contained in $A_E$ is symmetrical with respect to $H_{2i-1}$.
\end{cor}

It is easy to see that for every $i$, replacing a chamber $C$ by its image by the symmetry with respect to $H_{2i-1}$ is equivalent to switching the columns $i-1$ and $i$ in $M$. When $i-1$ and $i$ are both maximal for $\leq_M$, these columns are identical, hence $M$ is preserved. $\Box$

We can now prove that $\sum_{C\in Ch_{c,I',M}}f(C)=0$ the same way as when $I'$ is a singleton: let $i$ be an integer associated to $M$ by lemma \ref{immax}, and let $C,C'$ be the two chambers adjacent to $C_0$ and such that the geometric realization of one of their walls is contained in the geometric realization of $H_{2i-1}$ (these chambers exist by lemma \ref{tcc0}). With the help of corollary \ref{msymm}, we can now, by the same reasoning as in lemma \ref{zerosing}, obtain the desired result. Since this is true for every $M$, we obtain that $\sum_{C\in Ch_{c,I'}}f(C)=0$ when $I'$ is a single segment.

We finally use, with the help of lemma \ref{orth3}, the same reasoning applied to any one of the segments of $I'$ to prove that $\sum_{C\in Ch_{c,I'}}f(C)=0$ in the general case. $\Box$

Since by that proposition, $\sum_{C\in Ch_c}\phi_{C_0}(C)=\phi_{C_0}(C_0)\neq 0$, $\phi_{C_0}$ is a test vector for $\lambda$, and theorem \ref{th1} is now proved when $G$ is of type $A_{2n}$ and $q$ is large enough. $\Box$

Remark: in \cite{bc}, where $\lambda$ is defined in a similar way as in this subsection, since $E/F$ is unramified, the sum defining $\lambda$ converges because at every step, there are $q_E=q^2$ times more chambers on the building itself, which implies that for every $f\in\mathcal{H}(X_E)^\infty$, for chambers $C'$ located far away enough from the origin, at every step, $f(C)$ is divided by $q^2$ and we only have $q$ times more chambers to consider (see \cite[lemmas 4.3 and 4.4]{bc}). In the tamely ramified case, for the groups of type $A_{2n}$ we are considering here, there are only $q$ times more chambers on the building itself when the distance increases by $1$, but at every step, the distance increases by $2$ on average (lemma \ref{dbdist}), and the sum converges for that reason. We will see in the sequel that a similar argument applies to other types of groups as well.

\subsection{The other cases}

In this subsection, we assume that $\Phi$ is not of type $A_{2n}$ for any $n$. Let $\Sigma_a$ be a subset of $\Phi$ satisfying the conditions of proposition \ref{anismax}; we will prove that there exists a linear form $\lambda$ on $\mathcal{H}(X_E)^\infty$ with support in the $F$-anisotropy class $Ch_a$ of $Ch_E$ corresponding to $\Sigma_a$ and a test vector $f\in\mathcal{H}(X_E)^\infty$ such that $\lambda(f)=0$. Note that this time, our test vector will not be Iwahori-spherical.

Let $T$ be a $E$-split maximal $F$-torus of $G$ of $F$-anisotropy class $\Sigma_a$, let $A$ be the $\Gamma$-stable apartment of $X_E$ associated to $T$ and let $D$ be a facet of $A^\Gamma$ of maximal dimension. We assume that $D$ and $\Sigma_a$ have also been chosen in such a way that either proposition \ref{sigmaa} or (in cases $A_d$, $d$ odd, $D_d$, $d$ odd and $E_6$) proposition \ref{sigmaa2} is satisfied.

 As in the previous section, we denote by $\Phi_D$ the smallest Levi subsystem of $\Phi$ containing $\Sigma_a$; $\Phi_D$ is also the root system of $K_{D,E}/K_{D,E}^0$, where $K_{D,E}^0$ is the pro-unipotent radical of $K_{D,E}$.

Let $\mathcal{H}(Ch_D)$ be the space of harmonic cochains on $Ch_D$. First we prove that there actually exists an element of $\mathcal{H}(Ch_D)$ with support in $Ch_{D,a}$ which is stable by $K_A\cap G_{F,der}$ and not identically zero on $Ch_D$. Let $\phi_D$ be the function on $Ch_D$ defined the following way:
\begin{itemize}
\item the support of $\phi_D$ is $Ch_{D,a}$;
\item $\phi_D(C(1,\dots,1))=1$, and for every $\lambda_1,\dots,\lambda_r\in k_F^*$, $\phi_D(C(\lambda_1,\dots,\lambda_n))$ is either $1$ or $-1$, its values being chosen in such a way that, $f=\phi_D$ being viewed as a function on $H^1(\Gamma,K_{T\cap L_{E,der}})$, the relations of proposition \ref{eic} are all satisfied;
\item $\phi_D$ is $K_D\cap G_{F,der}$-stable.
\end{itemize}

First we check that the definition is consistent. The map $(\lambda_1,\dots,\lambda_r)\mapsto \phi_D(C(\lambda_1,\dots,\lambda_r))$ being a group morphism from $(k_F^*)^r$ to $\{\pm 1\}$, it is enough to prove the following lemma:

\begin{lemme}\label{gfdstab}
For every $g\in K_D\cap G_{F,der}$ such that $C'=gC(1,\dots,1)$ is of the form $C(\lambda_1,\dots,\lambda_r)$, we  have $\phi_D(C')=1$.
\end{lemme}

First we prove that we can assume $g$ is an element of $T_0\cap G_{F,der}$. Let $F'$ be the unique quadratic unramified extension of $F$; we deduce from lemma \ref{lfdc} that there exists $t\in L_{F',der}$ such that $tC'=C$, and $t$ obviously must belong to $K_{T_0,F'}$. Set $g'=gt$; $g'$ is then an element of $K_{D,F'}\subset G_{F'}$ such that $g'C=C$. On the other hand, such an element must satisfy $g'\gamma(C)=\gamma(C)$ as well, hence is contained in $K_{C\cap\gamma(C),F'}=K_{T,F'}K_{D,F'}^0\subset L_{F'}\cap K_{D,F'}$, and since $f_D(\beta)\in\frac 12+{\mth{Z}}$ for every $\beta\in\Sigma_a$, we have $U_\beta\cap K_{D,F'}\subset K_{D,F'}^0$ for every $\beta$, from which we deduce that $L_{F'}\cap K_{D,F'}\subset T_{0,F'}K_{D,F'}^0$. Hence we can assume $g'\in T_0$, which implies that $g\in T_0$ as well.

We now assume $g$ is an element of $T_0\cap G_{F,der}$, and even that $g$ is of the form $\alpha^\vee(c)$, with $\alpha$ being a simple root in $\Phi^+$ and  $c$ being an element of $\mathcal{O}_F^*$ which is not a square.

First we remark that when $\alpha\in-\Sigma_a$, say $\alpha=\beta_1$ for example, we have:
\[\alpha^\vee(c)C(1,\dots,1)=C(c^2,1,\dots,1)\]
and since we obviously have $\phi_D(C(c^2,1,\dots,1))=1$, the result follows.

Now we deal with the other simple roots with the help of a case-by-case analysis. Notations are the same as in proposition \ref{eic}.

\begin{itemize}
\item Assume $\Phi$ is of type $A_d$, with $d=2n-1$ being odd. Then the simple roots $\alpha_{2i-1}$, $i=1,\dots,n$, are all contained in $-\Sigma_a$, and when $i$ is even, for every $j$, setting $\beta_j=\alpha_{2j-1}$, $<\beta_j,\alpha_i^\vee>$ is $-1$  if $j$ is either $\frac i2$ or $\frac i2+1$, and $0$ in the other cases; we then have:
\[\alpha_i^\vee(c)C(1,\dots,1)=C(1,\dots,c^{-1},c^{-1},\dots,1),\]
the $c^{-1}$ being in $j$-th and $j+1$-th position; hence in $H^1(\Gamma,K_{T\cap L_{E,der}})$, we obtain $\alpha_i^\vee(c)=e_je_{j+1}$. By proposition \ref{eic}, for every $\lambda_1,\dots,\lambda_n$, we have:
\[\phi_D(C(\lambda_1,\dots,\lambda_n))=(-1)^s\phi_D(C(1,\dots,1))=(-1)^s,\]
where $s$ is the number of $\lambda_i$ which are not squares; the result follows immediately.
\item Assume $\Phi$ is of type $B_d$. Then the simple roots $\alpha_i$, with $i$ odd, are all contained in $-\Sigma_a$. On the other hand, when $i$ is even and strictly smaller than $d$, $\alpha_i$ has already been dealt with in propostion \ref{sract}. I will explicit  what it means in this case, the other cases being treated similarly. By the relations we have found in proposition \ref{sract}, for every such $i$, we have, in $H^1(\Gamma,K_{T\cap L_{E,der}})$, $\alpha_i^\vee(c)1=e_{i-1}e_ie_{i+1}e_{i+2}$, and we deduce immediately from proposition \ref{eic} that $\phi_D(\alpha_i^\vee(c)C(1,\dots,1))=1$; which is the expected result. When $d$ is odd, the result is now proved, and when $d$ is even, it only remains to consider $\alpha_d^\vee$. We have $<\beta_i,\alpha_d^\vee>=-2$ if $i$ is either $d-1$ or $d$ and $0$ in the other cases, hence:
\[\alpha_d^\vee(c)C(1,\dots,1)=C(1,\dots,1,c^{-2},c^{-2}).\]
The result follows immediately.
\item Assume $\Phi$ is of type $C_d$. The only simple root contained in $-\Sigma_a$ is then $\alpha_d$, and for every $i<d$, $<\alpha_i,\beta_j>$ is $-2$ if $j$ is either $i$ or $i+1$ and $0$ else, hence we have:
\[\alpha_i^\vee(c)C(1,\dots,1)=C(1,\dots,c^{-2},c^{-2},\dots,1).\]
The result follows.
\item Assume $\Phi$ is of type $D_d$. The simple roots contained in $-\Sigma_a$ are the $\alpha_i$, with $d-i$ odd, and $\alpha_d$. The $\alpha_i$, with $d-i$ even and $1<i<d$, have already been dealt with in proposition \ref{sract}, and when $d$ is odd, we have in $H^1(\Gamma,K_{T\cap L_{E,der}})$, by proposition \ref{sract}, $\alpha_1^\vee(c)1=e_1e_2$. On the other hand, by proposition \ref{eic}, we have:
\[\phi_D(e_1e_2)=\phi_D(1).\]
The result follows.
\item Assume $\Phi$ is of type $E_6$. The simple roots contained in $-\Sigma_a$ are $\alpha_2$, $\alpha_3$ and $\alpha_5$, and $\alpha_4$ has already been dealt with in proposition \ref{sract}. Now consider $\alpha_1$; we have:
\[\alpha_1^\vee(c)C(1,1,1,1)=C(c^{-1},c^{-1},1,1).\]
hence in $H^1(\Gamma,K_{T\cap L_{E,der}})$, we have $\alpha_1^\vee(c)1=e_1e_2$. On the other hand, by proposition \ref{eic}, we have $\phi_D(e_1e_2)=\phi_D(1)$.  The case of $\alpha_6$ being symmetrical, the result follows.
\item Assume $\Phi$ is of type $F_4$. The only simple root contained in $-\Sigma_a$ is $\alpha_2$, and $\alpha_1$ has already been dealt with in proposition \ref{sract}. On the other hand, we have:
\[\alpha_3^\vee(c)C(1,1,1,1)=C(1,1,c^2,c^{-2});\]
\[\alpha_4^\vee(c)C(1,1,1,1)=C(1,c^2,c^{-2},1).\]
The result follows immediately.
\item In the three remaining cases ($E_7$, $E_8$ and $G_2$), every simple root either belongs to $-\Sigma_a$ or has been dealt with in proposition \ref{sract}; these cases then follow immediately from that proposition.
\end{itemize}
The lemma is now proved. $\Box$

Now we check that $\phi_D$ satisfies the harmonicity condition.

\begin{prop}\label{phidhc}
 Let $D_1$ be any codimension $1$ facet of $X_E$ containing $D$; the sum of the values of $\phi_D$ on the chambers containing $D_1$ is zero.
\end{prop}

If $D_1$ is not contained in any element of $Ch_{D,a}$, the harmonicity condition is trivially satisfied; we can thus assume that $D_1$ is contained in some $C\in Ch_{D,a}$, and even, by eventually conjugating it, in some $C\in Ch_{D,a,L,C_0}$. Let $D'$ be the unique codimension $1$ facet of $C_0$ of the same type as $D_1$, or in other words the only one which is $G_{E,der}$-conjugated to $D_1$. Let $\alpha$ be the corresponding simple root in $\Phi_D^+$; assume first there exists a conjugate $\Sigma'$ of $\Sigma_a$ in $\Phi_D$ containing $\alpha$. Since $\alpha$ is a simple root, by definition of $f_{C_0}$, we have $f_{C_0}(-\alpha)= 1\in{\mth{Z}}$.

Let $\Phi'_D{}^+$ be any set of positive roots of $\Phi_D$ such that $\alpha$ is a simple root in $\Phi'_D{}^+$, and let $C'_0$ be the unique chamber of $A_{0,E}$ containing $D$ such that $-\Phi'_D{}^+$ is the set of roots of the Borel subgroup of $K_{D,E}/K_{D,E}^0$ corresponding to it. For every $\lambda_1,\dots,\lambda_r\in\mathcal{O}_E^*$, we define the chamber $C'(\lambda_1,\dots,\lambda_r)\in Ch_{D,a,L,C'_0}$ in a similar way as $C(\lambda_1,\dots,\lambda_r)$. Since $f_{C_0}(-\alpha)$ is an integer, by proposition \ref{chdalconj}, there exist $\lambda_1,\dots,\lambda_r$ such that $C'(\lambda_1,\dots,\lambda_r)$ is $K_D\cap G_{F,der}$-conjugated to $C$.

Let $D'_1$ be the codimension $1$ facet of $C'(\lambda_1,\dots,\lambda_n)$ of the same type as $D'$; $D'_1$ and $D'$ are then $G_{F,der}$-conjugates, which implies that every chamber of $X_E$ containing $D'_1$ is then $G_{F,der}$-conjugated to some chamber of $X_E$ containing $D'$; and that these conjugations induce a bijection between these two set of chambers; the harmonicity condition for the chambers containing $D'_1$, which follows from lemma \ref{sroot}, then implies the hamonicity condition for those containing $D_1$.

On the other hand, two roots of the same length are always conjugates, hence the condition on $\alpha$ holds as soon as $\Sigma_a$ contains roots of every length. This is trivially true when $\Phi$ is simply-laced, and we see from proposition \ref{sigmaa} that it is also true for types $B_d$, $d$ odd, and $G_2$.

Assume now we are in one of the remaining cases ($B_d$ with $d$ even, $C_d$ for any $d$ and $F_4$); $\Sigma_a$ then contains only long roots, and the above proof still works when $\alpha$ is long. Assume now $\alpha$ is short, and let $\beta$ be a long root belonging to $\Phi^+$ and not orthogonal to $\alpha$; $\alpha$ and $\beta$ then generate a subystem of $\Phi$ of type $B_2$, hence either $\beta+2\alpha$ or $\beta-2\alpha$ is also a long root, and that root must also belong to $\Phi^+$ (it is obvious for $\beta+2\alpha$; for $\beta-2\alpha$, as in lemma \ref{shoev}, it comes from the fact that $\beta$ contains at least one simple root different from $\alpha$ in its decomposition, hence $\beta-2\alpha$ cannot be negative). In both cases, $\alpha$ is the half-difference of two long roots belonging to $\Phi^+$, and we are then in the situation of lemma \ref{cedeux}; the harmonicity condition for $D_1$ then follows immediately from the expression (\ref{bedeux}) in the proof of that lemma. $\Box$

Now we check that $\phi_D$ is compatible with the Prasad character $\chi$, or in other words thet $\phi_D(gC)=\chi(g)\phi_D(C)$ for every $g\in G_F$ and every $C\in Ch_D$. Let $K_{T_0,F}$ be the maximal compact subgroup of $(T_0)_F$ and let $X_{T_0,F}$ be the subgroup of $(T_0)_F$ generated by the $\xi(\varpi_F)$, where $\xi$ runs over the one-parameter subgroups of $T_0$.

Remember that we have a decomposition $G_F=G_{F,der}K_{T_0,F}X_{T_0,F}$, and also that the $\varpi_F$ we have chosen is the norm of some element of $E$. The character $\chi$ is trivial on $G_{F,der}$ and on $X_{T_0,F}$; the compatibility of $\phi_D$ with $\chi$ is then an immediate consequence of  the following proposition:

\begin{prop}
Let $t$ be any element of the maximal compact subgroup $K_{T_0}$ of $T_{0,F}$. Then for every $C\in Ch_D$ and every $f\in\mathcal{H}(X_E)^{G_{F,der}}$, we have $f(tC)=\chi(t)f(C)$.
\end{prop}

Let $C$ be any element of $Ch_D$. If $C$ does not belong to $Ch_{D,a}$, then neither does $tC$ and we then have $f(tC)=\chi(t)f(C)=0$. We thus may assume that $f\in Ch_{D,a}$, and by eventually conjugating it, we can even assume that $C$ belongs to $Ch_{D,a,L,C_0}$.

We already know from lemma \ref{gfdstab} that if $t\in G_{F,der}$, $f(tC)=f(C)=\chi(t)f(C)$ since $\chi(t)=1$. Moreover, if $t$ is a square, then its image in $H^1(\Gamma,K_{T\cap L_{E,der},F})$ is a square too, hence trivial by proposition \ref{clcoh}, and since $\chi$ is quadratic, $\chi(t)$ is trivial too. Hence we only have to prove the result when $t$ belongs to some set of representatives in $T_0\cap G_F$ of some set of generators of the finite abelian group $Y/Y^2$, where $Y=(T_0\cap G_F)/(T_0\cap G_{F,der})$.

Let $\rho$ be the half-sum of the elements of $\Phi^+$; by \cite[\S 1, proposition 29]{bou} and \cite[lemma 3.1]{cou2}, for every $t\in T_0\cap G_F$, $\chi(t)=1$ if and only if $2\rho(t)$ is the norm of some element of $E^*$. We refer to \cite[plates I to IX]{bou} for the expressions of $Y$ and $\rho$ we use during the case-by-case analysis below. In the sequel, once again, $c$ is an element of $\mathcal{O}_F^*$ which is not a square.

Note first that the cases $E_8$, $F_4$ and $G_2$ are trivial since we then have $G_F=G_{F,der}$. We examine the other cases.

\begin{itemize}
\item Assume $\Phi$ is of type $A_{2n-1}$. Then $Y$ is cyclic of order $2n$, and with a slight abuse of notation, the element $t=Diag(c,1,\dots,1)$ of $GL_{2n}(F)$ is, for any choice of $c$, a representative of the unique nontrivial element of $Y/Y^2$, hence can be used to compare two quadratic characters of $Y$. Since $\beta_1=\alpha_1$, for every $h\in H^1(\Gamma,K_{T\cap L_{E,der},F})$ and every $C\in Ch(h)$, by proposition \ref{sract}, the chamber $tC$ belongs to $Ch(e_1h)$, and we deduce from proposition \ref{eic} that $\phi_D(tC)=-\phi_D(C)$. On the other hand, we have $2\rho=\sum_{i=1}^di(d+1-i)\alpha_i$, hence $2\rho(t)=c^d$, hence $\chi(t)=(-1)^d$. We thus obtain $\phi_D(tC)=\chi(t)\phi_D(C)$, as desired.

\item Assume now $\Phi$ is of type $B_d$; $Y$ is then of order $2$ and its nontrivial element admits $t=Diag(c,1,\dots,1,c^{-1})\in GSO'_{2d+1}(F)$, where $GSO'_{2d+1}$ is the split form of $GSO_{2d+1}$, as a representative. We denote by $n$ the largest integer such that $2n\leq d$.

By proposition \ref{sract}, for every $h\in H^1(\Gamma,K_{T\cap L_{E,der},F})$ and every $C\in Ch(h)$, $tC$ belongs to $Ch(e_1e_2h)$, hence by proposition \ref{eic}, $\phi_D(tC)=-\phi_D(C)$, On the other hand, we have $2\rho=\sum_{i=1}^di(2d-i)\alpha_i$; we then obtain $2\rho(t)=c^{2d-1}$, hence $\chi(t)=-1$ and the result follows.

\item Assume $\Phi$ is of type $C_d$; $Y$ is then of order $2$,  and with a slight abuse of notation, its nontrivial element admits $t=Diag(c,\dots,c,1,\dots,1)\in GSp_{2d}(F)$ as a representative. By proposition \ref{sract}, for every element $h$ of $H^1(\Gamma,K_{T\cap L_{E,der},F})$ and every $C\in Ch(h)$, $tC$ belongs to $Ch(e_1\dots e_dh)$, hence by proposition \ref{eic}, we have:
\[\phi_D(tC)=(\prod_{i=1}^d(-1)^{d+1-i})\phi_D(C)\]
\[=(\prod_{i=1}^d(-1)^i)\phi_D(C.)=(-1)^{\frac{d(d+1)}2}\phi_D(C).\]
On the other hand, we have:
\[2\rho=\sum_{i=1}^{d-1}i(2d+1-i)\alpha_i+\frac{d(d+1)}2\alpha_d;\]
we then obtain $2\rho(t)=c^{\frac{d(d+1)}2}$, hence $\chi(t)=(-1)^{\frac{d(d+1)}2}$. We finally get $\phi_D(tC)=\chi(t)\phi_D(C)$ once again.

\item Assume now $\Phi$ is of type $D_d$. When $d$ is even, $Y$ is isomorphic to $({\mth{Z}}/2{\mth{Z})}^2$ and is generated by the elements admitting respectively $t=Diag(c,1,\dots,1,c^{-1})$ and $t'=Diag(c,\dots,c,1,\dots,1)$, both belonging to $GSO'_{2d}(F)$, as representatives; when $d$ is odd, $Y$ is cyclic of order $4$ and one of its generators admits $t'$ as a representative. In both cases, we denote by $n$ the largest integer such that $2n\leq d$.

When $d$ is even, by proposition \ref{sract}, for every $h\in H^1(\Gamma,K_{T\cap L_{E,der},F})$ and every $C\in Ch(h)$, $tC$ belongs to $Ch(e_1e_2h)$, hence by proposition \ref{eic}, $\phi_D(tC)=\phi_D(C)$. On the other hand, we have:
\[2\rho=\sum_{i=1}^{d-2}i(2d-1-i)\alpha_i+\frac{d(d-1)}2(\alpha_{d-1}+\alpha_d);\]
 hence $2\rho(t)=c^{2d-2}$, from which we obtain that $\chi(t)=1$ and that $\phi_D(tC)=\chi(t)\phi_D(C)$, as desired.

Now we consider $t'$, $d$ being either odd or even. By proposition \ref{sract}, we have $t'Ch(h)=Ch(e_1e_3\dots e_{2n-1}h)$, hence by proposition \ref{eic}, $\phi_D(t'C)=(-1)^n\phi_D(C)$; on the other hand, using the same expression as above for $2\rho$, we obtain $\chi(t')=(-1)^{\frac{d(d-1)}2}$.

To prove the result, we thus only have to check that $n$ and $\frac{d(d-1)}2$ have the same parity. When $d$ is even, then $d=2n$, and $\frac{d(d-1)}2=n(d-1)$ and $n$ have the same parity. When $d$ is odd, then $d-1=2n$, and $\frac{d(d-1)}2=nd$ and $n$ also have the same parity. The result follows.

\item Assume now $\Phi$ is of type $E_6$. The character $\xi=\frac{\alpha_1^\vee-\alpha_3^\vee+\alpha_5^\vee-\alpha_6^\vee}3$ is then an element of $X_*(T)$, and if $t=\xi(c)$, we have $tCh(h)=Ch(h)$ for every $h\in H^1(\Gamma,K_{T\cap L_{E,der}})$, hence $\phi_D(tC)=\phi_D(C)$. On the other hand, since the group $Y$ is of order $3$ and $\chi$ is quadratic, it must be trivial, hence $\chi(\xi(c))=1$, and the result follows.

\item Assume now $\Phi$ is of type $E_7$. The group $Y$ is then of order $2$; moreover, the character $\xi=\frac{\alpha_2^\vee+\alpha_5^\vee+\alpha_7^\vee}2$ is an element of $X_*(T)$, and $t=\xi(c)$ is a representative of the nontrivial element of $Y$. By proposition \ref{sract}, we have, for every $h$, $\xi(c)Ch(h)=Ch(e_4e_6e_7h)$; hence, by proposition \ref{eic}, $\phi_D(tC))=-\phi_D(C)$. On the other hand, since by \cite[\S 1, proposition 29 (ii)]{bou}, $<\rho,\alpha_i^\vee>=1$ for every $i$, we obtain $<2\rho,\xi>=3$, hence $2\rho(t)=c^3$ and $\chi(t)=-1$, and the result follows.
\end{itemize}
$\Box$

Now we can define our linear form.
For every $C\in Ch_{D,a}$, let $O_C$ be the $G_F$-orbit of $Ch_E$ containing $C$, and let $R$ be a system of representatives of the $G_F$-orbits in $Ch_{D,a}$. Set:
\[\lambda:f\in\mathcal{H}(X_E)^\infty\longmapsto\sum_{C\in R}\sum_{C'\in O_C}f(C')\phi_D(C').\]
Since $R$ is a finite set, proposition \ref{elun} implies that the double sum always converges.

%

Now that we have a linear form on $\mathcal{H}(X_E)^\infty$, we want to find a test vector for it.
We start by the following propositions:

\begin{prop}
Let $C$ be any element of $Ch_E$. There exists a unique element of $Ch_D$ contained in the closure of $C\cup D$.
\end{prop}

Since $C$ is a chamber, by \cite[2.4.4]{bt}, the closure $cl(C\cup D)$ is a union of chambers of $X_E$. Hence $D$ is contained in some chamber $C'$ of that closure, which is then obviously an element of $Ch_D$.

On the other hand, let $A_C$ be an apartment of $X_E$ containing both $C$ and $D$; it then contains $cl(C\cup D)$. Consider the connected components of the complementary in $R(A_C)$ of the union of the walls containing $R(D)$; each one of them contains the geometric realization of a unique element of $Ch_D$. Let $\mathcal{S}$ be the one containing $R(C)$; its closure contains $R(D)$, hence also the geometric realization of $cl(C\cup D)$, which proves the unicity of $C'$. $\Box$


\begin{prop}\label{extharm}
Let $f_0$ be a function on $Ch_D$ satisfying the harmonicity condition, and let $f$ be the function on $Ch_E$ defined the following way: for every $C\in Ch_E$, if $C_0$ is the only element of $Ch_D$ contained in the closure of $C\cup D$, $f(C)=(-q)^{-d(C,C_0)}f_0(C_0)$. Then $f\in\mathcal{H}(X_E)^\infty$.
\end{prop}

Let $K_{D,E}^0$ be the pro-unipotent radical of $K_{D,E}$; $K_{D,E}^0$ fixes every element of $Ch_D$ pointwise. For every $C\in Ch_E$ and every $k\in K_{D,E}^0$, we then have:
\[f(kC)=(-q)^{-d(kC,kC_0)}f_0(kC_0)=(-q)^{-d(C,C_0)}f_0(C_0)=f(C);\]
since $K_{D,E}^0$ is an open compact subgroup of $G_E$, the smoothness of $f$ is proved. Now we check the harmonicity condition. Let $D'$ be any codimension $1$ facet of $X_E$. Assume first that the closure of $D'\cup D$ contains at least one chamber $C_1$ of $X_E$; it then contains exactly one element $C_0$ of $Ch_D$, namely the one whose geometric realization is contained in the  same connected component as $R(C_1)$ of the complementary of the union of the walls containing $R(D)$ in any $\Gamma$-stable apartment containing $R(C_1)$; on the other hand, that closure also contains exactly one chamber $C$ admitting $D'$ as a wall. Set $\delta=d(C,C_0)$; if $C'$ is any other chamber of $X_E$ admitting $D'$ as a wall, the closure of $C'\cap D$ contains $D$ and $C$, hence contains also $C_0$, and we have $d(C',C_0)=\delta+1$. Since there are $q$ such chambers, we obtain:
\[\sum_{C'\supset D'}f(C')=(-q)^{-\delta}f_0(C_0)+q(-q)^{-\delta-1}f_0(C_0)=0.\]
Assume now that the closure of $D'\cup D$ does not contain any chamber. It then contains a unique facet $D_0$ of $X_E$ of codimension $1$ containing $D$; moreover, if $C$ is a chamber of $X_E$ admitting $D'$ as a wall, the only element $C_0$ of $Ch_D$ contained in the closure of $C\cup D$ must admit $D_0$ as a wall. On the other hand, the group $K_{D'\cup D}$ permutes transitively the elements of $Ch_D$ admitting $D_0$ as a wall; since there are $q+1$ such chambers, and $q+1$ chambers of $X_E$ admitting $D'$ as a wall as well, the restriction to the second ones of the application $C\mapsto C_0$ must be a bijection, and all of them are at the same distance $\delta$ from $Ch_D$. We then have:
\[\sum_{C\supset D'}f(C)=(-q)^{-\delta}\sum_{C_0\supset D_0}f_0(C_0).\]
Since $f_0$ satisfies the harmonicity condition as a function on $Ch_D$, the right-hand side is zero, hence the left-hand side must be zero as well. Hence $f$ satisfies the harmonicity condition and the proposition is proved. $\Box$

Now let $\phi$ be the function on $Ch_E$ derived from $\phi_D$ by the previous proposition. We say that $\phi$ is the {\em extension by harmonicity} of $\phi_D$.

\begin{prop}\label{tvoth}
The function $\phi$ belongs to $\mathcal{H}(X_E)^\infty$, and is a test vector for $\lambda$.
\end{prop}

The fact that $\phi\in\mathcal{H}(X_E)^\infty$ is an immediate consequence of propositions \ref{phidhc} and \ref{extharm}. Now we prove that $\phi$ is a test vector for $\lambda$. First assume $D$ is a single vertex $x$; we then write $Ch_x$, $Ch_{x,a}$, $\phi_x$ instead of $Ch_D$, $Ch_{D,a}$, $\phi_D$. We first prove the following lemma:

\begin{lemme}
Let $C$ be an element of $Ch_a^0$ such that $\phi(C)\neq 0$. Then $C\in Ch_{x,a}$.
\end{lemme}

Assume $C\not\in Ch_{x,a}$; there exists then another vertex $x'$ of $X_E$ whose geometric realization is in $\mathcal{B}_F$, belonging to $C$ and such that $C\in Ch_{x',a}$. Let $C_0$ be the only element of $Ch_x$ contained in the closure of $C\cup\{x\}$; the closure of $C_0$ must then contain a facet of dimension at least $1$ of the closure of $\{x,x'\}$, whose geometric realization is contained in $\mathcal{B}_F$. Hence $C_0$ cannot belong to $Ch_{x,a}$, which implies that $\phi(C_0)$ must be zero, and $\phi(C)$ is then also zero by definition of $\phi$. $\Box$

According to this lemma, we have:
\[\lambda(\phi)=\sum_{C\in Ch_{x,a}}\phi_x(C)\phi(C)=\sum_{C\in Ch_{x,a}}1=\#(Ch_{x,a}).\]
Since $Ch_{x,a}$ is nonempty, $\lambda(\phi)\neq 0$ and the proposition is proved.

Now we deal with the cases where $D$ is of nonzero dimension. As before, we denote by $\Phi_D$ the root system of ${\mth{G}}_D$, which is also  the Levi subsystem of $\Phi$ generated by $\Sigma_a$, or equivalently the set of elements of $\Phi$ which are linear combinations witn coefficients in ${\mth{Q}}$ of the elements of $\Sigma_a$.

Remember that $Ch_a^0$ is the set of chambers of anisotropy class $\Sigma_a$ containing a $\Gamma$-fixed facet of the same dimension as $D$.

\begin{lemme}
Let $C$ be an element of $Ch_a^0$ such that $\phi_D(C)\neq 0$, and let $D'$ be the $\Gamma$-fixed facet of $C$ of maximal dimension. There exists a $\Gamma$-stable apartment $A$ of $X_E$ containing both $D$ and $C$, hence also $D'$, and $D$ and $D'$ are then facets of maximal dimension of $A^\Gamma$.
\end{lemme}

Let $C_0$ be the only element of $Ch_D$ contained in the closure of $C\cup D$; by definition of $\phi_D$, we must have $C_0\in Ch_{D,a}$, which implies that the intersection of $C_0$ and $\gamma(C_0)$ is $D$. Moreover, $C$ is also the only element of $Ch_{D'}$ contained in the closure of $C_0\cup D'$, hence $\gamma(C)$ is the only element of $Ch_{D'}$ contained in the closure of $\gamma(C_0\cup D')=\gamma(C_0)\cup D'$.

Consider now the closure of $\gamma(C_0)\cup C$; it contains both $\gamma(C_0)\cup D'$ and $C\cup D$, and by the previous remarks it must contain $C_0$ and $\gamma(C)$ as well, hence also the closure of $C_0\cup\gamma(C)$; by symmetry, these two closures are then equal. We have thus obtained a $\Gamma$-stable subset of $X_E$ which is the closure of the union of two facets; by \cite[proposition 2.3.1]{bt}, that set is contained in some apartment $A'$ of $X_E$, and by the same inductive reasoning as in proposition \ref{apst}, we obtain a $\Gamma$-stable apartment $A$ containing it, which must then satisfy the required conditions. $\Box$

Let $A$ be a $\Gamma$-stable apartment of $X_E$ containing at least one chamber belonging to $Ch_D$, and let $D',D''$ be facets of maximal dimension of $A^\Gamma$. We denote by $d_\Gamma(D',D'')$ the combinatorial distance between $D'$ and $D''$ inside the subcomplex $A^\Gamma$ of $X_E$.

\begin{lemme}\label{ming}
Let $D',D''$ be two distinct facets of maximal dimension of $A^\Gamma$, and let $C'$ be a chamber of $A$ containing $D'$.
\begin{itemize}
\item The parahoric subgroups $K_{D',E}$ and $K_{D'',E}$ of $G_E$ fixing respectively $D'$ and $D''$ are strongly associated (in the sense of \cite[definition 3.1.1]{deb}).
\item  There exists a unique chamber $C''$ of $A$ containing $D''$ and such that no minimal gallery between $C'$ and $C''$ contains any other chamber containing either $D'$ or $D''$.
\end{itemize}
\end{lemme}

Since $D'$ and $D''$ both generate $A^\Gamma$ as an affine subcomplex of $X_E$, the finite reductive groups $K_{D',E}/K_{D',E}^0$ and $K_{D'',E}/K_{D'',E}^0$ are both canonically isomorphic to $K_{A^\Gamma,E}/K_{A^\Gamma,E}^0$ , and we have:
\[K_{D',E}=K_{A^\Gamma,E}K_{D',E}^0,\]
from which we deduce:
\[(K_{D',E}\cap K_{D'',E})K_{D',E}^0\supset K_{A^\Gamma,E}K_{D',E}^0=K_{D',E}.\]
The other inclusion being obvious, we obtain in fact an equality. By switching $D'$ and $D''$ in the previous reasoning, we also obtain:
\[(K_{D',E}\cap K_{D'',E})K_{D'',E}^0=K_{D'',E}.\]
Hence $K_{D',E}$ and $K_{D'',E}$ are strongly associated, as desired. It implies in particular that $K_{D',E}/K_{D',E}^0$ and $K_{D'',E}/K_{D'',E}^0$ are canonically isomorphic to each other.

Now we prove the second assertion. We first observe that the image of $K_{C,E}\subset K_{D',E}$ in $K_{D',E}/K_{D',E}^0$ is a Borel subgroup of $K_{D',E}/K_{D',E}^0$. Let now $C''$ be the chamber of $A$ containing $D''$ and such that the image of $K_{C'',E}$ in $K_{D'',E}/K_{D'',E}^0$ is (up to the aforementioned canonical isomorphism) that same Borel subgroup. Assume there exists a minimal gallery $(C'_0=C',C'_1,\dots,C'_r=C'')$ between $C'$ and $C''$ such that $C'_i$ contains either $D'$ or $D''$, say for example $D'$, for some $i\in\{1,\dots,r-1\}$. Then $K_{C'_i,E}$ is contained in $K_{D',E}$, and its image in $K_{D',E}/K_{D',E}^0$ is a Borel subgroup which must be different from $K_{C',E}/K_{D',E}^0$ since $C'_i\neq C'$; hence $C'$ and $C'_i$ are separated by at least one hyperplane $H$ of $A$ containing $D'$. Such a hyperplane must then contain the whole subcomplex $A^\Gamma$, and in particular $D''$, and since $H$ then also separates $C'_i$ from $C''$, the gallery has to cross it at least twice, which contradicts its minimality.

Now let $C'''$ be another chamber satisfying the conditions of the second assertion. Since $D'$ and $D''$ are distinct, we must have $C'''\neq C'$. On the other hand, let $H$ be an hyperplane separating $C''$ from $C'''$. Since both $C''$ and $C'''$ contain $D''$, $H$ must contain $D''$ as well, hence $K_{C''',E}/K_{D'',E}^0$ is a Borel subgroup of $K_{D'',E}/K_{D'',E}^0$ which is different from $K_{C'',E}/K_{D'',E}^0\simeq K_{C',E}/K_{D',E}^0$; we deduce from this that there must exist a minimal gallery between $C'$ and $C'''$ containing $C''$, and this is possible only if $C'''=C''$. The lemma is now proved. $\Box$

\begin{lemme}\label{R1}
Let $D',D''$ be two facets of maximal dimension of $A^\Gamma$, let $C'$ be a chamber of $A$ containing $D'$ and let $C''$ be the only chamber of $A$ containing $D''$ and contained in the closure of $C'\cup D''$. Then $\frac{d(C',C'')}{d_\Gamma(D',D'')}$ is a positive integer $r_1$ which does not depend on the choice of $D'$, $D''$ and $C'$.
\end{lemme}

It is easy to prove (by for example \cite[lemma 4.2]{bc} and an obvious induction) that $q^{d(C',C'')}=[K_{C',E}:K_{C'\cup C'',E}]=[K_{C',E}:K_{C',E}\cap K_{C'',E}]$; moreover, we deduce immediately from the first assertion of lemma \ref{ming} that $[K_{C',E}:K_{C',E}\cap K_{C'',E}]=[K_{D',E}:K_{D',E}\cap K_{D'',E}]$. We thus only have to relate that last quantity to $d_\Gamma(D',D'')$.

As usual, we can without loss of generality assume that $A^\Gamma$ is contained in $A_{0,E}$. Assume first $D'$ and $D''$ are adjacent. Let $\Phi_{D'}$ be the Levi subsystem of $\Phi$ corresponding to the root system of $K_{D',E}/K_{D',E}^0$, which we can without loss of generality assume to be standard, and let $\alpha$ be any positive element of $\Phi$ corresponding to an hyperplane of $\mathcal{A}_0$ separating $D'$ from $D''$; the set of such hyperplanes is then precisely the set of elements of $\Phi^+$ contained in $\alpha+X_{D'}$, where $X_{D'}$ is the subgroup of $X^*(T_0)$ generated by $\Phi_{D'}$. We thus only have to check that the cardinality of $\Phi_{D',D''}=\Phi\cap(\alpha+X_{D'})$ is always the same.
\begin{itemize}
\item When $\Phi$ is of type $A_{2n-1}$, the simple roots contained in $\Phi_{D'}$ are the $\alpha_i$ with $i$ odd. We then have $\Phi_{D',D''}=\{\alpha_{2i},\alpha_{2i-1}+\alpha_{2i},\alpha_{2i}+\alpha_{2i+1},\alpha_{2i-1}+\alpha_{2i}+\alpha_{2i+1}\}$ for some $i$, and in particular $\Phi_{D',D''}$ always has $4$ elements. 
\item When $\Phi$ is of type $D_{2n+1}$, every simple root in $\Phi^+$ except $\alpha_1$ is contained in $\Phi_{D'}$. The set $\Phi_{D',D''}$ is then the full set of the elements of $\Phi^+$ which do not belong to $\Phi_{D'}$; there are $4n$ such roots, which are precisely the roots of the form $\varepsilon_1\pm\varepsilon_i$, $2\leq i\leq 2n+1$.
\item When $\Phi$ is of type $E_6$, the simple roots contained in $\Phi_{D'}$ are the $\alpha_i$ with $2\leq i\leq 5$. The set $\Phi_{D',D''}$ then contains every positive element of the Levi subsystem of $\Phi$ generated by $\Phi_{D'}$ and $\alpha_j$, with $j$ being either $1$ or $6$, which do not belong to $\Phi_{D'}$. Since in both cases this Levi subsystem is of type $D_5$, we are reduced to the previous case with $n=2$, and we obtain in particular that the cardinality of $\Phi_{D',D''}$ is always $8$.
\end{itemize}
In all these cases, the cardinality of $\Phi_{D',D''}$ is an integer $r_1$ which does not depend on the choice of $D'$ and $D''$.

Now we prove the general case by induction on $d_\Gamma(D',D'')$. Assume $d_\Gamma(D',D'')>1$ and let $D'''$ be a facet of maximal dimension of $A^\Gamma$ distinct from $D'$ and $D''$ and such that $d_\Gamma(D',D'')=d_\Gamma(D',D''')+d_\Gamma(D''',D'')$; $D'''$ is then contained in the closure of $D'\cap D''$, hence also in the closure of $C'\cap C''$, and that closure must then contain an element $C'''$ of $Ch_{D'''}$, which implies that $d(C',C'')=d(C',C''')+d(C''',C'')$. By induction hypothesis we have $d(C',C''')=r_1d_\Gamma(D',D'')$ and $d(C''',C'')=r_1d_\Gamma(D''',D'')$, hence $d(C',C'')=r_1d_\Gamma(D',D'')$ and the lemma is proved. $\Box$

\begin{lemme}\label{R2}
Let $D'$ be a facet of maximal dimension of $\mathcal{A}^\Gamma$. There exists an integer $r_2$ such that for every facet of maximal dimension $D''$ of $\mathcal{A}^\Gamma$, the number of $K_{D',F}$-conjugates of $D''$ is precisely $q^{r_2d_\Gamma(D,'D'')}$. Moreover, we have $r_2<r_1$.
\end{lemme}

The number of $K_{D',F}$-conjugates of $D''$ is precisely equal to $[K_{D',F}:K_{D',F}\cap K_{D'',F}]$, which cannot be greater than $[K_{D',E}:K_{D',E}\cap K_{D'',E}]=q^{r_1d_\Gamma(D,'D'')}$. Hence we already know that if $r_2$ exists, then $r_2\leq r_1$.

By the same induction as in lemma \ref{R1} we are reduced to the case where $D'$ and $D''$ are adjacent. We define $\Phi_{D',D''}$ the same way as in that lemma. Let $f_{D'}$ be the concave function on $\Phi$ associated with $D'$; $K_{D',F}/(K_{D',F}\cap K_{D'',F})$ is then generated by the images of the root subgroups of $K_{D',F}$ corresponding to elements $\alpha$ of $\Phi_{D',D''}$ such that $f_{D'}(\alpha)$ is an integer, which by definition of $C_0$ and $D'$ is true if and only if $\alpha$ is the sum of an even number of simple roots of $\Phi^+$. We thus only have to examine the different cases:
\begin{itemize}
\item when $\Phi$ is of type $A_{2n-1}$, $\Phi_{D',D''}$ always contains two such elements (either $\alpha_{2i-1}+\alpha_{2i},\alpha_{2i}+\alpha_{2i+1}$ or $\alpha_{2i},\alpha_{2i-1}+\alpha_{2i}+\alpha_{2i+1}$, depending on $D'$);
\item when $\Phi$ is of type $D_{2n+1}$, the elements of $\Phi_{D',D''}$ satisfying that condition are the $\varepsilon_1\pm\varepsilon_i$ with $i$ being of some given parity (which depends on $D'$), and there are $2n$ such roots;
\item when $\Phi$ is of type $E_6$, we are once again reduced to the case $D_5$ and $\Phi_{D',D''}$ then contains $4$ elements satisfying the required condition.
\end{itemize}
Hence in all these cases, $r_2$ exists and is strictly smaller than $r_1$, as required. $\Box$

Remark: in all cases, we have $r_1=2r_2$, which is a predictable result since the ramification index of $[E:F]$ is $2$. We will not use this fact in the sequel, though.

Now we prove proposition \ref{tvoth}. By lemma \ref{R1}, for every $D'$ and every $C'\in Ch_{D'}$, if $C$ is the only element of $Ch_D$ contained in the closure of $D\cup C'$, we have:
\[\phi_D(C')=(-q)^{-r_1d_\Gamma(D,D')}\phi_D(C),\]
hence:
\[\sum_{C'\in Ch_{D'}}\phi_D(C')=(-q)^{-r_1d_\Gamma(D,D')}\sum_{C\in Ch_{D}}\phi_D(C),\]
Let $d'$ be the dimension of $D$; we have the following lemma:
\begin{lemme}
Let $W'$ be the affine Weyl group of $G$ relative to $T_0$; the sum of the Poincar\'e series for a group of type $A_{d'}$ is:
\[\sum_{w\in W'}x^{l(w))}=\frac{1-x^{d'+1}}{(1-x)^{d'+1}}.\]
\end{lemme}

According to a formula given in the proof of \cite[corollary 3.4]{mcdo}, we have:
\[\sum_{w\in W'}x^{l(w))}=\prod_{i=1}^{d'}\frac{1-x^{m_i+1}}{(1-x)(1-x)^{m_i}},\]
where $m_1,\dots,m_{d'}$ are the exponents of $W'$ (see \cite[\S 6.2]{bou5}). On the other hand, according to \cite[plate I (X)]{bou}, we have $m_i=i$ for every $i$. The lemma follows then by an easy computation. $\Box$

We see immediately from this lemma that as soon as $|x|<1$, the sum in the left-hand side cannot be zero. Denote by $s(x)$ that sum.

By proposition \ref{elun}, the sum:
\[\sum_{C\in Ch_a^0}|\phi_D(C)|\]
converges, and we obtain, using lemmas \ref{R1} and \ref{R2} and taking into account the fact that $r_1$ and $r_2$ happen to be always even:
\[\sum_{C\in Ch_a^0}\phi_D(C)=\#(Ch_{D,a})s(q^{r_2-r_1}).\]
Since $r_2<r_1$,
the right-hand side is obviously nonzero. The proposition is now proved. $\Box$

\subsection{An Iwahori-spherical test vector}

In this last section, we prove that it is always possible to use a suitably chosen Iwahori-spherical vector as a test vector. We start by the following lemma:

\begin{lemme}\label{hctest}
For every codimension $1$ facet $D$ of $X_E$, we have:
\[\sum_{C\supset D}\phi_C=0.\]
\end{lemme}

Let $D$ be such a facet, and let $C'$ be any chamber of $X_E$: we have:
\[\sum_{C\supset D}\phi_C(C')=\sum_{C\supset D}(-q)^{-d(C,C')}.\]
Consider the closure $cl(D\cup C')$; by \cite[I, proposition 2.3.1]{bt}, it is contained in an apartment $A$ of $X_E$, and even in one of the two half-apartments of $A$ delimited by the wall containing $D$. Hence there exists exactly one chamber $C''$ containing $D$ and contained in $cl(D\cup C')$. Set $\delta=d(C'',C')$; if $C''''$ is another chamber of $X_E$ containing $D$, the closure of $C'''\cup C'$ must then contain $C''$, and since $C'''$ is neighboring $C''$, we must have $d(C''',C')=\delta+1$. Hence we have:
\[\sum_{C\supset D}\phi_C(C')=(-q)^{-\delta}+q((-q)^{-\delta-1})=0.\]
The lemma is then proved. $\Box$

Now we check that we can use some well-chosen Iwahori-spherical vector as a test vector when $G$ is not of type $A_{2n}$. In the case of type $A_{2n}$, we already know by proposition \ref{tva2n} that it is true.

\begin{prop}\label{testv}
Assume $G$ is not of type $A_{2n}$. Let $\lambda$ be any nonzero element of $\mathcal{H}(X_E)^{G_F,\chi}$, viewed as a linear form on $\mathcal{H}(X_E)^\infty$.
Let $C_0$ be any element of $Ch_a^0$ and let $\phi_{C_0}$ be the Iwahori-spherical vector associated to $C_0$. Then $\phi_{C_0}$ is a test vector for $\lambda$.
\end{prop}

We use the same argument as in \cite[proposition 6.2]{cou2}: since $St_E$ is an irreducible representation, it is generated by any of its nonzero vectors, for example an Iwahori-spherical vector $\phi$. We deduce from this that $\mathcal{H}(X_E)^\infty$ is generated as a ${\mth{C}}$-vector space by the $G_E$-conjugates of $\phi$, which are the Iwahori-spherical vectors $\phi_C$ attached to every chamber $C$ of $X_E$. 

By lemma \ref{hctest}, the Iwahori-spherical vectors satisfy relations between each other which are similar to the harmonicity condition.  Let $\lambda$ be a nonzero $(G_F,\chi)$-equivariant linear form on $\mathcal{H}(X_E)^\infty$. Assume $\lambda(f_{C_0})=0$. Then we prove in a similar way as for elements of $\mathcal{H}(X_E)^{G_{F,der}}$, using corollary \ref{suppaut} and propositions \ref{cclink}, \ref{chdalconj} and \ref{uniccha}, that $\lambda(\phi_C)=0$ for every $C\in Ch_E$ as well, which implies $\lambda=0$, and we thus reach a contradiction. Hence $\phi_{C_0}$ is a test vector for $\lambda$ and the corollary holds. $\Box$


\begin{thebibliography}{99}
\bibitem{ar} U.N. Anandavardhanan and C.S. Rajan, {\em Distinguished representations, base change and reducibility for unitary groups.}  Int. Math. Res. Not., vol. 14 (2005), pp. 183--192.
\bibitem{bot} A. Borel, J. Tits, {\em Groupes r\^eductifs.} Publications Math\^ematiques de l'IHES, vol. 27 (1965), pp. 55--151.
\bibitem{bou5} N. Bourbaki, {\em Groupes et alg\`ebres de Lie, chapitre 5: Groupes engendr\'es par des r\'eflexions}, Hermann.
\bibitem{bou} N. Bourbaki, {\em Groupes et alg\`ebres de Lie, chapitre 6: Syst\`emes de racines}, Hermann.
\bibitem{bc} P. Broussous, F. Court\`es, {\em Distinction of the Steinberg representation.} Int. Math. Res. Not. 2014, no. 11, pp. 3140--3157.  
\bibitem{br} K.S. Brown, \emph{Buildings}, Springer, 1996.
\bibitem{bt} F. Bruhat and J. Tits, {\em Groupes r\'eductifs sur un corps local. I. Donn\'ees radicielles valu\'ees.} Publ. Math. Inst. Hautes Etudes Sci. vol. 41 (1972), pp. 5--251.
\bibitem{car} R.W. Carter, {\em Finite groups of Lie type.} John Wiley \& Sons, 1985.
\bibitem{cou1} F; Court\`es, {\em Parametrization of tamely ramified maximal tori using bounded subgroups}: Annali di Matematica, vol. 188 (2009), pp. 1--33.
\bibitem{cou2} F. Court\`es, {\em Distinction of the Steinberg representation II: an equality of characters.}, Forum Mathematicum, vol. 27 (2015), pp. 3461-3475.
\bibitem{deb} S. DeBacker, {\em Parametrizing conjugacy classes of maximal unramified tori with Bruhat-Tits theory.}  Michigan Math. J., vol. 54 (2006), no. 1, pp. 157--178. 
\bibitem{ds} P. Delorme, V. S\'echerre, {\em An analogue of the Cartan decomposition for p-adic symmetric spaces of split p-adic reductive groups.} Pacific J. Math, vol. 251 (2011), pp. 1--21. 
\bibitem{dyn} E. B. Dynkin, {\em Semisimple subalgebras of semisimple Lie algebras}, American Mathematical Translations, vol. 6 (2), 1955, pp. 111--244.
\bibitem{mcdo} D. Macdonald, {\em The Poincar\'e series of a Coxeter group.} Mathematische Annalen, vol. 199 (1972), pp. 161--174
\bibitem{mat} N. Matringe, {\em Distinction of the Steinberg representation for inner forms of $GL_n$}, preprint, arXiv:1602.05101, 2016.
\bibitem{pr0} D. Prasad, {\em Invariant forms for representations of $GL_2$ over a local field.} Amer. J. Math, vol. 114 (1992), pp. 1317--1363.
\bibitem{pr} D. Prasad, {\em On a conjecture of Jacquet about distinguished representations of $GL_n$.} Duke Math J., vol. 109 (2001), pp. 67--78.
\bibitem{pr2} D. Prasad, {\em A "relative" local Langlands correspondence}, preprint, arXiv:1512.04347, 2016.
\bibitem{sha}J.A. Shalika, {\em On the space of cusp forms of a p-adic Chevalley group.} Ann. of Math., vol. 92 (1970), no 2, pp. 262--278
\bibitem{spr} T.A. Springer, {\em Linear algebraic groups.} Progress in Mathematics, Birkh\"auser, vol. 9, 2001 (2nd edition).
\bibitem{tits} J. Tits. {\em Classification of algebraic semisimple groups.} Proc. Sympos. Pure Math., vol. 26 (1965), pp. 33--62.

\end{thebibliography}
\end{document}